\documentclass[a4paper,10pt]{article}
\usepackage{amssymb}
\usepackage{latexsym}
\usepackage{amsmath}
\usepackage{amsthm}
\usepackage{amsfonts}
\usepackage{amssymb}
\usepackage{graphicx}

\newtheorem{Th}{Theorem}

\newtheorem{Co}{Corollary}
\newtheorem{Lm}{Lemma}
\newtheorem{Lma}{Lemma}[section]
\newtheorem{Dfi}{Definition}
\newtheorem{Rm}{Remark}

\newcommand{\be}{\begin{equation}}
\newcommand{\ee}{\end{equation}}
\newcommand{\bes}{\begin{equation*}}
\newcommand{\ees}{\end{equation*}}

\newcommand{\R}{\mathbb{R}}
\newcommand{\N}{\mathbb{N}}

\newcommand{\Z}{\mathbb{Z}}

\newcommand\res{\mathop{\hbox{\vrule height 7pt width .5pt depth 0pt
\vrule height .5pt width 6pt depth 0pt}}\nolimits}

\def\theequation{\thesection.\arabic{equation}}
\def\theTh{\Roman{section}.\arabic{Th}}
\def\theProp{\Roman{section}.\arabic{Prop}}
\def\theCo{\Roman{section}.\arabic{Co}}

\def\theRm{\Roman{section}.\arabic{Rm}}
\newcommand{\reset}{\setcounter{equation}{0}\setcounter{Th}{0}\setcounter{Prop}{0}\setcounter{Co}{0}
\setcounter{Lm}{0}\setcounter{Rm}{0}}

\def\ti{\tilde}
\def\lf{\left}
\def\rg{\right}

\def\al{\alpha}
\def\la{\lambda}

\def\ep{\varepsilon}

\def\ds{\displaystyle}
\def\ov{\overline}
\def\Om{\Omega}
\def\om{\omega}
\def\p{\partial}

\def\vP{\vec{\Phi}}

\def\gP{g_{\vec{\Phi}}}

\def\res{\mathop{\hbox{\vrule height 7pt width .5pt 
depth 0pt\vrule height .5pt width 6pt depth 0pt}}\nolimits}

\begin{document}

\title{ Willmore Minmax Surfaces and the Cost of the Sphere Eversion}

\author{ Tristan Rivi\`ere\footnote{Department of Mathematics, ETH Zentrum,
CH-8093 Z\"urich, Switzerland.}}

\maketitle

\medskip

{\bf Abstract :} {\it We develop a general Minmax procedure in Euclidian spaces for constructing Willmore surfaces of non zero indices. We implement
this procedure to the Willmore Minmax Sphere Eversion in the 3 dimensional euclidian space. We compute
the cost of the Sphere eversion in terms of Willmore energies of Willmore Spheres in ${\R}^3$.}
\medskip

\noindent{\bf Math. Class. 49Q10, 53A30, 53A05, 58E15, 58E30, 35J35, 35J48}

\section{Introduction}
\subsection{The search for Willmore Minmax Surfaces.}
Finding optimal shapes is a search probably as old as mathematics and whose motivation goes beyond the exclusive quest for beauty.
It is often closely related to the understanding of deep mathematical structures and ultimately to natural phenomena happening to be governed by these pure ideas.

The existence of closed geodesics on arbitrary manifolds as well as its higher dimensional counterpart, the existence question  of  minimal surfaces, belongs to this search and has been since
the XIXth century a very active area of research diffusing in other areas of mathematics and science in general much beyond the field of Differential Geometry.

The theory of Willmore surfaces, introduced by Wilhelm Blaschke around 1920, grew up out of the attempt to merge  {\it minimal surface theory} and {\it conformal invariance}.

 For an arbitrary immersion $\vec{\Phi}$ of a given oriented abstract
surface $\Sigma$ into an euclidian space ${\R}^m$ Wilhelm Blaschke introduced the lagrangian
\[
W(\vec{\Phi})=\int_{\Sigma}|\vec{H}_{\vec{\Phi}}|^2\ dvol_{g_{\vec{\Phi}}}
\]
where $\vec{H}_{\vec{\Phi}}$ and $dvol_{g_{\vec{\Phi}}}$ are respectively the mean curvature vector and the volume form of the metric
induced by the immersion. He proved in particular that for a closed surface $\Sigma$ this lagrangian is  invariant under conformal transformations.

The critical points to $W$ are called  ''Willmore Surfaces\footnote{This denomination has spread, and is now generally used, after the work \cite{Wi} which relaunched the  
study of these surfaces that Blaschke originally named ''conformal-minimal surfaces''.}''
The known set of critical points to the Willmore Lagrangian has been for a long time reduced to the minimal surfaces
and their conformal transformations. This maybe explains why the study of it's variations has been more or less stopped  during several decades following the seminal work of Wilhelm Blaschke which was slowly sinking into oblivion.

 After the reviving work of Tom Willmore the first main contribution to  ''Willmore Surfaces'' has been brought by Robert Bryant in \cite{Br}. Using algebraic geometric techniques he succeeded in describing all the immersed ''Willmore Spheres'' in ${\R}^3$ as being given exclusively by the images by  inversions of simply connected complete non compact minimal surfaces with planar ends. The Willmore energy of the immersed Willmore spheres was consequently proved to be equal to $4\pi$ times the number of planar ends.
Due to the non-triviality of the space of holomorphic quartic forms on any other riemann surface this approach has been restricted to the sphere exclusively.
Other algebraic geometric approaches for studying critical points to the Willmore lagrangian include ``spectral curve methods'' and integrable system theory, but these rather abstract methods are addressing issues which are mostly local and, until now, hardly translatable into ``back to earth'' results exhibiting new complete Willmore surfaces or characterizing 
the space of Willmore critical points in a decisive way.

Beside algebraic geometric methods a natural strategy for producing new Willmore surfaces would consist in developing the fundamental principles of the calculus of variations applied to the Willmore Lagrangian. Since $\vec{H}_{\vec{\Phi}}=2^{-1}\, \Delta_{g_{\vec{\Phi}}}\vec{\Phi}$ the Willmore Lagrangian
is nothing but $1/4$-th of the $L^2$ norm of the Laplacian of the immersion and is showing in that sense it's 4th order elliptic nature. This coercive structure gives some hope for the success of the variational methods. The pioneered work
studying the variations of $W$ has been written by Leon Simon (see \cite{Sim}) in which he was proving the existence of a torus minimizing the Willmore energy. This existence result was also motivated by the conjecture formulated by Willmore in \cite{Wi} according to which the torus obtained by rotating around the vertical $z$
axis of the vertical circle of radius 1 centered at $(\sqrt{2},0,0)$ and included in the plane $y=0$ would be the unique minimizer modulo the conformal transformations. This conjecture has been finally proved some years ago by Fernando Cod\'a Marques and Andre Neves in \cite{MN1}.
In \cite{BK} Matthias Bauer and Ernst Kuwert succeeded in proving a succession of strict inequalities excluding possible degeneracies and the splitting of the
underlying surface which was still left open in Leon Simon's argument for arbitrary genus. As a consequence the authors proved the existence of
a minimizer of $W$ for any arbitrary closed orientable two dimensional manifold $\Sigma$.

Leon Simon's approach to the minimization of Willmore is based on energy comparison arguments and local bi-harmonic graph approximation procedures
and in that sense is shaped for studying the ground states of index 0. This approach is mostly considering the image of the immersion $\vec{\Phi}(\Sigma)$ and not the immersion {\it per se} and can be called ''ambiant'' for that reason. In \cite{Ri3} the author of the present work has introduced an alternative proof to Leon Simon's existence result  using an approach called ``parametric''. In this approach the study of the variations of the immersion is made possible by local extraction of ``Coulomb Gauges'' (isothermic parametrization) and the use of the conservation laws
issued from the application of Noether theorem combined with the {\it integrability by compensation theory} (see also a systematic presentation of this theory
and it's application to Willmore in \cite{Ri1}). Since this approach is not making use of comparison arguments and since it is based on a weak formulation
of the Willmore Euler Lagrange equation discovered in \cite{Ri2}, it was giving the author good hope to apply it to more diversified calculus of variation arguments than strict minimizations. This is the main achievement of the present work.
Precisely, the purpose of the paper is to present a Minmax method for producing critical points to the Willmore energy of non zero indices.

\subsection{''Smoothers'' based on the second fundamental form.}

As already mentioned the Willmore energy is invariant under the action of the M\"obius group of conformal transformations of ${\R}^m$ which is known
to be non compact.  For that reason in particular it does not satisfy the Palais Smale condition. This is an obstruction for applying Minmax variational
principles such as the {\it mountain pass lemma} directly. We shall then adopt a {\it viscosity approach} and add to the Willmore energy what we call a {``smoother''} times a small ``viscosity parameter'' $\sigma^2$
\[
\mbox{Full Energy}(\vec{\Phi}):=W(\vec{\Phi})+\sigma^2\,\mbox{Smoother}(\vec{\Phi})
\] that makes
for the new energy the Palais-Smale condition being satisfied.  One can then apply the mountain path arguments to such energy and produce 
Minmax critical points. In the second part of the procedure one  makes $\sigma$ tend to zero
and one studies the converging process hopefully to a Willmore Minmax Surface.

In a first approach, following \cite{Ri4}, one could think adding to $W$ a term of the form
\[
\mbox{Smoother}(\vec{\Phi})=\int_\Sigma(1+|\vec{\mathbb I}_{\vec{\Phi}}|^2)^2
\] 
where $\vec{\mathbb I}_{\vec{\Phi}}$ is the second fundamental form of the immersion $\vec{\Phi}$. This will make the new Willmore relaxed energy satisfy the {\it Palais-Smale condition} (as proved in \cite{KLL}, see also \cite{La}) but this will bring us to the study
of $p-$harmonic systems which makes the analysis of the convergence rather involved - in particular the energy quantization - when the small viscosity parameter $\sigma$ tends to zero. From that perspective $p-$harmonic versus harmonic systems, as observed below and as it is also used in \cite{BR4},
replacing the full second fundamental form by it's trace $\vec{H}_{\vec{\Phi}}$ and consider instead
\[
\mbox{Smoother}(\vec{\Phi})=\int_\Sigma(1+|\vec{H}_{\vec{\Phi}}|^2)^2
\] 
has the surprising effect to make the highest order term in the Euler Lagrange to be $\Delta$ and not $\Delta_2:=div( (1+|H|^2) \nabla)$ if one makes use
of the various conservation laws issued from Noether theorem, following the main lines of \cite{Ri2}. The drawback however is that $\int_\Sigma(1+|\vec{H}_{\vec{\Phi}}|^2)^2$ fails to satisfy the {\it Palais-Smale condition} and cannot be a {\it smoother} by itself and has to be ``reinforced''.

\subsection{Frame Energies}

In the portofolio of surface energies, the author, in collaboration with Andrea Mondino, introduced in \cite{MR3} the notion of {\it frame energy} for an arbitrary immersion of a torus $\vec{\Phi}$ : $T^2\ \longrightarrow\ {\R}^m$ equipped by an orthonormal tangent frame $\vec{e}\ :\ T^2\,\rightarrow  S^{m-1}\times S^{m-1}$  where $\vec{e}\,(x)$ is realizing 
an orthonormal basis of $\vec{\Phi}_\ast T_xT^2$. The frame energy is then simply given by
\be
\label{I.1-1}
{\mathcal F}(\vec{\Phi},\vec{e}):=\frac{1}{4}\int_{T^2}|d\vec{e}|^2_{g_{\vec{\Phi}}}\ dvol_{g_{\vec{\Phi}}}\ge W(\vec{\Phi})
\ee
If one considers the ``Coulomb frame'' associated to the conformal immersion of a fixed flat torus\footnote{By the uniformization theorem such a parametrization always exists}, the {\it frame energy} $F(\vec{\Phi})$ is then defining an energy of the 
immersion $\vec{\Phi}$  that hapens to be more coercive than the Willmore energy itself. One could ask whether it can be naturally extended to any other immersion of an arbitrary surface $\Sigma$. This happens to be indeed the case as we explain in the sequel.

Let $\vec{\Phi}$ be an arbitrary immersion of a closed surface $\Sigma$ and denote by $g_0$ a constant scalar curvature metric of volume one on $\Sigma$ for which there exists $\al\ : \Sigma\ \rightarrow \ {\R}$
\be
\label{I.1}
\gP=e^{2\al}\ g_0\quad.
\ee
For $\Sigma\ne S^2$ the function $\al$ is defined without ambiguity, whereas in the sphere case we have to count with the action of a ''gauge group'',
the space ${\mathcal M}^+(S^2)$ of positive conformal transformations of $S^2$ and $\al$ is uniquely defined modulo the action of this gauge group. For the case of the torus one proves in this paper that
\[
F(\vec{\Phi})=W(\vec{\Phi})+\frac{1}{2}\int_{T^2}|d\al|^2_{\gP}\ dvol_{\gP}
\] 
We generalize the frame energy for any surface of genus larger than 1 as follows : $F:= W+{\mathcal O}$ where
\be
\label{I.0}
 {\mathcal O}(\vec{\Phi})= \frac{1}{2}\int_{\Sigma}|d\al|^2_{\gP}\ dvol_{\gP}+K_{g_0}\,\int_{\Sigma}\al\ e^{-2\al}\ dvol_{\gP}-2^{-1}\, K_{g_0}\,\log\int_\Sigma dvol_{\gP}
\ee
where $K_{g_0}$ is the constant scalar curvature metric of $g_0$. The reason why adding to the Dirichlet energy of $\al$ the term $$K_{g_0}\,\int_{\Sigma}\al\ e^{-2\al}\ dvol_{\gP}$$
comes from the fact that the first variation can be expressed miraculously by \underbar{local} quantities ! (although the operation which to $g_{\Phi}$ assigns $\al$ is highly
non-local - see for instance \cite{Tro}). Finally the reason why adding the third term is twofold :  it makes the energy scaling invariant and non negative as a direct application of Jensen when $K_{g_0}<0$.

Finally, when $\Sigma=S^2$ we define identically the {\it Frame energy} to be
\[
F(\vec{\Phi}):=W(\vec{\Phi})+\frac{1}{2}\int_{S^2}|d\al|^2_{\gP}\ dvol_{\gP}+ 4\,\pi\ \int_{S^2}\al\ e^{-2\al}\ dvol_{\gP}-\, 2\,\pi\ \log \int_{S^2}dvol_{\gP}
\]
Beside the fact that it is naturally generalizing to $S^2$ the Dirichlet energy of {\it Coulomb frames} on tori, there are 4 main reasons why considering this special expression with these particular coefficients
\begin{itemize}
\item[i)] $F-W$ is the well known {\it Onofri Energy} of $\al$ (see \cite{On} and a more recent presentation in \cite{DEJ}) and satisfies\footnote{This inequality is not the direct consequence of Jensen inequality when $K_{g_0}>0$ and requires more elaborated arguments.}
\be
\label{I.2}
{\mathcal O}(\vec{\Phi}):=\frac{1}{2}\int_{S^2}|d\al|^2_{\gP}\ dvol_{\gP}+ 4\,\pi\ \int_{S^2}\al\ e^{-2\al}\ dvol_{\gP}-\, 2\,\pi\ \log \int_{S^2}dvol_{\gP}\ge 0
\ee
\item[ii)] The first variation of $F$ is explicit and can be expressed using \underbar{local} quantities.
\item[iii)] The energy ${\mathcal O}(\vec{\Phi})$ is gauge invariant with respect to the action of the M\"obius group on $S^2$ and is independent of the choice
of $\al$ and $g_0$ satisfying (\ref{I.1}) and depends only on $\vec{\Phi}$
\item[iv)] The $F$-energy is dilation invariant : $F( e^t\ \vec{\Phi})=F(\vec{\Phi})$ for any $t\in{\R}$.
\end{itemize}


\medskip

\noindent{\bf Open Problem 1.} {\it It would be interesting to study the existence of minimizers of the frame energy $F$ in each regular homotopy class of immersions of spheres in ${\R}^4$. Since the work of Stephen Smale \cite{Sm} we know that there exists countably many of them given by
the $\pi_2(V_{4,2}({\R}))={\Z}$, the second homotopy class of the Stiefel manifold of 2-frames in ${\R}^4$. It would also be 
interesting to study the asymptotic dependence with respect to the class of the infimum of the $F-$energy as the class goes to infinity (if it is linear or sub-linear). }

\medskip 

Finally, for a given immersion $\vec{\Phi}$ it would be interesting to study the minimal Dirichlet energy  of any bundle map from $T\Sigma$ into ${\R}^3\times{\R}^3$  which is an isometry from each fiber $(T_x\Sigma,g_{\vec{\Phi}})$ into $\vec{\Phi}_\ast T_{x}\Sigma\subset G_2({\R}^3)$ and which projects onto $\vec{\Phi}$ - the map $(\vec{\Phi},d\vec{\Phi})$ is one of such maps of course. Starting from  $(\vec{\Phi},d\vec{\Phi})$ such a bundle map is just given by the choice of an $S^1$ rotation at each point.


In the case of $\Sigma=T^2$ this coincides with the Dirichlet energy of an optimal global frame and is equal to $F(\vec{\Phi})$. 
\subsection{A Viscous Approximation of Willmore.}

Inspired by the discussion above we propose to consider the following approximation of the Willmore energy
\[
F^\sigma(\vec{\Phi}):=W(\vec{\Phi})+\sigma^2\,\int_{\Sigma}(1+|\vec{H}|^2)^2\ dvol_{\vP}+\frac{1}{\log\frac{1}{\sigma}}{\mathcal O}(\vec{\Phi})
\]
where ${\mathcal O}(\vec{\Phi})$ is given by (\ref{I.2}) for $\Sigma= S^2$ or  by the expression (\ref{I.0}) otherwise.

We prove in section {V} that $F^\sigma$ satisfies the Palais Smale condition. One can then apply Minmax arguments to $F^\sigma$
for any admissible family.
One of the main achievement of the present work is the proof an $\epsilon-$regularity \underbar{independent of $\sigma$} (see lemma~\ref{lm-ep-reg}). It  is making use of the special choice we made of the logarithmic dependence  of the small parameter with respect to the viscosity parameter $\sigma$ in front of ${\mathcal O}$ . This $\epsilon-$regularity permits to pass to the limit in the equation for well chosen sequences of Minmax critical points  of $F^\sigma$ as $\sigma$ goes to zero.
The last main lemma in the paper is an energy quantization result when $\Sigma=S^2$ (see lemma~\ref{lm-energy-quantiz}). It roughly says that no energy can be dissipated in neck regions.

We then have the main tools for performing Minmax procedures for the Willmore energy of spheres. To that aim we introduce the space of $W^{2,4}(S^2,{\R}^m)$ immersion of $S^2$ into ${\R}^m$ that we denote ${\mathcal E}_{\Sigma,2}({\R}^m)$. This space is equipped with the $W^{2,4}$ topology. Our main result is the following.

\begin{Th}
\label{th-minmax}{\bf [Willmore Minmax Procedures for Spheres]}
Let $m\ge 3$ and $k\ge 1$. Let ${\mathcal A}$ be an admissible family of $W^{2,4}$ immersions of the sphere $S^2$ into ${\R}^m$, that is a family of images by continuous maps from a fixed topological space into the space ${\mathcal E}_{\Sigma,2}({\R}^m)$ which is invariant under the action of homeomorphism
of ${\mathcal E}_{S^2,2}({\R}^m)$ isotopic to the identity. Assume
\[
\beta_0:=\inf_{A\in {\mathcal A}}\max_{\vec{\Phi}\in A}W(\vec{\Phi})>0
\]
then there exists finitely many  Willmore immersions of $S^2$ minus finitely many points, $\vec{\xi}_1\cdots\vec{\xi}_n$,  such that
\be
\label{bilan}
\beta_0=\sum_{i=1}^n W(\vec{\xi}_i)-\,4\pi\, N
\ee
where $N\in {\N}$.\hfill $\Box$
\end{Th}
\begin{Rm}
\label{rm-I.1}
The theorem~\ref{th-minmax} goes with a ``bubble tree convergence'' of a sequence of Minmax critical points of $F^{\sigma^k}$ for some well chosen sequence $\sigma^k$ tending to zero. This convergence produces asymptotically a ``bubble tree'' of Willmore spheres some of them being shrunk to zero.
Among the ones which shrunk to zero there might be non compact (after asymptotic rescaling) simply connected Willmore surfaces with ends at infinity that we have to
inverse in order to make them being Willmore sphere. This operation in producing energy given by an integer multiple times $4\pi$. This is why such
a quantity is subtracted in (\ref{bilan}).
One of the hard parts in the whole proof is to show that between two successive asymptotic Willmore Spheres in that ``bubble tree'' no energy is lost at the limit.
This is the so called ``no neck energy'' property.   \hfill $\Box$
\end{Rm}
\noindent{\bf Open Problem 2.}
{\it It is left open at this stage to know whether among the maps $\vec{\xi}_i$ the ones  obtained after inversions in the bubble tree define or not ``true'' Willmore possibly branched immersions. By ``true'' we mean that the first residue\footnote{In 3 dimension this residue is also a multiple of the one that can be deduced from the integration of the one form (4.5) in \cite{KS}.} defined in 
\cite{Ri2}
\[
\int_{\Gamma} \p_\nu\vec{H}-3\, \pi_{\vec{n}}\lf( \p_\nu\vec{H}\rg)+\star\, \p_\tau\vec{n}\wedge\vec{H}
\]
is zero for any closed curve $\Gamma$ avoiding the center of the inversion. The inversion of the Catenoid is not a ``true'' Willmore surface in that sense whereas the inversion of the Enneper surface is a ``true'' Willmore Sphere with a branch point of multiplicity 3 at the origin (see \cite{BR2}).}

\medskip

\noindent{\bf Open Problem 3.} {\it Extend the previous result to general surfaces. The ``only'' obstruction comes from the Energy Quantization
result which is missing when the conformal class of the Minmax sequence of $F^{\sigma^k}$ possibly degenerates. The recent progresses made in
\cite{LR3} and \cite{LR4} should be of great help in solving this difficulty.\hfill $\Box$}

One consequence of the previous result is the following corollary. One considers the family ${\mathcal A}$ of loops into Imm$(S^2,{\R}^3)\simeq_{{hom}} SO(3)\times \Om^2(SO(3))$ realizing a non trivial element of
$\pi_1(\mbox{Imm}(S^2,{\R}^3))\simeq {\Z}_2\times{\Z}$. It is proved in \cite{BM} that for instance  the Froisart-Morin sphere eversion followed by the mirror image of the time reversed  of the same eversion is generating $\pi_1(\mbox{Imm}(S^2,{\R}^3))$. In order to avoid uninteresting loops coming from the action
of $\mbox{Diff}(S^2)$ one should rather work modulo the action of reparametrization of the sphere and consider the infinite orbifold\footnote{While $\mbox{Imm}(S^2,{\R}^3)/ \mbox{Diff}^+(S^2)$ the quotient
of $\mbox{Imm}(S^2,{\R}^3))$  by the group of \underbar{positive} diffeomorphisms of $S^2$ has a nice bundle structure due to the free action of $\mbox{Diff}^+(S^2)$
on $\mbox{Imm}(S^2,{\R}^3))$, this is not the case anymore for $\mbox{Imm}(S^2,{\R}^3)/ \mbox{Diff}(S^2)$. The space  $\mbox{Imm}(S^2,{\R}^3)/ \mbox{Diff}(S^2)$ is an infinite orbifold obtained by the quotient of $\mbox{Imm}(S^2,{\R}^3)/ \mbox{Diff}^+(S^2)$ by the map $x\rightarrow -x$.   This induces a 2 sheets covering away from the subspace of singular orbits which happens to have infinite codimension (see section 3 of \cite{CMM}). Because of the smallness of the size of singular orbits, using transversality arguments, one can  compute homotopy groups of $\mbox{Imm}(S^2,{\R}^3)/ \mbox{Diff}(S^2)$ as if this covering map would be without singularities. }
 $\mbox{Imm}(S^2,{\R}^3)/ \mbox{Diff}(S^2)$ instead of $\mbox{Imm}(S^2,{\R}^3))$ which is an open subspace of the Banach space $W^{2,4}(S^2,{\R}^3)$. 

 One can then take ${\mathcal A}$ to be 
the canonical projections onto $\mbox{Imm}(S^2,{\R}^3)/ \mbox{Diff}(S^2)$ of
paths from $[0,1]$ into $\mbox{Imm}(S^2,{\R}^3)$
  homotopic
 to a non trivial element in  $\pi_1(\mbox{Imm}(S^2,{\R}^3)/ \mbox{Diff}(S^2))={\Z}$. The projection of one time the Froisart-Morin sphere eversion gives such a loop for instance\footnote{The canonical projection $\pi$ of $\mbox{Imm}(S^2,{\R}^3)$ onto the infinite orbifold $\mbox{Imm}(S^2,{\R}^3)/ \mbox{Diff}(S^2)$ induces
  a morphism between the two first homotopy groups $\pi_\ast$ : $\pi_1(\mbox{Imm}(S^2,{\R}^3))=  {\Z}_2\times{\Z}\rightarrow \pi_1(\mbox{Imm}(S^2,{\R}^3)/ \mbox{Diff}(S^2))={\Z}$ equal to the multiplication by 2}.
\begin{Co}
\label{co-I.1} {\bf [The Cost of the Sphere Eversion]}
Let ${\Omega}$ be the space of continuous paths of $C^2$ immersions into ${\R}^3$ joining the standard sphere $S^2$ with the two opposite orientations and homotopic to the Froisart-Morin Sphere  eversion. Define $\beta_0$ to be the ``cost of the sphere eversion''  by
\[
\beta_0:=\inf_{\om\in {\Om}}\max_{\vec{\Phi}\in \om}W(\vec{\Phi})
\]
then there exists finitely many  Willmore immersions of $S^2$ minus finitely many points, $\vec{\xi}_1\cdots\vec{\xi}_n$,  such that
\[
\beta_0=\sum_{i=1}^n W(\vec{\xi}_i)-\,4\pi\, N
\]
where $N\in {\N}$.
\hfill $\Box$
\end{Co}
\begin{Rm}
\label{rm-I.7}
Performing the Willmore Minmax Sphere Eversion has been originally proposed by Robert Kusner.
\hfill $\Box$
\end{Rm}
\noindent{\bf Open Problem 4. } {\it Is it true that in corollary~\ref{co-I.1}
\[
\beta_0=16\pi\quad ?
\]
A topological result, see \cite{BM} (see also the enlightening proof in \cite{Hug}), asserts that any element in $\Omega$ has to contain at least one immersion with a point of self intersection of order 4 (i.e. a quadruple point). Hence using 
Li-Yau's result we deduce that $\beta_0\ge 16\pi$. In \cite{FSKBHC} a candidate for the realization of $\beta_0$ is proposed. It is the inversion of a simply connected complete minimal surface with 4 planar ends. Hence the Willmore energy of this candidate is $16\pi$. Interesting computations reinforcing this conjecture are performed in this work. Establishing upper-bound of the lowest energy Minmax sphere eversion is making difficult by the fact that
producing concrete sphere eversions is highly challenging and has been at the origin of many rigorous works, computer simulations and videos too starting from the first example given by Arnold Shapiro (see for instance \cite{MP}, \cite{FM}, \cite{Phi}, \cite{Ape}, \cite{FSH}, \cite{Sul}, \cite{Ait}...)  }
\begin{Rm}
\label{rm-I.3}
If one could prove that the $\vec{\xi}_i$ are true Willmore spheres in the sense clarified above and if one could extend Bryant's result to general ``true'', possibly branched Willmore sphere, then one would know that the cost of the sphere eversion is $4\pi$ times an integer.\hfill $\Box$
\end{Rm}
 \begin{Rm}
 \label{rm-I.4}
 An interesting  upper-bound of the  cost of the Total Curvature Minmax Sphere Eversion  is presented in \cite{Ek}. In this paper it is proved in particular that
\[
\inf_{\om\in {\Om}}\max_{\vec{\Phi}\in \om}\int_{S^2}|K_{\vec{\Phi}}|\ dvol_{g_{\vec{\Phi}}}\le 8\pi
\]
where $K_{\vec{\Phi}}$ is the Gauss curvature of $\vec{\Phi}$.
\hfill $\Box$
\end{Rm}

\medskip

\noindent{\bf Open Problem 5. }
 {\it It would be interesting to study the cost of the Frame Energy $W+{\mathcal O}$ Minmax Sphere Eversion. As we saw above this energy is closely related to the minimizing Dirichlet energy among the  bundle maps injections induced by $(\vec{\Phi},d\vec{\Phi})$ - which are used by Smale to compute the homotopy type
 of the space of immersions.}

 \medskip
 
 \noindent{\bf Open Problem 6. }{\it It would be interesting, beside the study of the Wilmore Minmax Sphere Eversion exclusively, to explore also the Willmore Minmax
 for ${\mathcal A}$ corresponding to other non trivial classes in $\pi_k( \mbox{Imm}(S^2,{\R}^3))=\pi_k( SO(3)\times \Om^2(SO(3)))$ for arbitrary $k$. How would their indices be related to $k$ ?}
 
 \medskip
 
 \noindent{\bf Open Problem 7. }{\it Explore the topology of $\mbox{Imm}(S^2,{\R}^3)$ using $W$ as some kind of ``quasi\footnote{``quasi'' because we know that $\p^2W$ has at least the null directions given by the action of M\"obius transformations.} Morse function'' : For instance we could ask if  all the Willmore immersions of $S^2$ into ${\R}^3$ - described by Robert Bryant - are related to a Minmax procedure
 involving the various classes of the various groups $\pi_k( \mbox{Imm}(S^2,{\R}^3))$ ? In order to start ``slow'' this far-reaching question a first certainly instructive step would consist in computing the indices of the Willmore Sphere Immersions in ${\R}^3$. 
 Of course one would have also to complete the space of immersions by considering possibly branched immersions and, up to now, there is no general result known about the extension of Bryant classification of Willmore Immersed Spheres in ${\R}^3$  to  Willmore Branched Spheres beside some very special cases treated in \cite{LN}.}
 
 \medskip
 
Most of the proofs below are presented in the particular case $m=3$ in order to make the presentation of them more accessible. The general case $m\ge 3$
is very similar but requires the use of the conservation laws in arbitrary codimensions introduced in \cite{Ri2} whose formulation involves the use of multi-vectors instead of vectors and are a bit more tedious but do not bring any new fundamental difficulties.

\section{The space of immersions into ${\R}^3$ with $L^q$ bounded second fundamental form.}
\reset

 For $k\in {\N}$ and $1\le q\le +\infty$ We recall the definition of $W^{k,q}$ Sobolev function on a closed smooth surface $\Sigma$ (i.e. $\Sigma$ is compact without boundary). To that aim we take
some reference smooth metric $g_0$ on $\Sigma$
\[
W^{k,q}(\Sigma,{\R}):=\lf\{f\ \mbox{ measurable } \; s.t \quad \nabla_{g_0}^kf\in L^q(\Sigma,g_0)\rg\}
\]
where $\nabla^k_{g_0}$ denotes the $k-$th iteration of the Levi-Civita connection associated to $\Sigma$. Since the surface is closed the space defined in this way is independent of $g_0$.

For $p\ge1$,  following \cite{Ri3} in which the case $p=1$ was considered, we define the space ${\mathcal E}_{\Sigma,p}$ of {\it weak immersions of $\Sigma$ into ${\R}^3$ with $L^{2p}$ bounded second fundamental form as follows}.
\[
{\mathcal E}_{\Sigma,p}:=\lf\{
\begin{array}{l}
\vec{\Phi}\in W^{1,\infty}(\Sigma, {\R}^3) \quad\mbox{ s. t. }\exists \ C>1\\[3mm]
C^{-1}\ g_0\le \vec{\Phi}^\ast g_{{\R}^3}\le C\ g_0\\[3mm]
\mbox{ the Gauss map }\vec{n}_{\vec{\Phi}}\in W^{1,2p}(\Sigma, Gr_2({\R}^{3}))
\end{array}
\rg\}
\]
For any $\vec{\Phi}\in{\mathcal E}_{\Sigma,p}$, starting from the equation $\Delta_{g_{\vec{\Phi}}}\vec{\Phi}=2\,\vec{H}$, classical elliptic estimates permit to bootstrap  (in the case $p>1$)
 and get that $\vec{\Phi}$ is in fact in $W^{2,2p}(\Sigma,{\R}^3)$ (see \cite{Ri1}). It is then not difficult to observe that ${\mathcal E}_{\Sigma,p}$ is in fact, \underbar{for $p>1$} an open
subset of the Banach space $W^{2,2p}(\Sigma,{\R}^3)$. The border line case $p=1$ is more delicate, it was introduced first in \cite{Ri3}  as being of prior interest for studying
the variations of Willmore and is extensively presented also in \cite{Ri1}.

\section{Frame Energies.} 
\reset
For any weak immersion $\vec{\Phi}$ in ${\mathcal E}_{\Sigma,p}$ 
we denote $A_{\gP}$  the connection 1-form - which lives on $\Sigma$ since the bundle is abelian - equal to the difference between the Levi-Civita connection $\nabla^{\gP}$ defined by $\gP$ and the Levi-Civita connection $\nabla^{g_0}$
\[
A_{\gP}:=\nabla^{\gP}-\nabla^{g_0}
\]
where $g_0$ is a constant scalar curvature metric of volume one on $\Sigma$ for which there exists $\al\ : \Sigma\ \rightarrow \ {\R}$
\[
\gP=e^{2\al}\ g_0\quad.
\]
For $\Sigma\ne S^2$ the function $\al$ is defined without ambiguity, whereas in the sphere case we have to count with the action of a ''gauge group'',
the space ${\mathcal M}^+(S^2)$ of positive conformal transformations of $S^2$ and $\al$ is uniquely defined modulo the action of this gauge group. We shall now express $|A_{\gP}|_{g_{\vec{\Phi}}}^2$ locally using moving frame. Let $(e_1,e_2)$ be an orthonormal local frame for the metric $\gP$. We have
\[
|\nabla^{\gP}-\nabla^{g_0}|^2_{\gP}=\sum_{i,j=1}^2\lf|\lf(\lf(\nabla^{\gP}-\nabla^{g_0}\rg)e_i,e_j\rg)\rg|^2_{\gP}
\]
Since $e_i$ is a unit vector field for $\gP$ we have $(\nabla^{\gP}e_i,e_i)=0$ and since $f_i:=e^\al\,e_i$ is a unit vector field for $g^0$ we have 
$$
(\nabla^{g_0}e_i,e_i)=e^{-2\al}\, (\nabla^{g_0}f_i,f_i)_{\gP}-\, e^{-2\al}\ de^\la\ (e_i,f_i)_{\gP}=-\ d\al
$$
Hence
\[
\sum_{i=1}^2\lf|\lf(\lf(\nabla^{\gP}-\nabla^{g_0}\rg)e_i,e_i\rg)\rg|^2_{\gP}=\,2\ |d\al|^2_{\gP}
\]
In order to compute $\lf|\lf(\lf(\nabla^{\gP}-\nabla^{g_0}\rg)e_1,e_2\rg)\rg|^2_{\gP}$ we choose local conformal coordinates $(x_1,x_2)$ for $g_{\vec{\Phi}}$ and we have then the existence of $\mu$ locally such that respectively
$\gP=e^{2\, (\al+\mu)}\ [dx_1^2+dx_2^2]$ and $g_0=e^{2\mu}\ [dx_1^2+dx_2^2]$. We choose $e_i:=e^{-\al-\mu}\ \p_{x_i}$ and thus $f_i= e^{-\mu}\ \p_{x_i}$. Using classical computation of the Levi-Civita connection of a metric
in conformal charts (see \cite{Ri1}) we have 
\[
\lf(\nabla^{\gP}e_1,e_2\rg)_{\gP}= \ast\ d(\al+\mu)   \quad\mbox{ and }\quad \lf(\nabla^{g_0}f_1,f_2\rg)_{g_0}=\ast\ d\mu
\]
Since $e_1$ and $e_2$ are orthogonal to each other with respect to $g_0$, we have
\[
\lf(\nabla^{g_0}e_1,e_2\rg)_{\gP}= e^{2\al}\, \lf(\nabla^{g_0}e_1,e_2\rg)_{g_0}=\lf(\nabla^{g_0}f_1,f_2\rg)_{g_0}=\ast\ d\mu
\]
Hence
\[
\lf(\lf(\nabla^{\gP}-\nabla^{g_0}\rg)e_1,e_2\rg)_{\gP}=\ast d\al\quad.
\]
Combining the previous we obtain
\be
\label{V.1}
|\nabla^{\gP}-\nabla^{g_0}|^2_{\gP}=\sum_{i,j=1}^2\lf|\lf(\lf(\nabla^{\gP}-\nabla^{g_0}\rg)e_i,e_j\rg)\rg|^2_{\gP}=4\, |d\al|^2_{\gP}\quad.
\ee
For any $C^1$ function $f$ we consider the following ``Frame energy''.
\[
F_f(\vec{\Phi}):= \int_\Sigma \lf[f(H_{\gP})+2^{-3}\,|A_{\gP}|_{g_{\vec{\Phi}}}^2\rg]\ dvol_{\gP}
\]
Observe then that in the case of $\Sigma=T^2$ and $f(t)=t^2$, for some global Coulomb frame $\vec{e}$ (see \cite{MR3}), we have
\be
\label{V.2}
F_{t^2}(\vec{\Phi})=\int_{T^2}[|H_{\gP}|^2+2^{-3}\,|A_{\gP}|_{g_{\vec{\Phi}}}^2]\ dvol_{\gP}=\frac{1}{4}\int_{T^2}|d\vec{e}\,|^2\ dvol_{\gP}\quad.
\ee
which is noting but the Dirichlet energy of the frame and justifies the denomination ``frame energy''.

\subsection{The first variation of Frame energies.}

We shall now compute the first variation of $F_f$. We we shall first concentrate  on the second part
\[
C(\vec{\Phi}):=\int_\Sigma|A_{\gP}|_{g_{\vec{\Phi}}}^2\ dvol_{\gP}
\]
that we call {\it connection energy}.
We observe that locally for any unit vector field $e$ for the metric $\gP$ one has
\[
|d\al|^2_{\gP}=|d\,|e|_{g_0}|^2_{g_0}
\]
We consider a perturbation $\vec{\Phi}_t:=\vec{\Phi}+t\,\vec{w}$.  Recall the Leibnitz formula
\be
\label{V.2-a}
-\Delta_{g_{\vec{\Phi}_t}}\al_t= K_{g_{\vec{\Phi}_t}}-e^{-2\al_t}\, K_{g_0}
\ee
where $\Delta_{\vec{\Phi}_t}$ is the \underbar{negative} Laplace Beltrami operator. Observe that $K_{g_0}$ is independent of $t$. We have then
\be
\label{V.3}
\frac{d}{dt}\lf[\Delta_{g_{\vec{\Phi}_t}}\al_t\rg]+\,2\, K_{g_0}\ e^{-2\al}\frac{d\al}{dt}=-\frac{dK_{g_{\vec{\Phi}_t}}}{dt}
\ee
Hence
\be
\label{V.4}
\Delta_{g_{\vec{\Phi}_t}}\frac{d\al}{dt}+\,2\, K_{g_0}\ e^{-2\al}\frac{d\al}{dt}=-\frac{d}{dt}\lf[\Delta_{g_{\vec{\Phi}_t}}\rg]\al-\frac{dK_{g_{\vec{\Phi}_t}}}{dt}
\ee
In other words, we have
\be
\label{V.3-a}
\Delta_{g_0}\frac{d\al}{dt}+\,2\, K_{g_0}\ \frac{d\al}{dt}=-\,e^{2\al}\,\frac{d}{dt}\lf[\Delta_{g_{\vec{\Phi}_t}}\rg]\al-e^{2\al}\,\frac{dK_{g_{\vec{\Phi}_t}}}{dt}
\ee

\medskip

In a local chart we have
\be
\label{V.5}
K_{g_{\vec{\Phi}_t}}\ dvol_{g_{\vec{\Phi}_t}}=\vec{n}\cdot\p_{x_1}\vec{n}\times\p_{x_2}\vec{n}\ \ dx_1\wedge dx_2
\ee
Since
\be
\label{V.6}
\ds\frac{d\vec{n}}{dt}(0)=-\lf<\vec{n}\cdot d\vec{w}\, ,\, d\vec{\Phi}\rg>_{\gP}\quad.
\ee
we deduce
\be
\label{V.7}
\begin{array}{l}
\ds\frac{d(K_{g_{\vec{\Phi}_t}}\ dvol_{g_{\vec{\Phi}_t}})}{dt}=-\lf<\vec{n}\cdot d\vec{w}\, ,\, d\vec{\Phi}\rg>_{\gP}\cdot\p_{x_1}\vec{n}\times\p_{x_2}\vec{n}\ \ dx_1\wedge dx_2\\[5mm]

\ds-\lf[\vec{n}\cdot\p_{x_1}\lf<\vec{n}\cdot d\vec{w}\, ,\, d\vec{\Phi}\rg>_{\gP}\times\p_{x_2}\vec{n} +\vec{n}\cdot\p_{x_1}\vec{n}\times\p_{x_2}\lf<\vec{n}\cdot d\vec{w}\, ,\, d\vec{\Phi}\rg>_{\gP}\rg]\ dx_1\wedge dx_2
\end{array}
\ee
We have
\be
\label{V.8}
-\lf<\vec{n}\cdot d\vec{w}\, ,\, d\vec{\Phi}\rg>_{\gP}\cdot\p_{x_1}\vec{n}\times\p_{x_2}\vec{n}=0
\ee
We choose a chart in which $\vec{\Phi}$ is conformal and we denote
\[
g_{\vP}=e^{2\la}\ [dx_1^2+dx_2^2]\quad\mbox{ and }\quad g_0=e^{2\mu}\ [dx_1^2+dx_2^2]
\] 
thus $\al=\la-\mu$. We have in one hand
\be
\label{V.9}
\begin{array}{l}
\ds-\vec{n}\cdot\p_{x_1}\lf<\vec{n}\cdot d\vec{w}\, ,\, d\vec{\Phi}\rg>_{\gP}\times\p_{x_2}\vec{n} =-\sum_{i=1}^2\p_{x_1}\lf(e^{-2{\la}}\, \vec{n}\cdot\p_{x_i}\vec{w}\rg)\ \vec{n}\cdot\p_{x_i}\vec{\Phi}\times\p_{x_2}\vec{n}\\[5mm]
\ds\quad\quad\quad-\sum_{i=1}^2e^{-2{\la}}\, \vec{n}\cdot\p_{x_i}\vec{w}\ \ \vec{n}\cdot\p^2_{x_ix_1}\vec{\Phi}\times\p_{x_2}\vec{n}\\[5mm]
\ds\quad\quad={\mathbb I}_{22}\ \p_{x_1}\lf(e^{-2{\la}}\, \vec{n}\cdot\p_{x_1}\vec{w}\rg)-{\mathbb I}_{21}\ \p_{x_1}\lf(e^{-2{\la}}\, \vec{n}\cdot\p_{x_2}\vec{w}\rg)\\[5mm]
\ds\quad\quad-\lf<\vec{n}\cdot d\vec{w},d{\la}\rg>_{\gP}\ \vec{n}\cdot\p_{x_1}\vec{\Phi}\times\p_{x_2}\vec{n}-\sum_{i=1}^2 e^{-4{\la}}\ \vec{n}\cdot\p_{x_i}\vec{w}\ \ \p^2_{x_ix_1}\vec{\Phi}\cdot\p_{x_2}\vec{\Phi}\ {\mathbb I}_{12}
\end{array}
\ee
and in the other hand
\be
\label{V.10}
\begin{array}{l}
\ds-\vec{n}\cdot\p_{x_1}\vec{n}\times\p_{x_2}\lf<\vec{n}\cdot d\vec{w}\, ,\, d\vec{\Phi}\rg>_{\gP}=\sum_{i=1}^2\p_{x_2}\lf(e^{-2{\la}}\, \vec{n}\cdot\p_{x_i}\vec{w}\rg)\ \vec{n}\cdot\p_{x_i}\vec{\Phi}\times\p_{x_1}\vec{n}\\[5mm]
\ds\quad\quad\quad+\sum_{i=1}^2e^{-2{\la}}\, \vec{n}\cdot\p_{x_i}\vec{w}\ \ \vec{n}\cdot\p^2_{x_ix_2}\vec{\Phi}\times\p_{x_1}\vec{n}\\[5mm]
\ds\quad\quad={\mathbb I}_{11}\ \p_{x_2}\lf(e^{-2{\la}}\, \vec{n}\cdot\p_{x_2}\vec{w}\rg)-{\mathbb I}_{21}\ \p_{x_2}\lf(e^{-2{\la}}\, \vec{n}\cdot\p_{x_1}\vec{w}\rg)\\[5mm]
\ds\quad\quad+\lf<\vec{n}\cdot d\vec{w},d{\la}\rg>_{\gP}\ \vec{n}\cdot\p_{x_2}\vec{\Phi}\times\p_{x_1}\vec{n}-\sum_{i=1}^2 e^{-4{\la}}\ \vec{n}\cdot\p_{x_i}\vec{w}\ \ \p^2_{x_ix_2}\vec{\Phi}\cdot\p_{x_1}\vec{\Phi}\ {\mathbb I}_{12}
\end{array}
\ee
Summing (\ref{V.8}), (\ref{V.9}) and (\ref{V.10}) gives
\be
\label{V.11}
\begin{array}{l}
\ds\frac{d(K_{g_{\vec{\Phi}_t}}\ dvol_{g_{\vec{\Phi}_t}})}{dt}=\lf[{\mathbb I}_{22}\ \p_{x_1}\lf(e^{-2\la}\, \vec{n}\cdot\p_{x_1}\vec{w}\rg)+{\mathbb I}_{11}\ \p_{x_2}\lf(e^{-2\la}\, \vec{n}\cdot\p_{x_2}\vec{w}\rg)\rg] dx_1dx_2\\[5mm]
\ds\quad -{\mathbb I}_{12}\ \lf[ \p_{x_1}\lf(e^{-2\la}\, \vec{n}\cdot\p_{x_2}\vec{w}\rg)+ \p_{x_2}\lf(e^{-2\la}\, \vec{n}\cdot\p_{x_1}\vec{w}\rg)\rg]\ dx_1dx_2\\[5mm]
\ds\quad +2\ H\ e^{2\la}\ \lf<\vec{n}\cdot d\vec{w},d\la\rg>_{\gP}\ dx_1dx_2
\end{array}
\ee
Recall Codazzi
\be
\label{V.12}
\lf\{
\begin{array}{l}
\ds \p_{x_1}{\mathbb I}_{22}-\p_{x_2}{\mathbb I}_{12}=H\ \p_{x_1}e^{2\la}\\[5mm]
\ds \p_{x_2}{\mathbb I}_{11}-\p_{x_1}{\mathbb I}_{12}=H\ \p_{x_2}e^{2\la}
\end{array}
\rg.
\ee
Hence we have proved the following lemma
\begin{Lm}
\label{lm-der-Gauss}
Under the above notations we have
\be
\label{V.13}
\begin{array}{l}
\ds\frac{d(K_{g_{\vec{\Phi}_t}}\ dvol_{g_{\vec{\Phi}_t}})}{dt}=\lf[\p_{x_1}\lf({\mathbb I}_{22}\ e^{-2\la}\, \vec{n}\cdot\p_{x_1}\vec{w}\rg)+\p_{x_2}\lf({\mathbb I}_{11}\ e^{-2\la}\, \vec{n}\cdot\p_{x_2}\vec{w}\rg)\rg] dx_1\wedge dx_2\\[5mm]
\ds\quad -\lf[ \p_{x_1}\lf({\mathbb I}_{12}\ e^{-2\la}\, \vec{n}\cdot\p_{x_2}\vec{w}\rg)+ \p_{x_2}\lf({\mathbb I}_{12}\ e^{-2\la}\, \vec{n}\cdot\p_{x_1}\vec{w}\rg)\rg]\ dx_1\wedge dx_2
\end{array}
\ee
\hfill $\Box$
\end{Lm}
Recall
\be
\label{III.8-1}
\frac{d}{dt}(dvol_{g_{\vec{\Phi}}})(0)=\lf[\sum_{i=1}^2\p_{x_i}\vec{\Phi}\cdot\p_{x_i}\vec{w}\rg]\ dx_1\wedge dx_2=\lf<d\vec{\Phi}\,;d\vec{w}\rg>_{g_{\vec{\Phi}}}\ dvol_{g_{\vec{\Phi}}}\quad.
\ee
Hence
\be
\label{V.14}
\begin{array}{l}
\ds e^{2\la}\ \frac{dK_{g_{\vec{\Phi}_t}}}{dt}=\p_{x_1}\lf({\mathbb I}_{22}\ e^{-2\la}\, \vec{n}\cdot\p_{x_1}\vec{w}\rg)+\p_{x_2}\lf({\mathbb I}_{11}\ e^{-2\la}\, \vec{n}\cdot\p_{x_2}\vec{w}\rg)\\[5mm]
\ds\quad -\p_{x_1}\lf({\mathbb I}_{12}\ e^{-2\la}\, \vec{n}\cdot\p_{x_2}\vec{w}\rg)- \p_{x_2}\lf({\mathbb I}_{12}\ e^{-2\la}\, \vec{n}\cdot\p_{x_1}\vec{w}\rg)-K_{\gP}\,\sum_{i=1}^2\p_{x_i}\vec{\Phi}\cdot\p_{x_i}\vec{w}
\end{array}
\ee
Recall
\be
\label{V.15}
\Delta_gf:=(det(g_{kl}))^{-1/2}\sum_{i,j=1}^2\p_{x_i}\lf((det(g_{kl}))^{1/2}\,g^{ij}\, \p_{x_j}f\rg)
\ee
and
\be
\label{V.16}
\frac{dg_{ij}}{dt}(0)=\p_{x_i}\vec{w}\cdot\p_{x_j}\vec{\Phi}+\p_{x_j}\vec{w}\cdot\p_{x_i}\vec{\Phi}\quad.
\ee
Since $\sum_ig_{ki}\,g^{ij}=\delta_{kj}$ and $g_{ki}=e^{2\la}\, \delta_{ki}$, we have
\be
\label{V.17}
\frac{d g^{ij}}{dt}(0)=-e^{-4\la}\ \lf[\p_{x_i}\vec{\Phi}\cdot\p_{x_j}\vec{w}+\p_{x_j}\vec{\Phi}\cdot\p_{x_i}\vec{w}\rg]\quad.
\ee
Thus in particular
\be
\label{V.18}
\begin{array}{l}
\ds\frac{d}{dt}\lf(det(g_{ij})\rg)^{1/2}=2^{-1}\lf(det(g_{ij})\rg)^{-1/2}\,e^{2\la}\,\lf[\frac{dg_{11}}{dt}+\frac{dg_{22}}{dt}\rg]\, \\
\ds\quad\quad=\sum_{i=1}^2\p_{x_i}\vec{\Phi}\cdot\p_{x_i}\vec{w}
\end{array}
\ee
and
\be
\label{V.19}
\ds\frac{d}{dt}\lf(det(g_{ij})\rg)^{-1/2}=-e^{-4\la}\ \sum_{i=1}^2\p_{x_i}\vec{\Phi}\cdot\p_{x_i}\vec{w}
\ee
Combining (\ref{V.15})...(\ref{V.19}) we obtain
\begin{Lm}
\label{lm-der-laplace}
Under the previous notations, for any function $f$ independent of $t$ on $\Sigma$ we have
\be
\label{V.20}
\begin{array}{l}
\ds\frac{d(\Delta_{g_{\vec{\Phi}_t}})f}{dt}=\lf<d<d\vec{\Phi};d\vec{w}>_{\gP},df\rg>_{\gP}-\ast_{\gP} d\ast\lf[[d\vec{\Phi}\dot{\otimes}d\vec{w}+d\vec{w}\dot{\otimes}d\vec{\Phi}]\res_{\gP}df  \rg]

\end{array}
\ee
where we have explicitly in conformal coordinates
\[
\ast_{\gP} d\ast\lf[[d\vec{\Phi}\dot{\otimes}d\vec{w}+d\vec{w}\dot{\otimes}d\vec{\Phi}]\res_{\gP}df  \rg]=e^{-2\la}\sum_{i,j=1}^2\p_{x_i}\lf(e^{-2\la}\lf(\p_{x_i}\vec{\Phi}\cdot\p_{x_j}\vec{w}+\p_{x_j}\vec{\Phi}\cdot\p_{x_i}\vec{w}\rg)\ \p_{x_j}f\rg)
\]
\hfill $\Box$
\end{Lm}
We have
\be
\label{V.21}
\begin{array}{l}
\ds\frac{d}{dt}\lf[\int_\Sigma |d\al_t|^2_{g_{\vec{\Phi}_t}}\ dvol_{g_{\vec{\Phi}_t}}\rg]=\sum_{i,j=1}^2\int_\Sigma\frac{dg^{ij}}{dt}\,\p_{x_i}\al\,\p_{x_j}\al\ dvol_{\gP}+\int_{\Sigma}|d\al|^2_{\gP}\ \frac{d dvol_{g_{\vec{\Phi}_t}}}{dt}\\[5mm]
\ds\quad\quad\quad-2\ \int_\Sigma\al\ \Delta_{\gP}\frac{d\al_t}{dt}\ dvol_{\gP}
\end{array}
\ee
We first have using (\ref{V.17})
\be
\label{V.22}
\sum_{i,j=1}^2\int_\Sigma\frac{dg^{ij}}{dt}\,\p_{x_i}\al\,\p_{x_j}\al\ dvol_{\gP}=-\,2\,\int_\Sigma\lf<d\vec{\Phi},d\al\rg>_{\gP}\cdot \lf<d\vec{w},d\al\rg>_{\gP}\ dvol_{\gP}
\ee
We then have using (\ref{III.8-1})
\be
\label{V.23}
\int_{\Sigma}|d\al|^2_{\gP}\ \frac{d dvol_{g_{\vec{\Phi}_t}}}{dt}=\int_{\Sigma}|d\al|^2_{\gP}\ \lf<d\vec{\Phi}\,;d\vec{w}\rg>_{g_{\vec{\Phi}}}\ dvol_{g_{\vec{\Phi}}}
\ee
Using now (\ref{V.4}) we have
\be
\label{V.24}
\begin{array}{l}
\ds -2\ \int_\Sigma\al\ \Delta_{{\gP}}\frac{d\al_t}{dt}\ dvol_{\gP}=-2\ \int_\Sigma\al\ \Delta_{g_0}\frac{d\al_t}{dt}\ dvol_{g_0}\\[5mm]
\ds\quad=4\ K_{g_0}\int_\Sigma\al\,\frac{d\al}{dt}\ dvol_{g_0}+2\ \int_\Sigma \al\ \frac{d (\Delta_{g_{\vec{\Phi}_t}})\al}{dt}\ dvol_{\gP}+2\int_\Sigma\al\, \frac{d K_{g_{\vec{\Phi}_t}}}{dt}\ dvol_{\gP}
\end{array}
\ee
Using lemma~\ref{lm-der-laplace} we obtain
\be
\label{V.25}
\begin{array}{l}
\ds 2\ \int_\Sigma \al\ \frac{d (\Delta_{g_{\vec{\Phi}_t}})\al}{dt}\ dvol_{\gP}=-2\ \int_{\Sigma}|d\al|^2_{\gP}\ \lf<d\vec{\Phi}\,;d\vec{w}\rg>_{g_{\vec{\Phi}}}\ dvol_{g_{\vec{\Phi}}}\\[5mm]
\ds-\,2\int_{\Sigma}\al\,\Delta_{\gP}\al\ \ \lf<d\vec{\Phi}\,;d\vec{w}\rg>_{g_{\vec{\Phi}}}\ dvol_{g_{\vec{\Phi}}}+\,4\, \int_\Sigma\lf<d\vec{\Phi},d\al\rg>_{\gP}\cdot \lf<d\vec{w},d\al\rg>_{\gP}\ dvol_{\gP}
\end{array}
\ee
Using (\ref{V.14}) we obtain
\be
\label{V.26}
\begin{array}{l}
 \ds 2\int_\Sigma\al\ \frac{d K_{g_{\vec{\Phi}_t}}}{dt}\ dvol_{\gP}=-2\sum_{i=1}^2\int_{\Sigma}{\mathbb I}_{ii}\ e^{-2\la}\ \vec{n}\cdot\p_{x_{i+1}}\vec{w}\ \p_{x_{i+1}}\al\ dx^2\\[5mm]
 \ds\quad+2\int_{\Sigma}{\mathbb I}_{12}\ e^{-2\la}\ \vec{n}\cdot\lf[\p_{x_1}\vec{w}\ \p_{x_2}\al+\p_{x_2}\vec{w}\ \p_{x_1}\al\rg] dx^2\\[5mm]
 \ds\quad-2\ \int_\Sigma \al\ K_{\gP}\ \lf<d\vec{\Phi}\,;d\vec{w}\rg>_{g_{\vec{\Phi}}}\ dvol_{g_{\vec{\Phi}}}
 \end{array}
\ee
Combining (\ref{V.24}), (\ref{V.25}) and (\ref{V.26}) gives, using that $\Delta_{\gP}\al+K_{\gP}=e^{-2\al}\,K_{g_0}$
\be
\label{V.27}
\begin{array}{l}
\ds -2\ \int_\Sigma\al\ \Delta_{\gP}\frac{d\al_t}{dt}\ dvol_{\gP}=-\,2\ \int_{\Sigma}|d\al|^2_{\gP}\ \lf<d\vec{\Phi}\,;d\vec{w}\rg>_{g_{\vec{\Phi}}}\ dvol_{g_{\vec{\Phi}}}\\[5mm]
\ds\quad\quad+\,4\, \int_\Sigma\lf<d\vec{\Phi},d\al\rg>_{\gP}\cdot \lf<d\vec{w},d\al\rg>_{\gP}\ dvol_{\gP}\\[5mm]
\ds\quad-2\sum_{i=1}^2\int_{\Sigma}{\mathbb I}_{ii}\ e^{-2\la}\ \vec{n}\cdot\p_{x_{i+1}}\vec{w}\ \p_{x_{i+1}}\al\ dx^2\\[5mm]
 \ds\quad+2\int_{\Sigma}{\mathbb I}_{12}\ e^{-2\la}\ \vec{n}\cdot\lf[\p_{x_1}\vec{w}\ \p_{x_2}\al+\p_{x_2}\vec{w}\ \p_{x_1}\al\rg] dx^2\\[5mm]
\ds\quad-2\ K_{g_0}\ \int_\Sigma \al\ e^{-2\al}\ \lf<d\vec{\Phi}\,;d\vec{w}\rg>_{g_{\vec{\Phi}}}\ dvol_{g_{\vec{\Phi}}}+4\ K_{g_0}\int_\Sigma\al\,\frac{d\al}{dt}\ dvol_{g_0}
\end{array}
\ee
Observe that
\be
\label{V.28}
\begin{array}{l}
\ds-2\ K_{g_0}\ \int_\Sigma \al\ e^{-2\al}\ \lf<d\vec{\Phi}\,;d\vec{w}\rg>_{g_{\vec{\Phi}}}\ dvol_{g_{\vec{\Phi}}}+4\ K_{g_0}\int_\Sigma\al\,\frac{d\al}{dt}\ dvol_{g_0}\\[5mm]
\ds\quad\quad=-2\, K_{g_0}\ \int_\Sigma\al\ \frac{d dvol_{g_0}}{dt}
\end{array}
\ee
and since $g_0$ is normalized to have volume 1, we have
\be
\label{V.29}
\int_\Sigma \frac{d dvol_{g_0}}{dt}=0
\ee
which is consistent with the fact that the addition of constants to $\al$ (i.e. dilations of $\vec{\Phi}$) are zero directions of the lagrangian $C$.
Combining now (\ref{V.21}), (\ref{V.22}), (\ref{V.23}) and (\ref{V.27}) we obtain the following 
\be
\label{V.30}
\begin{array}{l}
\ds\frac{d}{dt}\lf[\int_\Sigma |d\al_t|^2_{g_{\vec{\Phi}_t}}\ dvol_{g_{\vec{\Phi}_t}}\rg]=-\ \int_{\Sigma}|d\al|^2_{\gP}\ \lf<d\vec{\Phi}\,;d\vec{w}\rg>_{g_{\vec{\Phi}}}\ dvol_{g_{\vec{\Phi}}}\\[5mm]
\ds\quad\quad+\, 2\,\int_\Sigma\lf<d\vec{\Phi},d\al\rg>_{\gP}\cdot \lf<d\vec{w},d\al\rg>_{\gP}\ dvol_{\gP}+2\,\int_\Sigma\lf(\vec{\mathbb I}\res_{\gP}(\ast_{\gP}d\al)\rg)\dot{\wedge} d\vec{w}\\[5mm]
\ds\quad-2\ K_{g_0}\ \int_\Sigma \al\ e^{-2\al}\ \lf<d\vec{\Phi}\,;d\vec{w}\rg>_{g_{\vec{\Phi}}}\ dvol_{g_{\vec{\Phi}}}+4\ K_{g_0}\int_\Sigma\al\,\frac{d\al}{dt}\ dvol_{g_0}
\end{array}
\ee
where we have explicitly in positive conformal coordinates
\be
\label{V.31}
\begin{array}{l}
\ds2\,\vec{\mathbb I}\res_{\gP}(\ast_{\gP}d\al)\dot{\wedge} d\vec{w}=-2\,\lf[\sum_{i=1}^2{\mathbb I}_{ii}\ e^{-2\la}\ \vec{n}\cdot\p_{x_{i+1}}\vec{w}\ \p_{x_{i+1}}\al\rg]\ dx_1\wedge dx_2\\[5mm]
 \ds\quad+\,2\,{\mathbb I}_{12}\ e^{-2\la}\ \vec{n}\cdot\lf[\p_{x_1}\vec{w}\ \p_{x_2}\al+\p_{x_2}\vec{w}\ \p_{x_1}\al\rg] dx_1\wedge dx_2
 \end{array}
\ee
moreover $d\al/dt$ solves the following PDE 
\be
\label{V.32}
\begin{array}{l}
\ds \Delta_{g_0}\frac{d\al}{dt}+2\,K_{g_0}\, \frac{d\al}{dt}=-\, \lf<d\lf(<d\vec{\Phi};d\vec{w}>_{\gP}\rg),d\al\rg>_{g_0}\\[5mm]
\ds\quad\quad+\ast_{g_0}\, d\ast_{\gP}\lf[[d\vec{\Phi}\dot{\otimes}d\vec{w}+d\vec{w}\dot{\otimes}d\vec{\Phi}]\res_{\gP}d\al  \rg]\\[5mm]
\ds\quad\quad+ \ast_{g_0}\,d\lf[{\mathbb I}\res_{\gP}(\vec{n}\cdot\ast_{\gP}d\vec{w})\rg]+K_{\gP}\ \lf<d\vec{\Phi},d\vec{w}\rg>_{g_0}
\end{array}
\ee
Observe that
\be
\label{V.33}
\begin{array}{l}
\ds \frac{d}{dt}\lf[-\frac{1}{4}\int_\Sigma\lf[2\al_t\ e^{-2\al_t}+e^{-2\al_t}\rg]\ dvol_{\vec{\Phi}_t}\rg]=\int_\Sigma\al\, \frac{d\al}{dt}\ dvol_{g_0}\\[5mm]
\ds\quad\quad-\frac{1}{4}\int_\Sigma\lf[2\al\ e^{-2\al}+e^{-2\al}\rg]\ <d\vec{\Phi};d\vec{w}>_{\gP}\ dvol_{\gP}
\end{array}
\ee
Hence
\be
\label{V.34}
\begin{array}{l}
\ds4\,K_{g_0}\,\int_\Sigma\al\, \frac{d\al}{dt}\ dvol_{g_0}-2\,K_{g_0}\,\int_\Sigma\al\ e^{-2\al}\ <d\vec{\Phi};d\vec{w}>_{\gP}\ dvol_{\gP}\\[5mm]
\ds\quad=\frac{d}{dt}\lf[-\,K_{g_0}\,\int_\Sigma\lf[2\,\al_t\ e^{-2\al_t}+e^{-2\al_t}\rg]\ dvol_{\vec{\Phi}_t}\rg]\\[5mm]
\ds\quad+K_{g_0}\,\int_\Sigma\al\ e^{-2\al}\ <d\vec{\Phi};d\vec{w}>_{\gP}\ dvol_{\gP}
\end{array}
\ee
Observe that
\be
\label{V.35}
\int_\Sigma e^{-\al_t}\ dvol_{\vec{\Phi}_t}=\int_\Sigma dvol_{g_{0,t}}\equiv 1
\ee
Combining (\ref{V.30}), (\ref{V.34}) and (\ref{V.35}) we obtain
\begin{Lm}
\label{lm-der-energ-connection}
Under the previous notations we have
\be
\label{V.36}
\begin{array}{l}
\ds\frac{d}{dt}\lf[\int_\Sigma \lf[|d\al_t|^2_{g_{\vec{\Phi}_t}}+2\, K_{g_0}\,\al_t\ e^{-2\al_t}\rg]\, \ dvol_{g_{\vec{\Phi}_t}}\rg]\\[5mm]
\ds\quad= \int_{\Sigma}\lf(|d\al|^2_{\gP}\ \ast_{\gP}d\vec{\Phi}\rg)\dot{\wedge} d\vec{w}-\,2\, \int_\Sigma\lf(\lf<d\vec{\Phi},d\al\rg>_{\gP}\ \ast_{\gP}d\al\rg)\dot{\wedge} d\vec{w}\\[5mm]
\ds\quad+\,2\,\int_\Sigma\lf(\vec{\mathbb I}\res_{\gP}(\ast_{\gP}d\al)\rg)\dot{\wedge} d\vec{w}-K_{g_0}\,\int_\Sigma\lf(\al\ e^{-2\al}\ \ast_{\gP}d\vec{\Phi}\rg)\dot{\wedge}d\vec{w}
\end{array}
\ee
\hfill $\Box$
\end{Lm}
This lemma justifies the introduction of the following modified {\it frame energy}.
\[
\begin{array}{l}
\ds\ti{F}_f(\vec{\Phi}):= \int_\Sigma[f(H_{\gP})+2^{-3}\,|A_{\gP}|_{g_{\vec{\Phi}}}^2+K_{g_0}\,\al\ e^{-2\al}]\ dvol_{\gP}\\[5mm]
\ds\quad=\int_\Sigma[f(H_{\gP})+2^{-1}\,|d\al|_{g_{\vec{\Phi}}}^2+K_{g_0}\,\al\ e^{-2\al}]\ dvol_{\gP}
\end{array}
\]
\begin{Lm}
\label{lm-gauge-inv}
Let $\vec{\Phi}$ be an immersion of the sphere $S^2$ and $g_0$ be a metric of constant curvature equal to $4\pi$ and volume $1$ such that
there exists a function $\al$ satisfying
\[
g_{\vec{\Phi}}=e^{2\al}\, g_0
\]
then the Lagrangian
\[
L(\vec{\Phi},g_0):=\int_{S^2}2^{-1}|d\al|^2_{g_{\vec{\Phi}}}\ dvol_{g_{\vec{\Phi}}}+4\pi\ \int_{S^2}\al\ dvol_{g_0}
\]
is independent of the choice of $g_0$ and in this sense is gauge invariant for the gauge group given by the M\"obius group of positive conformal
transformations ${\mathcal M}^+(S^2)$.\hfill $\Box$
\end{Lm}
\noindent{\bf Proof of lemma~\ref{lm-gauge-inv}.}
Let $\al(t)$ and $g_0(t)$ be smooth functions such that
\[
g_{\vec{\Phi}}= e^{2\al(t)}\, g_0(t)
\]
We have
\[
\begin{array}{l}
\ds\frac{d}{dt}L(\vec{\Phi},g_0(t))=:-\int_{S^2}\Delta_{g_{\vec{\Phi}}}\al(t)\ \frac{d\al}{dt}\ dvol_{g_{\vec{\Phi}}}+4\pi\ \int_{S^2}\frac{d\al}{dt}\ [e^{-\, 2\,\al}-2\, \al\ e^{-\,2\,\al}]\ dvol_{g_{\vec{\Phi}}}\\[5mm]
\ds\quad=-\int_{S^2}\lf[\Delta_{g_0(t)}\frac{d\al}{dt}+8\pi\,\frac{d\al}{dt}\rg]\ \al(t)\ dvol_{g_{0}(t)}+4\pi\ \int_{S^2}\frac{d\al}{dt}\ dvol_{g_{0}(t)}
\end{array}
\]
Since $\int_{S^2}dvol_{g_0(t)}\equiv 1$ we have in one hand
\[
0=\frac{d}{dt}\int_{S^2}e^{-\,2\,\al(t)}\ dvol_{g_{\vec{\Phi}}}=-\, 2\ \int_{S^2}\frac{d\al}{dt}\ dvol_{g_{0}(t)}
\]
The Liouville equation gives in the other hand
\[
0=\Delta_{g_{\vec{\Phi}}}\al(t)+\, K_{g_{\vec{\Phi}}}-\, 4\,\pi\, e^{-2\,\al(t)}\ 
\]
Taking the derivative gives
\[
0=\Delta_{g_{\vec{\Phi}}}\frac{d\al}{dt}+\,8\pi\, e^{-2\,\al(t)}\, \frac{d\al}{dt}=e^{-2\,\al(t)}\ \lf[\Delta_{g_0(t)}\frac{d\al}{dt}+8\pi\,\frac{d\al}{dt}\rg]
\]
all the previous says that $L(\vec{\Phi},g_0(t))$ is independent of $t$ and lemma~\ref{lm-gauge-inv} is proved.\hfill $\Box$
We have the following definition
\begin{Dfi}
\label{aubin-gauge} Let $\vec{\Phi}$ be a weak immersion in ${\mathcal E}_{S^2}$ out of all $\al$ such that there exists $g_0$ of constant
curvature equal to $4\pi$ and volume 1 such that
\be
\label{V.36-a}
g_{\vec{\Phi}}=e^{2\al}\, g_0
\ee
we call a {\it Aubin Gauge} a choice of $\al$ and $\Psi\in $Diff$(S^2)$ such that
\be
\label{V.36-b}
\Psi^\ast g_0=\frac{g_{S^2}}{4\pi}\quad\mbox{ and }\quad \forall\ j\in\{1,2,3\}\quad\quad \int_{S^2}\ x_j\ e^{2\,\al\circ\Psi(x)}\ dvol_{S^2}=0\quad.
\ee
where $g_{S^2}$ is the standard metric on $S^2$. \hfill $\Box$
\end{Dfi}
We have the following theorem by E.Onofri. 
\begin{Th}
\label{onofri} \cite{On}
For any weak immersion $\vec{\Phi}$ of $S^2$ and any $\al$ satisfying (\ref{V.36-a}) then the following inequality holds
\be
\label{V.36-c}
\int_{S^2}2^{-1}|d\al|^2_{g_{\vec{\Phi}}}\ dvol_{g_{\vec{\Phi}}}+4\pi\ \int_{S^2}\al\ dvol_{g_0}\ge 2\pi \log\int_{S^2}e^{2\al}\ dvol_{g_0}
\ee
Moreover for any $\vec{\Phi}$ there exists a Aubin Gauge $(\Psi,\al)$ satisfying (\ref{V.36-b}).
\hfill $\Box$
\end{Th}
We are going to use the following result proved by N.Ghoussoub and C.S.Lin.
\begin{Th}
\label{ghoussoub}\cite{GhL}
For any weak immersion $\vec{\Phi}$ of $S^2$ and any $\al$ satisfying (\ref{V.36-a}) and (\ref{V.36-b}) then the following inequality holds
\be
\label{V.36-d}
\int_{S^2}3^{-1}|d\al|^2_{g_{\vec{\Phi}}}\ dvol_{g_{\vec{\Phi}}}+4\pi\ \int_{S^2}\al\ dvol_{g_0}\ge 2\pi \log\int_{S^2}e^{2\al}\ dvol_{g_0}
\ee
\hfill $\Box$
\end{Th}
It is suggested by A.Chang and P.Yang ( \cite{CY} section 3) that the constant $3^{-1}$ could be replaced by $4^{-1}$ in (\ref{V.36-d}).
\subsection{The variation of the mean curvature}
We have
\[
\frac{d}{dt}H=\frac{1}{2}\sum_{i,j}\frac{dg^{ij}}{dt}\,{\mathbb I}_{ij}+\frac{1}{2}\sum_{i,j}g^{ij}\,\frac{d{\mathbb I}_{ij}}{dt}
\]
We have
\[
\frac{dg_{ij}}{dt}=\frac{d}{dt}\lf(\p_{x_i}\vec{\Phi}_t\cdot\p_{x_j}\vec{\Phi}_t\rg)=\p_{x_i}\vec{w}\cdot\p_{x_j}\vec{\Phi}+\p_{x_i}\vec{\Phi}\cdot\p_{x_j}\vec{w}
\]
Since
\[
e^{-2\la}\, \frac{dg_{ik}}{dt} \delta_{kj}+e^{2\la}\, \frac{dg^{kj}}{dt} \delta_{ik}=0
\]
we have
\[
\frac{dg^{ij}}{dt}=-e^{-4\la}\, \lf[\p_{x_i}\vec{w}\cdot\p_{x_j}\vec{\Phi}+\p_{x_i}\vec{\Phi}\cdot\p_{x_j}\vec{w}\rg]\quad.
\]
So
\[
\frac{1}{2}\sum_{i,j}\frac{dg^{ij}}{dt}\,{\mathbb I}_{ij}=e^{-4\la}\ \frac{1}{2}\sum_{i,j}\lf[\p_{x_i}\vec{w}\cdot\p_{x_j}\vec{\Phi}+\p_{x_i}\vec{\Phi}\cdot\p_{x_j}\vec{w}\rg]\, \p_{x_i}\vec{n}\cdot\p_{x_j}\vec{\Phi}=e^{-2\la}\nabla\vec{w}\cdot\nabla\vec{\Phi}
\]
We have moreover
\[
\frac{d{\mathbb I}_{ij}}{dt}=-\frac{d}{dt}\lf(\p_{x_i}\vec{n}_t\cdot\p_{x_j}\vec{\Phi}_t\rg)=-\p_{x_i}\frac{d\vec{n}_t}{dt}\cdot\p_{x_j}\vec{\Phi}-\p_{x_i}\vec{n}\cdot\p_{x_j}\vec{w}
\]
So combining the previous assertions we have
\[
\frac{dH}{dt}=-\frac{e^{-2\la}}{2}\nabla\vec{n}\cdot\nabla\vec{w}-\frac{e^{-2\la}}{2}\nabla\frac{d\vec{n}_t}{dt}\cdot\nabla\vec{\Phi}+e^{-2\la}\nabla\vec{n}\cdot\nabla\vec{\Phi}
\]
Since 
\[
\frac{d\vec{n}}{dt}\cdot\Delta\vec{\Phi}=0
\]
we obtain 
\[
\frac{dH}{dt}=\frac{e^{-2\la}}{2}\nabla\vec{n}\cdot\nabla\vec{w}-\frac{e^{-2\la}}{2}\mbox{div}\lf(\frac{d\vec{n}_t}{dt}\cdot\nabla\vec{\Phi}\rg)
\]
One has
\[
\frac{d\vec{n}}{dt}=-\ e^{-2\la}\ \vec{n}\cdot\nabla\vec{w}\ \nabla\vec{\Phi}
\]
Hence
\[
\mbox{div}\lf(\frac{d\vec{n}_t}{dt}\cdot\nabla\vec{\Phi}\rg)=-\mbox{div}\lf(\vec{n}\cdot\nabla\vec{w}\rg)\quad.
\]
Combining the previous gives
\be
\label{V.37}
\begin{array}{l}
\ds\frac{dH}{dt}=\frac{e^{-2\la}}{2}\ \lf[\mbox{div}\lf(\vec{n}\cdot\nabla\vec{w}\rg)+\nabla\vec{n}\cdot\nabla\vec{w}\rg]\\[5mm]
\ds\quad=-\,2^{-1}d^{\ast_{\gP}}\lf(\vec{n}\cdot d\vec{w}\rg)+2^{-1}\,\lf<d\vec{w}; d\vec{n}\rg>_{\gP}
\end{array}
\ee
We then have for any $C^1$ function $f$ 
\be
\label{V.38}
\begin{array}{l}
\ds\frac{d}{dt}\lf[\int_\Sigma f(H_{\vec{\Phi}_t})\ dvol_{\vec{\Phi}_t}\rg]=-\, 2^{-1}\, \int_\Sigma f'(H)\ d^{\ast_{\gP}}\lf(\vec{n}\cdot d\vec{w}\rg)\, dvol_{\gP}\\[5mm]
\ds\quad\quad+\, 2^{-1}\,\int_\Sigma f'(H)\,\lf<d\vec{w}; d\vec{n}\rg>_{\gP}\, dvol_{\gP}+\int_\Sigma f(H)\ \lf<d\vec{w}; d\vec{\Phi}\rg>_{\gP}\, dvol_{\gP}
\end{array}
\ee
which gives
\begin{Lm}
\label{lm-der-H}
Let $f$ be a $C^1$ function, under the previous notations we have
\be
\label{V.39}
\begin{array}{l}
\ds\frac{d}{dt}\lf[\int_\Sigma f(H_{g_{\vec{\Phi}_t}})\ dvol_{g_{\vec{\Phi}_t}}\rg]= 2^{-1}\, \int_\Sigma\lf(\ast_{\gP} d[f'(H)\,\vec{n}]\rg)\dot{\wedge} d\vec{w}\\[5mm]
\ds\quad-\,\int_\Sigma\lf(f'(H)\,\ast_{\gP}d\vec{n}\rg)\dot{\wedge} d\vec{w}-\,\int_\Sigma\lf(f(H)\ \ast_{\gP} d\vec{\Phi}\rg)\dot{\wedge} d\vec{w}\quad.
\end{array}
\ee
\hfill $\Box$
\end{Lm}
Observe moreover that from (\ref{V.18}) we can also deduce the following elementary lemma
\begin{Lm}
\label{lm-area}
Under the previous notations we have
\be
\label{V.39-1}
\begin{array}{l}
\ds\frac{d}{dt}\log\lf[ \int_{\Sigma}e^{2\al}\ dvol_{g_0}\rg]=[A_{\vec{\Phi}}(\Sigma)]^{-1}\,\int_\Sigma <d\vec{\Phi},d\vec{w}>_{g_{\vec{\Phi}}}\ dvol_{g_{\vec{\Phi}}}\\[5mm]
\ds\quad\quad=-[A_{\vec{\Phi}}(\Sigma)]^{-1}\,\int_\Sigma \lf[\ast_{\gP} d\vec{\Phi}\rg]\dot{\wedge} d\vec{w}
\end{array}
\ee
where $A_{\vec{\Phi}}(\Sigma):=\int_{\Sigma}e^{2\al}\ dvol_{g_0}=\int_\Sigma dvol_{g_{\vec{\Phi}}}$ is the area of the immersion $\vec{\Phi}$. \hfill $\Box$
\end{Lm}

\subsection{The first variation of the Frame energies and conservation laws}
Combining lemma~\ref{lm-der-energ-connection} and lemma~\ref{lm-der-H} we obtain
\begin{Lm}
\label{lm-first-var-frame}
Let $\Sigma$ be a closed oriented two dimensional manifold. Let $f$ be a $C^1$ function on ${\R}$ , let $\vec{\Phi}$ be an immersion into ${\R}^3$, let $g_{\vec{\Phi}}$ be the induced metric on $\Sigma$. Let $g_0$ be the
unique constant Gauss curvature metric of volume 1 on $\Sigma$ such that there exists $\al$ with $\gP=e^{2\al}\,g_0$. The immersion is a critical point of 
\be
\label{V.40}
\ti{F}^\Lambda_f(\vec{\Phi}):=\int_\Sigma \lf[f(H)+\, 2^{-1}\, |d\al|^2_{\gP}+K_{g_0}\ \al\ e^{-2\al}\rg]\ dvol_{\gP}-\Lambda\,\log\lf(\int_{\Sigma} e^{2\al}\ dvol_{g_0}\rg)
\ee
if and only if the following conservation law (i.e. closeness of a one form) holds
\be
\label{V.41}
\begin{array}{l}
d\lf[\ast_{\gP} d[f'(H)\,\vec{n}]-\, 2\,f'(H)\, \ast_{\gP}d\vec{n}\rg.\\[5mm]
\ds\quad\quad\lf.+\, \lf[-\, 2\,f(H)+|d\al|_{\gP}^2-\,K_{g_0}\,\al\,e^{-2\al}+2\,\Lambda\,[A_{\vec{\Phi}}(\Sigma)]^{-1} \rg]\, \ast_{\gP} d\vec{\Phi}\rg.\\[5mm]
\ds\quad\quad\lf.-\,2\,\lf<d\vec{\Phi},d\al\rg>_{\gP}\ \ast_{\gP}d\al+\,2\,\vec{\mathbb I}\res_{\gP}(\ast_{\gP}d\al)\rg]=0
\end{array}
\ee
\hfill $\Box$
\end{Lm}
Assume $\vec{\Phi}$ is a critical point of $\ti{F}_f$ given by (\ref{V.40}) and denote locally
\be
\label{V.42}
\begin{array}{l}
\ds d\vec{L}:=\ast_{\gP} d[f'(H)\,\vec{n}]-\, 2\, f'(H)\, \ast_{\gP}d\vec{n}\\[5mm]
\ds\quad\quad+\, \lf[-\,2\, f(H)+|d\al|_{\gP}^2-\,K_{g_0}\,\al\,e^{-2\al}+2\,\Lambda\,[A_{\vec{\Phi}}(S^2)]^{-1}\rg]\, \ast_{\gP} d\vec{\Phi}\\[5mm]
\quad\quad\quad-\,2\,\lf<d\vec{\Phi},d\al\rg>_{\gP}\ \ast_{\gP}d\al+\,2\,\vec{\mathbb I}\res_{\gP}(\ast_{\gP}d\al)
\end{array}
\ee
In conformal coordinates this gives
\be
\label{V.43-a}
\begin{array}{l}
\ds\p_{x_1}\vec{L}=-\p_{x_2}f'(H)\,\vec{n}+f'(H)\,\p_{x_2}\vec{n}\\[5mm]
\ds\quad-\,[-\,2\, f(H)+|d\al|_{\gP}^2-\ K_{g_0}\ \al\ e^{-2\al}+2\,\Lambda\,[A_{\vec{\Phi}}(S^2)]^{-1}]\ \p_{x_2}\vec{\Phi}\\[5mm]
\ds\quad+\, 2\ \lf<d\vec{\Phi},d\al\rg >_{\gP}\ \p_{x_2}\al+\, 2\, e^{-2\la}\,\sum_{i=1}^2{\vec{\mathbb I}}_{1i}(\ast_{\gP}d\al)_i
\end{array}
\ee
and
\be
\label{V.43-b}
\begin{array}{l}
\ds\p_{x_2}\vec{L}=\p_{x_1}f'(H)\,\vec{n}-f'(H)\,\p_{x_1}\vec{n}\\[5mm]
\ds\quad+\,[-\,2\, f(H)+|d\al|_{\gP}^2-\ K_{g_0}\ \al\ e^{-2\al}+2\,\Lambda\,[A_{\vec{\Phi}}(S^2)]^{-1}]\ \p_{x_1}\vec{\Phi}\\[5mm]
\ds\quad-\, 2\ \lf<d\vec{\Phi},d\al\rg >_{\gP}\ \p_{x_1}\al+\, 2\, e^{-2\la}\,\sum_{i=1}^2{\vec{\mathbb I}}_{2i}(\ast_{\gP}d\al)_i
\end{array}
\ee

We have
\be
\label{V.44}
\begin{array}{l}
\ds d\vec{\Phi}\dot{\wedge}d\vec{L}:=\lf[\p_{x_1}\vec{\Phi}\cdot\p_{x_2}\vec{L}-\p_{x_2}\vec{\Phi}\cdot\p_{x_1}\vec{L}\rg]\ dx_1\wedge dx_2\\[5mm]
\ds\quad=2\, (f'(H)\, H-2\, f(H)-K_{g_0}\,\al\ e^{-2\al}+2\,\Lambda\,[A_{\vec{\Phi}}(S^2)]^{-1})\ dvol_{\gP}
\end{array}
\ee
We have also
\be
\label{V.45}
\begin{array}{l}
\ds d\vec{\Phi}{\wedge}d\vec{L}:=\lf[\p_{x_1}\vec{\Phi}\times\p_{x_2}\vec{L}-\p_{x_2}\vec{\Phi}\times\p_{x_1}\vec{L}\rg]\ dx_1\wedge dx_2\\[5mm]
\ds\quad=\p_{x_1}\vec{\Phi}\times\lf[\p_{x_1}f'(H)\,\vec{n}-f'(H)\,\p_{x_1}\vec{n}-\,2\,\lf<d\vec{\Phi},d\al\rg>_{\gP}\ \p_{x_1}\al\rg]\ dx_1\wedge dx_2\\[5mm]
\ds\quad+2\,e^{-2\la}\,\sum_{i=1}^2\p_{x_1}\vec{\Phi}\times{\vec{\mathbb I}}_{2i}(\ast_{\gP}d\al)_i\ dx_1\wedge dx_2\\[5mm]
\ds\quad+\p_{x_2}\vec{\Phi}\times\lf[\p_{x_2}f'(H)\,\vec{n}-f'(H)\,\p_{x_2}\vec{n}-\,2\,\lf<d\vec{\Phi},d\al\rg>_{\gP}\ \p_{x_2}\al\rg]\ dx_1\wedge dx_2\\[5mm]
\ds\quad-\,2\,e^{-2\la}\,\sum_{i=1}^2\p_{x_2}\vec{\Phi}\times{\vec{\mathbb I}}_{1i}(\ast_{\gP}d\al)_i\ dx_1\wedge dx_2
\end{array}
\ee
This gives
\be
\label{V.46}
\begin{array}{l}
\ds d\vec{\Phi}{\wedge}d\vec{L}:=[\p_{x_1}\vec{\Phi}\,\p_{x_2}f'(H)-\p_{x_2}\vec{\Phi}\,\p_{x_1}f'(H)]\ dx_1\wedge dx_2\\[5mm]
\ds\quad\quad-\,2\, e^{-2\la}\,\sum_{i,j=1}^2{\mathbb I}_{ij}(\ast_{\gP}d\al)_j\ \p_{x_i}\vec{\Phi}\ dx_1\wedge dx_2
\end{array}
\ee
We have
\be
\label{V.47}
\begin{array}{l}
\ds-\,2\, e^{-2\la}\,\sum_{i,j=1}^2{\mathbb I}_{ij}(\ast_{\gP}d\al)_j\ \p_{x_i}\vec{\Phi}\ =2\, \p_{x_1}\al \,\lf[e^{-2\la}\,{\mathbb I}_{12}\,\p_{x_1}\vec{\Phi}+e^{-2\la}\,{\mathbb I}_{22}\,\p_{x_2}\vec{\Phi}\rg]\\[5mm]
\ds\quad\quad\quad\quad\quad-2\, \p_{x_2}\al \,\lf[e^{-2\la}\,{\mathbb I}_{11}\,\p_{x_1}\vec{\Phi}+e^{-2\la}\,{\mathbb I}_{12}\,\p_{x_2}\vec{\Phi}\rg]
\end{array}
\ee
We compute, using Codazzi identity (\ref{V.12})
\be
\label{V.48-a}
\begin{array}{l}
\ds e^{2\la}\,\p_{x_1}\lf[e^{-2\la}\,{\mathbb I}_{12}\,\p_{x_1}\vec{\Phi}+e^{-2\la}\,{\mathbb I}_{22}\,\p_{x_2}\vec{\Phi}\rg]-\,e^{2\la}\,\p_{x_2}\lf[e^{-2\la}\,{\mathbb I}_{11}\,\p_{x_1}\vec{\Phi}+e^{-2\la}\,{\mathbb I}_{12}\,\p_{x_2}\vec{\Phi}\rg]\\[5mm]
\ds=\p_{x_1}\vec{\Phi}\ \lf[-2\,\p_{x_1}\la\ {\mathbb I}_{12}+\p_{x_1}{\mathbb I}_{12} + e^{-2\la}\,{\mathbb I}_{12}\, \p^2_{x_1^2}\vec{\Phi}\cdot\p_{x_1}\vec{\Phi}+e^{-2\la}\,{\mathbb I}_{22}\, \p^2_{x_1x_2}\vec{\Phi}\cdot\p_{x_1}\vec{\Phi}\rg]\\[5mm]
\ds+\p_{x_2}\vec{\Phi}\ \lf[-2\,\p_{x_1}\la\ {\mathbb I}_{22}+\p_{x_1}{\mathbb I}_{22}+e^{-2\la}\,{\mathbb I}_{22}\, \p^2_{x_1x_2}\vec{\Phi}\cdot\p_{x_2}\vec{\Phi}+e^{-2\la}\,{\mathbb I}_{12}\, \p^2_{x_1^2}\vec{\Phi}\cdot\p_{x_2}\vec{\Phi}\rg]\\[5mm]
\ds-\p_{x_1}\vec{\Phi}\ \lf[-2\,\p_{x_2}\la\ {\mathbb I}_{11}+\p_{x_2}{\mathbb I}_{11} + e^{-2\la}\,{\mathbb I}_{11}\, \p^2_{x_1 x_2}\vec{\Phi}\cdot\p_{x_1}\vec{\Phi}+e^{-2\la}\,{\mathbb I}_{12}\, \p^2_{x_2^2}\vec{\Phi}\cdot\p_{x_1}\vec{\Phi}\rg]\\[5mm]
\ds-\p_{x_2}\vec{\Phi}\ \lf[-2\,\p_{x_2}\la\ {\mathbb I}_{12}+\p_{x_2}{\mathbb I}_{12}+e^{-2\la}\,{\mathbb I}_{12}\, \p^2_{x_2^2}\vec{\Phi}\cdot\p_{x_2}\vec{\Phi}+e^{-2\la}\,{\mathbb I}_{11}\, \p^2_{x_1 x_2}\vec{\Phi}\cdot\p_{x_2}\vec{\Phi}\rg]
\end{array}
\ee
Making use of Codazzi identity (\ref{V.12}) we finally obtain
\be
\label{V.48}
\begin{array}{l}
\ds e^{2\la}\,\p_{x_1}\lf[e^{-2\la}\,{\mathbb I}_{12}\,\p_{x_1}\vec{\Phi}+e^{-2\la}\,{\mathbb I}_{22}\,\p_{x_2}\vec{\Phi}\rg]- e^{2\la}\,\p_{x_2}\lf[e^{-2\la}\,{\mathbb I}_{11}\,\p_{x_1}\vec{\Phi}+e^{-2\la}\,{\mathbb I}_{12}\,\p_{x_2}\vec{\Phi}\rg]\\[5mm]
\ds=\p_{x_1}\vec{\Phi}\ \lf[-H\, \p_{x_2}e^{2\la}+\p_{x_2}\la\ [{\mathbb I}_{22}+{\mathbb I}_{11}]\rg]+\p_{x_2}\vec{\Phi}\ \lf[H\, \p_{x_1}e^{2\la}-\p_{x_1}\la\ [{\mathbb I}_{22}+{\mathbb I}_{11}]\rg]\\[5mm]
\ds=0\quad.
\end{array}
\ee
Hence there exists locally $\vec{D}$ such that (see also lemma III.2 in \cite{MR3})
\be
\label{V.49}
\lf\{
\begin{array}{l}
\ds \p_{x_1}\vec{D}:= \lf[e^{-2\la}\,{\mathbb I}_{11}\,\p_{x_1}\vec{\Phi}+e^{-2\la}\,{\mathbb I}_{12}\,\p_{x_2}\vec{\Phi}\rg]\\[5mm]
\ds\p_{x_2}\vec{D}:= \lf[e^{-2\la}\,{\mathbb I}_{12}\,\p_{x_1}\vec{\Phi}+e^{-2\la}\,{\mathbb I}_{22}\,\p_{x_2}\vec{\Phi}\rg]
\end{array}
\rg.
\ee
Combining all the previous we obtain the following lemma
\begin{Lm}
\label{cons-laws}
Let $\Sigma$ be a closed two dimensional manifold. Let $f$ be a $C^1$ function on ${\R}$ , let $\vec{\Phi}$ be an immersion into ${\R}^3$, let $g_{\vec{\Phi}}$ be the induced metric on $\Sigma$. Let $g_0$ be the
unique constant Gauss curvature metric of volume 1 on $\Sigma$ such that there exists $\al$ with $\gP=e^{2\al}\,g_0$. Assume  the immersion $\vec{\Phi}$ is a critical point of 
\be
\label{V.50}
\ti{F}^\Lambda_f(\vec{\Phi}):=\int_\Sigma \lf[f(H)+\, 2^{-1}\, |d\al|^2_{\gP}+K_{g_0}\ \al\ e^{-2\al}\rg]\ dvol_{\gP}-\Lambda\,\log\lf(\int_{\Sigma} e^{2\al}\ dvol_{g_0}\rg)
\ee
and following lemma~\ref{lm-first-var-frame}  introduce locally $\vec{L}$
\be
\label{V.51}
\begin{array}{l}
d\vec{L}:=\ast_{\gP} d[f'(H)\,\vec{n}]-\, 2\, f'(H)\, \ast_{\gP}d\vec{n}\\[5mm]
\ds\quad\quad+\, \lf[-\,2\, f(H)+|d\al|_{\gP}^2-\,K_{g_0}\,\al\,e^{-2\al}+2\,\Lambda\,[A_{\vec{\Phi}}(S^2)]^{-1}\rg]\, \ast_{\gP} d\vec{\Phi}\\[5mm]
\ds\quad\quad\quad-\,2\,\lf<d\vec{\Phi},d\al\rg>_{\gP}\ \ast_{\gP}d\al+\,2\,\vec{\mathbb I}\res_{\gP}(\ast_{\gP}d\al)
\end{array}
\ee
then the following almost conservation law holds
\be
\label{V.52}
d\vec{\Phi}\dot{\wedge}d\vec{L}=2\, (f'(H)\, H-2\, f(H)-K_{g_0}\,\al\ e^{-2\al}+2\,\Lambda\,[A_{\vec{\Phi}}(S^2)]^{-1})\ dvol_{\gP}
\ee
and the following exact conservation law holds
\be
\label{V.53}
d\vec{\Phi}\wedge d\vec{L}=d\vec{\Phi}\wedge df'(H)+ 2\,d\al\wedge d\vec{D}
\ee
where $\vec{D}$ satisfies
\be
\label{V.54}
d\vec{D}= {\mathbb I}\res_g d\vec{\Phi}\quad.
\ee
\hfill $\Box$
\end{Lm}
\begin{Rm}
The three conservation laws or almost conservation laws (\ref{V.51}), (\ref{V.52}) and (\ref{V.53}) can be deduced from Noether theorem (see \cite{Be}). Precisely
the existence of $\vec{L}$ satisfying the first conservation law (\ref{V.51}) is due to the translation invariance of the Lagrangian $\ti{F}^\Lambda_f$, (\ref{V.52}) instead is related to the lack
of invariance of the lagrangian  under dilation whereas (\ref{V.53}) is related  to the rotation invariance of the Lagrangian.
\hfill $\Box$
\end{Rm}
\subsection{Various Bounds involving the Frame Energies.}
First of all we establish the following lemma.
\begin{Lm}
\label{lm-alpha-bounds}
Under the previous notations we have
\be
\label{ex-1}
\begin{array}{l}
\ds\lf|\frac{\ds\log\lf(\int_{S^2}e^{2\al}\ dvol_{g_0}\rg)}{\ds\log\lf(\frac{1}{\sigma}\rg)}\rg|\le 2+\frac{\ds\log\lf(\sigma^2\int_{S^2}(1+H^2)^2\ dvol_{\gP}\rg)}{\log\lf(\ds\frac{1}{\sigma}\rg)}\\[5mm]

\end{array}
\ee
\end{Lm}
\noindent{\bf Proof of lemma~\ref{lm-alpha-bounds}.}
We have obviously
\be
\label{ex-2}
\log\lf(\int_{S^2}e^{2\al}\ dvol_{g_0}\rg)\le \log\lf(\sigma^2\int_{S^2}(1+H^2)^2\ dvol_{\gP}\rg)+2\,\log\lf(\ds\frac{1}{\sigma}\rg)
\ee
We have also
\be
\label{V.60}
16\,\pi^2\le \lf(\int_{S^2}H_{\vec{\Phi}}^2\ dvol_{\gP}\rg)^2\le\int_{S^2}e^{2\al}\ dvol_{g_0}\ \int_{S^2}H_{\vec{\Phi}}^4\ dvol_{\gP}
\ee
Hence
\be
\label{ex-3}
2\,\log\lf(\ds\frac{1}{\sigma}\rg)+\log\lf(\sigma^2\int_{S^2}(1+H^2)^2\ dvol_{\gP}\rg)\ge-\log\lf(\int_{S^2}e^{2\al}\ dvol_{g_0}\rg)+\log 16\pi
\ee
Combining (\ref{ex-2}) and (\ref{ex-3}) gives (\ref{ex-1}) and lemma~\ref{lm-alpha-bounds} is proved.\hfill $\Box$

\medskip
A useful lemma is the following which is a direct consequence of theorem~\ref{ghoussoub}
\begin{Lm}
\label{dirichlet-bound}
Let $\vec{\Phi}$ be a weak immersion of $S^2$ and $(\al,\Psi)$ be an Aubin Gauge satisfying (\ref{V.36-b}) then the following inequality holds
\be
\label{V.62}
\begin{array}{l}
\ds 6^{-1}\int_{S^2}|d\al|^2_{g_{0}}\, dvol_{g_0}\\[5mm]
\ds\le  \int_{S^2}2^{-1}|d\al|^2_{g_{\vec{\Phi}}}\ dvol_{g_{\vec{\Phi}}}+4\pi\ \int_{S^2}\al\ dvol_{g_0}- 2\pi \log\int_{S^2}e^{2\al}\ dvol_{g_0}
\end{array}
\ee
\hfill $\Box$
\end{Lm}

\medskip

For $\Sigma$ being an arbitrary closed surface we denote
\be
\label{V.0-0}
\begin{array}{l}
\ds F^\sigma(\vec{\Phi}):= \lf(\log\frac{1}{\sigma}\rg)^{-1}\,\ti{F}^{K_{g_0}}_{f_\sigma}(\vec{\Phi})\\[5mm]
\ds\quad\quad:=\lf(\log\frac{1}{\sigma}\rg)^{-1}\,\int_{\Sigma}\lf[f_\sigma(H_{\gP})+\, 2^{-1}\, |d\al|^2_{\gP}+K_{g_0}\ \al\ e^{-2\al}\rg]\ dvol_{\gP}\\[5mm]
\ds\quad-\,2^{-1}\, K_{g_0}\  \lf(\log\frac{1}{\sigma}\rg)^{-1}\,\log\lf(\int_{\Sigma}e^{2\al}\ dvol_{g_0}\rg)
\end{array}
\ee
where $f_\sigma(t)=\log\sigma^{-1}[ t^2+\sigma^2 (1+t^2)^2]$. First, in the particular case $\Sigma=S^2$ we have
\begin{Lm}
\label{lm-L-infty-borne}
Let $\vec{\Phi}$ be a weak immersion of $S^2$ in ${\mathcal E}_{S^2,2}$ and $g_0$ be a constant Gauss curvature metric on $S^2$ of volume 1
such that
\[
g_{\vP}=e^{2\al}\, g_0
\]
We have for $\sigma\in (0,1)$
\be
\label{V.67}
W(\vec{\Phi})+\sigma^2\int_{S^2}(1+H^2)^2\ dvol_{g_{\vec{\Phi}}}\le  F^\sigma(\vec{\Phi})
\ee
and for any Aubin gauge we have
\be
\label{V.67-a1}
\lf(\log\frac{1}{\sigma}\rg)^{-1}\, \int_{S^2} |d\al|^2_{\gP}\le 6\ \lf[F^\sigma(\vec{\Phi})-W(\vec{\Phi})\rg]
\ee
We have also for any gauge
\be
\label{V.69}
\inf_{x\in S^2}\al\ge \frac{1}{2}\log A_{\vec{\Phi}}(S^2)- \frac{1}{2\pi}\,\int_{S^2}|d\al|^2_{g_0}\, dvol_{g_0}\ -C\ \lf[ 1+W(\vec{\Phi})\rg]
\ee
and
\be
\label{V.69-a}
\lf(\log\frac{1}{\sigma}\rg)^{-1}\,\|\al\|_{L^\infty(S^2)}\le C\, \lf(\log\frac{1}{\sigma}\rg)^{-1}\ \lf[ \int_{S^2} |d\al|^2_{\gP}+W(\vec{\Phi})\rg]+\lf|\frac{\log A_{\vec{\Phi}}(S^2)}{2\,\log\frac{1}{\sigma}}\rg|
\ee
moreover
\be
\label{V.69-b}
\frac{\log A_{\vec{\Phi}}(S^2)}{2\,\log\frac{1}{\sigma}}\le 1+\lf(\log\frac{1}{\sigma}\rg)^{-1}\log\lf[F^\sigma(\vec{\Phi})-W(\vec{\Phi})\rg]
\ee
where the constant $C$ is a positive universal constant.\hfill $\Box$
\end{Lm}
\noindent{\bf Proof of lemma~\ref{lm-L-infty-borne}.}
First of all, Onofri inequality (\ref{V.36-c}) implies
\be
\label{V.67-a}
F^\sigma(\vec{\Phi})\ge W(\vec{\Phi})+\sigma^2\int_{S^2}(1+H^2)^2\ dvol_{g_{\vec{\Phi}}}
\ee
We have also
\be
\label{V.67-b}
\begin{array}{l}
\ds F^\sigma(\vec{\Phi})-W(\vec{\Phi})\ge\lf(\log\frac{1}{\sigma}\rg)^{-1}\,\int_{S^2}\lf[\, 2^{-1}\, |d\al|^2_{\gP}+K_{g_0}\ \al\ e^{-2\al}\rg]\ dvol_{\gP}\\[5mm]
\ds\quad\quad\quad\quad-\,2\pi\ \lf(\log\frac{1}{\sigma}\rg)^{-1}\,\log\lf(\int_{S^2}e^{2\al}\ dvol_{g_0}\rg)
\end{array}
\ee
Then we obtain (\ref{V.67-a1}) directly from (\ref{V.62}).
\medskip

The uniformization theorem on $S^2$ gives the existence of a diffeomorphism $\Psi$ such that $\Psi^\ast g_0=g_{S^2}/4\pi$ where $g_{S^2}$ is 
the standard metric on $S^2$. We then consider the immersion $\vec{\Phi}\circ\Psi$. 

\medskip

Following \cite{Ri1} we have the existence of a constant $C(S^2)$ such that
\be
\label{V.70}
\|d\al\|_{L^{2,\infty}}\le C(S^2)\ \lf[1+W(\vec{\Phi})\rg]
\ee
where the norm is taken with respect to the metric $g_0$.    We cover $S^2$ by a finite, controlled family of geodesic convex balls for the  metric $g_0$ and each of these balls we choose a conformal chart
$\Psi_i\ :\ D^2\, \rightarrow \Sigma$ such that $\Psi_i(D^2_{1/2})$ still cover $S^2$. We consider any of these ball and we omit to write explicitly the composition of $\vec{\Phi}$ with $\Psi_i$. Denote $\vec{e}_j:=e^{-\la}\,\p_{x_j}\vec{\Phi}$ where $e^\la=\log|\p_{x_j}\vec{\Phi}|$. $(\vec{e}_1,\vec{e}_2)$ realizes a moving frame and
\[
|\nabla\vec{e}_j|^2=|\vec{n}\cdot\nabla\vec{e}_j|^2+|\nabla\la|^2\quad.
\]
We have
\[
-\Delta\al=\Delta\mu+(\nabla\vec{e}_1,\nabla\vec{e}_2)
\]
Using Wente together with classical elliptic estimates we get the existence of $\ov{\al}\in {\R}$ such that
\be
\label{V.71}
\|\al-\ov{\al}\|_{L^\infty(D^2_{3/4})}\le\, C\, \lf[\|\nabla\al\|_{L^{2,\infty}(D^2)}+\|\nabla\mu\|_{L^2(D^2)}+\|\nabla\vec{n}\|_{L^2(D^2)}^2\rg]+\frac{1}{2\pi}\,\int_{S^2}|d\al|^2_{g_0}\, dvol_{g_0}
\ee
Since $\int_\Sigma dvol_{g_0}=1$, there exists $x\in\Sigma$ such that
\be
\label{V.73}
\al(x):=\frac{1}{2}\log \mbox{Area}(\vec{\Phi}(S^2))
\ee
Combining (\ref{V.71}) and (\ref{V.73}) we get (\ref{V.69}) and lemma~\ref{lm-L-infty-borne} is proved.\hfill $\Box$
First we shall need the following lemma
\begin{Lm}
\label{lm-bound-bel-genus}
Let $\Sigma$ be a closed oriented surface and $\vec{\Phi}$ be a weak immersion in ${\mathcal E}_{\Sigma,2}$. Let $g_0$ be a constant Gauss curvature metric of volume 1 
on $\Sigma$ such that there exists $\al$ for which
\[
g_{\vP}=e^{2\al}\,g_0
\]
Then, if genus$(\Sigma)\ge 1$ we have
\be
\label{V.55-aa}
K_{g_0}\ \int_\Sigma\al \ dvol_{g_0}\ge 2^{-1}\,K_{g_0}\,\log\int_\Sigma\ dvol_{\gP}
\ee
\end{Lm}
\noindent{\bf Proof of lemma~\ref{lm-bound-bel-genus}}
 We have to bound from bellow the following quantity
\[
4\pi\, (1-\mbox{genus}(\Sigma))\, \int_\Sigma \ \al\ dvol_{g_0}\quad.
\]
Since we fixed $\int_\Sigma dvol_{g_0}=1$, using the convexity of $\exp$ and Jensen inequality we have
\be
\label{V.56}
\exp\lf(2\, \int_\Sigma \ \al\ dvol_{g_0}\rg)\le \int_{\Sigma}\ e^{2\al}\ dvol_{g_0}=\int_\Sigma\ dvol_{\gP}
\ee
We then have for $K_{g_0}<0$
\be
\label{V.57}
K_{g_0}\ \int_\Sigma\al \ dvol_{g_0}\ge 2\pi\ (1-\mbox{genus}(\Sigma))\,\log\int_\Sigma\ dvol_{\gP}
\ee
This concludes the proof of lemma~\ref{lm-bound-bel-genus}.\hfill $\Box$

\medskip
 
 Combining the previous lemma with the proof of lemma~\ref{lm-L-infty-borne} we obtain the following lemma
  \begin{Lm}
\label{lm-L-infty-borne-genre-neg}
Let $\Sigma$ be a closed surface of non zero genus and $\vec{\Phi}$ be a weak immersion of $\Sigma$ in ${\mathcal E}_{\Sigma,2}$ and $g_0$ be the constant Gauss curvature metric on $\Sigma$ of volume 1
such that
\[
g_{\vP}=e^{2\al}\, g_0
\]
We have for $\sigma\in (0,1)$
\be
\label{V.67-agen}
W(\vec{\Phi})+\sigma^2\int_{\Sigma}(1+H^2)^2\ dvol_{g_{\vec{\Phi}}}+\lf(2\,\log\frac{1}{\sigma}\rg)^{-1}\, \int_{\Sigma} |d\al|^2_{\gP}\le  F^\sigma(\vec{\Phi})
\ee
and
\be
\label{V.69-agen}
\lf(\log\frac{1}{\sigma}\rg)^{-1}\,\|\al\|_{L^\infty(\Sigma)}\le C\, \lf(\log\frac{1}{\sigma}\rg)^{-1}\ \lf[ \int_{\Sigma} |d\al|^2_{\gP}+W(\vec{\Phi})\rg]+\lf|\frac{\log A_{\vec{\Phi}}(\Sigma)}{2\,\log\frac{1}{\sigma}}\rg|
\ee
where the constant $C$ is a positive universal constant.\hfill $\Box$
\end{Lm}

\section{Uniform regularity for critical points to the Frame Energies approximating the Willmore Energy}
\reset
In this section we shall consider critical points of the following family of lagrangians where the parameter $\sigma$ belongs to $[0,\sigma_0)$.
\be
\label{VI.0-0}
\begin{array}{l}
\ds F^\sigma(\vec{\Phi}):= \lf(\log\frac{1}{\sigma}\rg)^{-1}\,\ti{F}^{K_{g_0}}_{f_\sigma}(\vec{\Phi})\\[5mm]
\ds\quad\quad:=\lf(\log\frac{1}{\sigma}\rg)^{-1}\,\int_{\Sigma}\lf[f_\sigma(H_{\gP})+\, 2^{-1}\, |d\al|^2_{\gP}+K_{g_0}\ \al\ e^{-2\al}\rg]\ dvol_{\gP}\\[5mm]
\ds\quad-\,2^{-1}\, K_{g_0}\  \lf(\log\frac{1}{\sigma}\rg)^{-1}\,\log\lf(\int_{\Sigma}e^{2\al}\ dvol_{g_0}\rg)
\end{array}
\ee
where $f_\sigma(t)=\log\sigma^{-1}[ t^2+\sigma^2 (1+t^2)^2]$. Recall that we denote for any domain $U$ of $S^2$ 
\[
A_{\vec{\Phi}}(U):=\mbox{Area}(\vec{\Phi}(U))\quad.
\]
The goal of this subsection is to establish the following result
\begin{Lm}
\label{lm-ep-reg} {\bf [uniform $\ep-$regularity]}
For any $C_0>0$, there exists $\ep>0$ and $\sigma_0>0$ such that for any $\sigma\in(0,\sigma_0)$ and any critical point $\vec{\Phi}$ of $F^\sigma$ satisfying
\be
\label{VI.0}
F^\sigma(\vec{\Phi})\le C_1\quad
\ee
denoting  $l_\sigma:=\lf(\log\frac{1}{\sigma}\rg)^{-1}$, assume moreover that we have a gauge $\al$ such that
\be
\label{VI.0-a}
\inf_{D^2}\al\ge\log\sigma- C_1
\ee
and
\be
\label{VI.2}
l_\sigma\,\|\al\ e^{2\,\mu}\|_{L^\infty(D^2)}+\int_{D^2}|\nabla\vec{n}|^2+\sigma^2\,H^4\ e^{2\la}\ dx^2<\ep
\ee
then for any $j\in {\N}$ one has
\be
\label{VI.3}
\begin{array}{l}
\ds|\nabla^{j+1}\vec{n}|^2(0)+| e^\la\,\nabla^j(\vec{H}\,(1+2\,\sigma^2\,(1+H^2))|^2(0)+\sigma^4\, H^2\, (1+H^2)^2\,e^{2\la}(0)\\[5mm]
\ds\quad\le   C_j\, \int_{D^2}|\nabla\vec{n}|^2\ dx^2+C_j\ \lf[\sigma^2\,\int_{D^2}H^4\,e^{2\la}\ dx^2\rg]^2+C_j\,\lf[l_\sigma\,\int_{D^2}|\nabla\al|^2\ dx^2 \rg]^2 \\[5mm]
\ds\quad\quad+C_j\, [l^2_\sigma\,|\ov{\al}|^2+ l^2_\sigma\,\|e^{4\mu}\|_{L^\infty(D^2)}]\,\|e^{4\mu}\|_{L^\infty(D^2)}+C_j\, l^2_\sigma\,\lf[\frac{A_{\vec{\Phi}}(D^2)}{A_{\vec{\Phi}}(\Sigma)}\rg]^2
\end{array}
\ee
where $\ov{\al}=|D^2_{1/2}|^{-1}\,\int_{D^2_{1/2}}\al$ and
\be
\label{VI.3-a}
\begin{array}{l}
\ds l_\sigma |\nabla^{j+1}\al|^2(0)\le\ C_j\,l_\sigma\,\int_{D^2}|\nabla\al|^2\ dx^2 +C_j\, \lf[\int_{D^2}|\nabla\vec{n}|^2\ dx^2\rg]^2\\[5mm]
\ds\quad\quad+C_j\ \lf[\sigma^2\,\int_{D^2}H^4\,e^{2\la}\ dx^2\rg]^4+C_j\, l_\sigma\,\|e^{4\mu}\|_{L^\infty(D^2)}+C_j\, l^4_\sigma\,\lf[\frac{A_{\vec{\Phi}}(D^2)}{A_{\vec{\Phi}}(\Sigma)}\rg]^4
\end{array}
\ee
where $C_j$ only depends on $C_1$ and $j$.
\hfill $\Box$
\end{Lm}
\noindent{\bf Proof of lemma~\ref{lm-ep-reg}.}
In the first part of the proof, following the original ideas of \cite{Ri2}, we derive from equation (\ref{V.44}) more conservation laws. We do it first formally, not worrying of
the regularity. In the second part of the proof we will revisit each step with estimates in relevant Banach spaces.

\medskip

\noindent{\bf Step 1}. Equations.
Let $l_\sigma:=\lf(\log\frac{1}{\sigma}\rg)^{-1}$
\be
\label{VI.3a}
\begin{array}{l}
\ds\nabla\vec{L}:=l_\sigma\nabla^\perp (f'_\sigma(H)\, \vec{n})-2\, l_\sigma\, f'_\sigma(H)\,\nabla^\perp\vec{n}-2\ e^{-2\la}\, l_\sigma\, \nabla\vec{\Phi}\cdot\nabla\al\ \nabla^\perp\al\\[5mm]
\ds\quad+l_\sigma\ \lf[-2\, f_\sigma(H)+e^{-2\la}\,|\nabla\al|^2-K_{g_0}\, \al\ e^{-2\al}+K_{g_0}\,A_{\vec{\Phi}}(\Sigma)^{-1}\rg]\ \nabla^\perp\vec{\Phi}+2\, l_\sigma\ e^{-2\la}\, \vec{\mathbb I}\res\nabla^\perp\al
\end{array}
\ee
 Equation (\ref{V.44}) gives
\be
\label{VI.4}
\nabla\vec{\Phi}\cdot\nabla^{\perp}\vec{L}=2\,l_\sigma\ e^{2\la}\, \lf(2\,f_\sigma(H)-\,H\, f'_\sigma(H)+K_{g_0}\, \al\, e^{-2\al}-K_{g_0}\,A_{\vec{\Phi}}(\Sigma)^{-1}\rg)
\ee
Let $Y$ be the solution of
\be
\label{VI.5-a}
\lf\{
\begin{array}{l}
\ds-\Delta Y=2\, l_\sigma\, e^{2\la}\, \lf(2\,f_\sigma(H)-\,H\, f'_\sigma(H)+K_{g_0}\, \al\, e^{-2\al}-K_{g_0}\,A_{\vec{\Phi}}(\Sigma)^{-1}\rg)\quad\mbox{ in }D^2\\[5mm]
\ds Y=0\quad\quad\mbox{ on }\p D^2
\end{array}
\rg.
\ee
Observe that $2\,f_\sigma(H)-\,H\, f'_\sigma(H)= 2\, \,l_\sigma^{-1}\, \sigma^2\,  (1-H^4)$. So $Y$ satisfies
\be
\label{VI.5}
\lf\{
\begin{array}{l}
\ds-\Delta Y=4\, e^{2\la}\, \sigma^2\, (1-H^4)+2\,l_\sigma\,K_{g_0}\, \al\, e^{2\mu}-2\, K_{g_0}\,l_\sigma\,e^{2\la}\,A_{\vec{\Phi}}(\Sigma)^{-1}\quad\mbox{ in }D^2\\[5mm]
\ds Y=0\quad\quad\mbox{ on }\p D^2
\end{array}
\rg.
\ee
Using Poincar\'e theorem we deduce the existence of a function $S$ such that
\be
\label{VI.6}
\nabla S=\vec{L}\cdot\nabla\vec{\Phi}+\nabla^\perp Y
\ee
The equation (\ref{V.53}) in conformal coordinates gives
\be
\label{VI.7}
\nabla\vec{\Phi}\times\nabla^\perp\vec{L}=-\,l_\sigma\,\nabla^\perp\vec{\Phi}\cdot\nabla f'_\sigma(H)+2\,\nabla\al\cdot\nabla^\perp\vec{D}
\ee
where
\be
\label{VI.7-a}
\nabla\vec{D}=\lf(l_\sigma\ e^{-2\la} \sum_{i=1}^2{\mathbb I}_{1i}\,\p_{x_i}\vec{\Phi}\ ,\ l_\sigma\ e^{-2\la} \sum_{i=1}^2{\mathbb I}_{2i}\,\p_{x_i}\vec{\Phi}\rg)
\ee
Using again Poincar\'e Lemma on $D^2$ we obtain the existence of $\vec{V}$ such that
\be
\label{VI.8}
\nabla\vec{V}:=\vec{L}\times\nabla\vec{\Phi}+\, l_\sigma\,f'_\sigma(H)\,\nabla\vec{\Phi}-2\, (\al-\ov{\al})\, \nabla\vec{D}
\ee
Using the explicit expression of $\nabla\vec{D}$ given by (\ref{V.54}) we obtain
\be
\label{VI.9}
\vec{n}\cdot\nabla\vec{V}=\vec{n}\cdot(\vec{L}\times\nabla\vec{\Phi})=\vec{L}\cdot\nabla^\perp\vec{\Phi}=\nabla^\perp S+\nabla Y
\ee
We have also
\be
\label{VI.10}
\vec{n}\times\nabla\vec{V}=-\, (\vec{L}\cdot\vec{n})\,\nabla\vec{\Phi}-\,l_\sigma\,f_\sigma'(H)\,\nabla^\perp\vec{\Phi}-2\,(\al-\ov{\al})\,\vec{n}\times\nabla\vec{D} 
\ee
Denote $\pi_T(\nabla^\perp \vec{V})$ the tangential projection of $\nabla^\perp\vec{V}$ , we have
\be
\label{VI.11}
\ds\pi_T(\nabla^\perp \vec{V})= (\vec{L}\cdot\vec{n})\,\nabla\vec{\Phi}+l_\sigma\,f'_\sigma(H)\,\nabla^\perp\vec{\Phi}-2\, (\al-\ov{\al})\, \nabla^\perp\vec{D}
\ee
Hence
\be
\label{VI.12}
\vec{n}\times\nabla\vec{V}=-\nabla^\perp \vec{V}-\,2\, (\al-\ov{\al})\ \lf( \nabla^\perp\vec{D}+ \vec{n}\times\nabla\vec{D}\rg)-\vec{n}\,(\nabla S-\nabla^\perp Y)
\ee
Let $\vec{v}$ be the unique solution to
\be
\label{VI.13}
\lf\{
\begin{array}{l}
\ds\Delta\vec{v}=\nabla^\perp Y\cdot\nabla\vec{n}\quad\quad\mbox{ in } D^2\\[5mm]
\ds\vec{v}=0 \quad\quad\mbox{ on }\p D^2
\end{array}
\rg.
\ee
Using one more time Poincar\'e lemma we obtain the existence of $\vec{u}$ such that
\be
\label{VI.14}
\vec{n}\,\nabla^\perp{Y}=\nabla\vec{v}+\nabla^\perp \vec{u}
\ee
Finally, let $\vec{R}:=\vec{V}-\vec{u}$. We have
\be
\label{VI.15}
\vec{n}\times\nabla\vec{V}=\vec{n}\times\nabla\vec{R}+\vec{n}\times\nabla\vec{u}=\vec{n}\times\nabla\vec{R}+\vec{n}\times\nabla^\perp\vec{v}
\ee
Hence (\ref{VI.12}) becomes
\be
\label{VI.16}
\vec{n}\times\nabla\vec{R}+\vec{n}\times\nabla^\perp\vec{v}=-\nabla^\perp\vec{R}+\nabla\vec{v}-\vec{n}\,\nabla S-\,2\, (\al-\ov{\al})\ \lf( \nabla^\perp\vec{D}+ \vec{n}\times\nabla\vec{D}\rg)
\ee
which gives
\be
\label{VI.17}
\lf\{
\begin{array}{l}
\ds\Delta S=-\nabla\vec{n}\cdot\nabla^\perp{\vec{R}}+\mbox{div}\lf(\vec{n}\cdot\nabla\vec{v}-\,2\,(\al-\ov{\al})\ \vec{n}\cdot\nabla^\perp\vec{D}\rg)\\[5mm]
\ds\Delta\vec{R}=\nabla\vec{n}\times\nabla^\perp\vec{R}+\nabla\vec{n}\cdot\nabla^\perp S+\mbox{div}\lf(-\,\vec{n}\times\nabla\vec{v}+2\,(\al-\ov{\al})\   (-\nabla\vec{D}+  \vec{n}\times \nabla^\perp\vec{D}\rg)
\end{array}
\rg.
\ee
Taking the vectorial product between (\ref{VI.8}) and $\nabla^\perp\vec{\Phi}$ we obtain
\be
\label{VI.17-a}
\begin{array}{l}
\ds\nabla\vec{V}\times\nabla^\perp\vec{\Phi}=(\vec{L}\cdot\nabla^\perp\vec{\Phi})\cdot\nabla\vec{\Phi}-2\,l_\sigma\, f'_\sigma(H)\, e^{2\la}\,\vec{n}-2\, (\al-\ov{\al})\, \nabla\vec{D}\times\nabla^\perp\vec{\Phi}\\[5mm]
\ds\quad\quad\quad=\nabla^\perp S\cdot\nabla\vec{\Phi}+\nabla Y\cdot\nabla\vec{\Phi}-2\,l_\sigma\, f'_\sigma(H)\, e^{2\la}\,\vec{n}-2\, (\al-\ov{\al})\, \nabla\vec{D}\times\nabla^\perp\vec{\Phi}\quad.
\end{array}
\ee
We have also
\be
\label{VI.18}
\begin{array}{l}
\ds\nabla\vec{V}\times\nabla^\perp\vec{\Phi}=\nabla\vec{R}\times\nabla^\perp\vec{\Phi}+\nabla\vec{u}\times\nabla^\perp\vec{\Phi}\\[5mm]
\ds\quad=\nabla\vec{R}\times\nabla^\perp\vec{\Phi}+\nabla^\perp\vec{v}\times\nabla^\perp\vec{\Phi}+\nabla Y\cdot (\vec{n}\times\nabla^\perp\vec{\Phi})\\[5mm]
\ds\quad=\nabla\vec{R}\times\nabla^\perp\vec{\Phi}+\nabla\vec{v}\times\nabla\vec{\Phi}+\nabla Y\cdot\nabla\vec{\Phi}\quad.
\end{array}
\ee
Combining (\ref{VI.17-a}) and (\ref{VI.18}) gives
\be
\label{VI.19}
2\, l_\sigma\, f'_\sigma(H)\, e^{2\la}\,\vec{n}=\nabla^\perp S\cdot\nabla\vec{\Phi}-2\, (\al-\ov{\al})\, \nabla\vec{D}\times\nabla^\perp\vec{\Phi}-\nabla\vec{R}\times\nabla^\perp\vec{\Phi}-\nabla\vec{v}\times\nabla\vec{\Phi}
\ee
We have explicitly $l_\sigma\, f'_\sigma(H)=2\, H\,(1+ 2\,\sigma^2 (1+H^2))$ moreover, a straightforward computation gives
\be
\label{VI.19-1}
-\, \nabla\vec{D}\times\nabla^\perp\vec{\Phi}=\,2\, l_\sigma\,\vec{H}\, e^{2\la}=l_\sigma\,\Delta\vec{\Phi}
\ee
Inserting (\ref{VI.19-1}) in (\ref{VI.19}) gives then
\be
\label{VI.19-2}
2\,(1+2\, \sigma^2\, (1+H^2)-l_\sigma\,(\al-\ov{\al}))\ \Delta\vec{\Phi}=\nabla^\perp S\cdot\nabla\vec{\Phi}-\nabla\vec{R}\times\nabla^\perp\vec{\Phi}-\nabla\vec{v}\times\nabla\vec{\Phi}
\ee
\medskip

\noindent{\bf step 2} : We now prove that, for $\ep$ small enough, $\nabla S$ and $\nabla\vec{R}$ are in $L^2$.

Since $\int_{D^2}|\nabla\vec{n}|^2\ dx^2<\ep$, for $\ep$ small enough, following H\'elein's construction of energy controlled moving frame ( see \cite{Ri1}), we get the existence of $\vec{e}_i$ such that
\be
\label{VI.19-a}
\vec{e}_1\times\vec{e}_2=\vec{n}\quad\mbox{ and }\quad\int_{D^2}|\nabla\vec{e}_i|^2\ dx^2\le C\ \int_{D^2}|\nabla\vec{n}|^2\ dx^2
\ee
Using the assumptions (\ref{VI.0}) we have
\be
\label{VI.19-ab} 
\begin{array}{l}
 \ds 2\, C_1\ge W(\vec{\Phi})+\sigma^2\,\int_{S^2}(1+H^2)^2\\[5mm]
 \ds\quad+l_\sigma\,\lf[\int_{S^2} 2^{-1}\, |d\al|^2_{\gP}\ dvol_{\gP}+K_{g_0}\ \al\ e^{-2\al}-2^{-1}\, K_{g_0}\,\log A_{\vec{\Phi}}(\Sigma)\rg]
 \end{array}
\ee
using Onofri inequality (\ref{V.36-c}) we deduce that
\be
\label{VI.19-ac}
W(\vec{\Phi})+\sigma^2\,\int_{S^2}(1+H^2)^2\le 2\, C_1
\ee
Since $W(\vec{\Phi})\le 2\, C_1$ using (\ref{V.70}), we obtain the existence of a constant depending only on $C_1$ such that
\be
\label{VI.19-b}
\|\nabla\la\|_{L^{2,\infty}(D^2)}\le C
\ee
Recall (see for instance \cite{Ri1}) that for $\vec{e}_i:=e^{-\la}\,\p_{x_i}\vec{\Phi}$ the Liouville equation giving the expression of the Gauss curvature in conformal coordinates is equivalent to
\be
\label{VI.21}
-\Delta\la=(\nabla\vec{e}_1;\nabla^\perp\vec{e}_2)
\ee
Let $\nu$ be the solution of
\be
\label{VI.21-a}
\lf\{
\begin{array}{l}
\ds-\Delta\nu=(\nabla\vec{e}_1;\nabla^\perp\vec{e}_2)\quad\quad\mbox{ in }D^2\\[5mm]
\ds\nu=0\quad\quad\quad\mbox{ on }\p D^2
\end{array}
\rg.
\ee
We get, using Wente inequality
\be
\label{V.21-b}
\|\nabla\nu\|_{L^2(D^2)}+\|\nu\|_{L^\infty(D^2)}\le C\, \sum_{i=1}^2\int_{D^2}|\nabla\vec{e}_i|^2\ dx^2\le C\, \int_{D^2}|\nabla\vec{n}|^2\ dx^2\le \ep
\ee
Since $\la-\nu$ is harmonic and since
\be
\label{V.21-c}
\|\nabla(\la-\nu)\|_{L^{2,\infty}(D^2)}\le C
\ee
we have
\be
\label{V.21-d}
\|\nabla(\la-\nu)\|_{L^{\infty}(D^2_{1/2})}\le C
\ee
Hence, using (\ref{V.21-b}), there exists $\ov{\la}=(\la-\nu)(0)\in{\R}$ such that
\be
\label{VI.22}
\|\la-\ov{\la}\|_{L^\infty(D^2_{1/2})}\le C
\ee
Since $\|\nabla\mu\|_\infty\le \, C$, we have the existence of $\ov{\mu}\in {\R}$ such that
\be
\label{VI.23}
\|\mu-\ov{\mu}\|_{L^\infty(D^2)}<C\,
\ee
Hence we deduce the existence of $\ov{\al}\in {\R}$ such that
\be
\label{VI.24}
\|\al-\ov{\al}\|_{L^\infty(D^2_{1/2})}\le C
\ee
We rescale the domain so that $D^2_{1/2}$ becomes $D^2$. 
We now proceed to the introduction of $\vec{L}$, $Y$, $\vec{V}$, $\vec{v}$, $\vec{u}$ ,$\vec{R}$, $S$ as in step 1 but on $D^2_\rho$.
First of all using classical elliptic estimates we deduce from (\ref{VI.5}), using the hypothesis (\ref{VI.2}) and for $\sigma$ small enough
\be
\label{VI.38}
\|\nabla Y\|_{L^{2,\infty}(D^2_\rho)}\le C\ \int_{D^2_\rho} \sigma^2\, [1+H^4]\ e^{2\la}\ dx^2+C\, l_\sigma\, |\ov{\al}|\int_{D^2}e^{2\mu}\ dx^2+C\,l_\sigma\,\frac{A_{\vec{\Phi}}(D^2)}{A_{\vec{\Phi}}(\Sigma)}\le C\ \ep
\ee
and using (\ref{VI.13}) we deduce by the mean of Wente estimates 
\be
\label{VI.39}
\|\nabla\vec{v}\|_{L^2(D^2_\rho)}\le C\, \|\nabla\vec{n}\|_{L^2(D^2_\rho)}\ \lf[\int_{D^2_\rho} \sigma^2\, H^4\ e^{2\la}\ dx^2+\, l_\sigma\ |\ov{\al}|\int_{D^2}e^{2\mu}\ dx^2+l_\sigma\rg]\quad.
\ee
We have also taking $\int_{D^2_\rho}\vec{L}\ dx^2=0$
\be
\label{VI.40}
\begin{array}{l}
\ds\|e^{\ov{\la}}\,\vec{L}\|_{L^{2,\infty}+\sigma^{-1/2}L^{4/3}(D^2_\rho)}\le C\, \|e^\la\, H\|_{L^{2}(D^2_\rho)}+C\,l_\sigma\,\|\nabla\al\|^2_{L^2(D^2_\rho)}\\[5mm]
\ds+C\, l_\sigma\, |\ov{\al}|\int_{D^2}e^{2\mu}\ dx^2+\|\nabla\vec{n}\|^2_{L^{2}(D^2_\rho)}+\lf[ \int_{D^2_\rho}\sigma^2 |\nabla\vec{n}|^4 e^{-2\la}  \rg]^{3/4}+C\ l_\sigma\,\frac{A_{\vec{\Phi}}(D^2)}{A_{\vec{\Phi}}(\Sigma)}
\end{array}
\ee
Using the assumptions of the lemma we have then
\be
\label{VI.40-a}
\begin{array}{l}
\ds\|e^{\ov{\la}}\,\vec{L}\|_{L^{2,\infty}+\sigma^{-1/2}L^{4/3}(D^2_\rho)}\le C\ \|\nabla\vec{n}\|_{L^2(D^2_\rho)}+C\,l_\sigma\,\|\nabla\al\|^2_{L^2(D^2_\rho)}\\[5mm]
\ds\quad\quad\quad+C\, l_\sigma\, |\ov{\al}|\int_{D^2}e^{2\mu}\ dx^2+\lf[\int_{D^2_\rho}\sigma^2 |\nabla\vec{n}|^4 e^{-2\la}\ dx^2\rg]^{3/4}+ C\ l_\sigma\,\frac{A_{\vec{\Phi}}(D^2)}{A_{\vec{\Phi}}(\Sigma)}
\end{array}
\ee
Hence, using the explicit expressions (\ref{VI.6}), (\ref{VI.8}) and (\ref{VI.14}) we obtain that $\nabla\vec{R}$ and $\nabla\vec{S}$ are uniformly bounded in $L^{2,\infty}+\sigma^{-1/2}L^{4/3}$ and we have respectively
\be
\label{VI.41}
\begin{array}{l}
\|\nabla\vec{R}\|+\|\nabla S\|_{L^{2,\infty}+\sigma^{-1/2}L^{4/3}(D^2_\rho)}\le C\, \|\nabla\vec{n}\|_{L^{2}(D^2_\rho)}+C\,l_\sigma\,\|\nabla\al\|^2_{L^2(D^2_\rho)}\\[5mm]
\ds\quad\quad\quad+C\, l_\sigma\, |\ov{\al}|\int_{D^2}e^{2\mu}\ dx^2+\lf[\int_{D^2_\rho}\sigma^2 |\nabla\vec{n}|^4 e^{-2\la}\ dx^2\rg]^{3/4}+C\ l_\sigma\,\frac{A_{\vec{\Phi}}(D^2)}{A_{\vec{\Phi}}(\Sigma)}
\end{array}
\ee
Let $\varphi$ and $\vec{\Psi}$ the unique solutions in $W^{1,2}_0(D^2_\rho)$ of the linear system
\be
\label{VI.42}
\lf\{
\begin{array}{l}
\ds\Delta \varphi=-\nabla\vec{n}\cdot\nabla^\perp{\vec{R}}+\mbox{div}\lf(\vec{n}\cdot\nabla\vec{v}-\,2\,(\al-\ov{\al})\ \vec{n}\cdot\nabla^\perp\vec{D}\rg)\\[5mm]
\ds\Delta\vec{\Psi}=\nabla\vec{n}\times\nabla^\perp\vec{R}+\nabla\vec{n}\cdot\nabla^\perp S+\mbox{div}\lf(-\,\vec{n}\times\nabla\vec{v}+2\,(\al-\ov{\al})\   (-\nabla\vec{D}+  \vec{n}\times \nabla^\perp\vec{D}\rg)
\end{array}
\rg.
\ee
Using  lemma~\ref{lm-wente} and in particular (\ref{A-2}), together with the fact that both $\nabla\vec{v}$ and $(\al-\ov{\al})\,\nabla\vec{D}$ are in $L^2$ we deduce using $l_\sigma\ \|\al-\ov{\al}\|_{L^\infty(D^2_\rho)}<\ep$
\be
\label{VI.45}
\begin{array}{l}
\ds\|\nabla\varphi\|_{L^{2,\infty}(D^2_\rho)}+\|\nabla\vec{\Psi}\|_{L^{2,\infty}(D^2_\rho)}\le \ep^{1/2}\ \|\nabla\vec{n}\|_{L^{2}(D^2_\rho)}+\,C\,\ep^{1/2}\,l_\sigma\,\|\nabla\al\|^2_{L^2(D^2_\rho)}\\[5mm]
\ds+\,C\, \ep^{1/2}\,l_\sigma\, |\ov{\al}|\int_{D^2}e^{2\mu}\ dx^2+\ep^{1/2}\ \lf[\int_{D^2_\rho}\sigma^2 |\nabla\vec{n}|^4 e^{-2\la}\ dx^2\rg]^{3/4}+C\, l_\sigma\,\ep^{1/2}\,\frac{A_{\vec{\Phi}}(D^2)}{A_{\vec{\Phi}}(\Sigma)}
\end{array}
\ee
Since $S-\varphi$ and $\vec{R}-\vec{\Psi}$ are harmonic, we have
\be
\label{VI.46}
\begin{array}{l}
\ds\|\nabla(S-\varphi)\|+\|\nabla(\vec{R}-\vec{\Psi})\|_{L^{2,\infty}(D^2_{\rho/2})}\le C\, \|\nabla(S-\varphi)\|+\|\nabla(\vec{R}-\vec{\Psi})\|_{L^{4/3}(D^2_\rho)}\\[5mm]
\ds\quad\le C\,\|\nabla S\|+C\,\|\nabla\vec{R}\|_{L^{2,\infty}+\sigma^{-1/2}L^{4/3}(D^2_\rho)}+C\,\|\nabla\varphi\|_{L^{2,\infty}(D^2_\rho)}+C\,\|\nabla\vec{\Psi}\|_{L^{2,\infty}(D^2_\rho)}
\end{array}
\ee
Hence
\be
\label{VI.47}
\begin{array}{l}
\|\nabla S\|+\|\nabla\vec{R}\|_{L^{2,\infty}(D^2_{\rho/2})}\le\|\nabla S\|+\|\nabla\vec{R}\|_{L^{2,\infty}+\sigma^{-1/2}L^{4/3}(D^2_\rho)}\\[5mm]
\quad\quad+2\,\|\nabla\varphi\|_{L^{2,\infty}(D^2_\rho)}+2\,\|\nabla\vec{\Psi}\|_{L^{2,\infty}(D^2_\rho)}\\[5mm]
\ds\quad\quad\le C\, \|\nabla\vec{n}\|_{L^{2}(D^2_\rho)}+C\,\sigma^2\,\|\nabla\al\|^2_{L^2(D^2_\rho)}\\[5mm]
\ds\quad\quad+C\, l_\sigma\ |\ov{\al}|\int_{D^2}e^{2\mu}\ dx^2+\lf[\int_{D^2_\rho}\sigma^2 |\nabla\vec{n}|^4 e^{-2\la}\ dx^2\rg]^{3/4}+C\, l_\sigma\,\frac{A_{\vec{\Phi}}(D^2)}{A_{\vec{\Phi}}(\Sigma)}

\end{array}
\ee
Let $\varphi_1$ and $\vec{\Psi}_1$ the unique solutions in $W^{1,2}_0(D^2_{\rho/2})$ of the linear system
\be
\label{VI.48}
\lf\{
\begin{array}{l}
\ds\Delta \varphi_1=-\nabla\vec{n}\cdot\nabla^\perp{\vec{R}}+\mbox{div}\lf(\vec{n}\cdot\nabla\vec{v}-\,2\,(\al-\ov{\al})\ \vec{n}\cdot\nabla^\perp\vec{D}\rg)\\[5mm]
\ds\Delta\vec{\Psi}_1=\nabla\vec{n}\times\nabla^\perp\vec{R}+\nabla\vec{n}\cdot\nabla^\perp S+\mbox{div}\lf(-\,\vec{n}\times\nabla\vec{v}+2\,(\al-\ov{\al})\   (-\nabla\vec{D}+  \vec{n}\times \nabla^\perp\vec{D}\rg)
\end{array}
\rg.
\ee
Wente estimates, combined with classical elliptic estimates give, using (\ref{VI.39})
\be
\label{VI.49}
\begin{array}{l}
\ds\|\nabla\varphi_1\|_{L^2(D^2_{\rho/2})}+\|\nabla\vec{\Psi}_1\|_{L^2(D^2_{\rho/2})}\le\ C\, \|\nabla\vec{n}\|_{L^{2}(D^2_{\rho/2})}\ \lf[   \|\nabla S\|_{L^{2,\infty}(D^2_{\rho/2})}+\|\nabla\vec{R}\|_{L^{2,\infty}(D^2_{\rho/2})} \rg]\\[5mm]
\ds +\,C\, \|\nabla\vec{n}\|_{L^2(D^2_\rho)}\ \lf[\int_{D^2_\rho} [|\nabla\vec{n}|^2+\sigma^2|\nabla\vec{n}|^4]\ dx^2+C\,l_\sigma\, |\ov{\al}|\int_{D^2}e^{2\mu}\ dx^2+C\, l_\sigma\,\frac{A_{\vec{\Phi}}(D^2)}{A_{\vec{\Phi}}(\Sigma)}
\rg]\\[5mm]
\ds + \,C\, l_\sigma\ \|\al-\ov{\al}\|_{L^\infty(D^2_\rho)}\ \|\nabla\vec{n}\|^2_{L^{2}(D^2_\rho)}
\end{array}
\ee
Since $S-\varphi$ and $\vec{R}-\vec{\Psi}$ are harmonic we obtain finally
\be
\label{VI.50}
\begin{array}{l}
\|\nabla S\|_{L^{2}(D^2_{\rho/4})}+\|\nabla\vec{R}\|_{L^{2}(D^2_{\rho/4})}\le\|\nabla S\|_{L^{2,\infty}(D^2_{\rho/2})}+\|\nabla\vec{R}\|_{L^{2,\infty}(D^2_{\rho/2})}\\[5mm]
\quad\quad+2\,\|\nabla\varphi_1\|_{L^{2}(D^2_\rho/2)}+2\,\|\nabla\vec{\Psi}_1\|_{L^{2}(D^2_{\rho/2})}\\[5mm]
\ds\le C\, \|\nabla\vec{n}\|_{L^{2}(D^2_\rho)}+C\,l_\sigma\,\|\nabla\al\|^2_{L^2(D^2_\rho)}+C\, l_\sigma\, |\ov{\al}|\int_{D^2}e^{2\mu}\ dx^2\\[5mm]
\ds\quad+\, C\,\lf[\int_{D^2_\rho}\sigma^2 |\nabla\vec{n}|^4 e^{-2\la}\ dx^2\rg]^{3/4}+C\, l_\sigma\,\frac{A_{\vec{\Phi}}(D^2)}{A_{\vec{\Phi}}(\Sigma)}
\end{array}
\ee
Combining lemma~\ref{lm-L-4} with (\ref{VI.0-a}) gives
\be
\label{VI.50-a}
\sigma^2\int_{D^2_{1/2}}|\nabla\vec{n}|^4\ e^{-2\la}\, dx^2\le C\ \sigma^2\,\int_{D^2}H^4\,e^{2\la}\ dx^2+C\, \lf[\int_{D^2}|\nabla\vec{n}|^2\ dx^2\rg]^2
\ee
Combining (\ref{VI.50}) and (\ref{VI.50-a}) gives then
\be
\label{VI.50-b}
\begin{array}{l}
\|\nabla S\|_{L^{2}(D^2_{1/8})}+\|\nabla\vec{R}\|_{L^{2}(D^2_{1/8})}\le C\, \|\nabla\vec{n}\|_{L^{2}(D^2)}+C\,l_\sigma\,\|\nabla\al\|^2_{L^2(D^2)}\\[5mm]
\ds\quad\quad+C\, l_\sigma\, |\ov{\al}|\int_{D^2}e^{2\mu}\ dx^2+C\ \sigma^2\,\int_{D^2}H^4\,e^{2\la}\ dx^2+C\, l_\sigma\,\frac{A_{\vec{\Phi}}(D^2)}{A_{\vec{\Phi}}(\Sigma)}
\end{array}
\ee

Hence $\nabla S$ and $\nabla\vec{R}$ are in $L^2(D^2_{1/8})$ and under the assumptions of the lemma, using also lemma~\ref{lm-L-infty-borne}, we have that $\|\nabla S\|_{L^2(D^2_{1/8})}+\|\nabla\vec{R}\|_{L^2(D^2_{1/8})}$
are bounded by a constant depending only on $C_0$. We rescale the domain in such a way that $\nabla S$ and $\nabla\vec{R}$ are in $L^2(D^2)$.

\medskip

\noindent{\bf Step 3.} Uniform Morrey decrease of the Willmore energy. Precisely we are going to prove the existence of $\gamma>0$ independent of $\sigma$ and independent
of the solution such that
\be
\label{VI.51}
\sup_{x_0\in D^2_{1/2}\ ; \ r<1/4}\quad r^{-\gamma}\int_{B^2_r(x_0)}H^2\, (1+\sigma^2 (1+H^2))^2\ e^{2\la}\ dx^2\le C
\ee
Following step 1 of the proof of the theorem the map $\vec{U}:=(e^{-\ov{\la}}\vec{\Phi},\vec{v},S,\vec{R})$ satisfies the following system on $B_r(x_0)$
\be
\label{VI.54}
\lf\{
\begin{array}{l}
\ds \Delta(e^{-\ov{\la}}\,\vec{\Phi})=h_\sigma\ e^{-\ov{\la}}\,\lf[\nabla^\perp S\cdot\nabla\vec{\Phi}-\nabla\vec{R}\times\nabla^\perp\vec{\Phi}-\nabla\vec{v}\times\nabla\vec{\Phi}\rg]\\[5mm]
\ds\Delta\vec{v}=\nabla^\perp Y\cdot\nabla\vec{n}\\[5mm]
\ds\Delta S=-\nabla\vec{n}\cdot\nabla^\perp{\vec{R}}+\mbox{div}\lf(\vec{n}\cdot\nabla\vec{v}-\,2\,(\al-\ov{\al})\ \vec{n}\cdot\nabla^\perp\vec{D}\rg)\\[5mm]
\ds\Delta\vec{R}=\nabla\vec{n}\times\nabla^\perp\vec{R}+\nabla\vec{n}\cdot\nabla^\perp S+\mbox{div}\lf(-\,\vec{n}\times\nabla\vec{v}+2\,(\al-\ov{\al})\   (-\nabla\vec{D}+  \vec{n}\times \nabla^\perp\vec{D}\rg)
\end{array}
\rg.
\ee
where $0\le h_\sigma:=2\,(1+2\, \sigma^2\, (1+H^2)-l_\sigma\,(\al-\ov{\al}))^{-1}\le1$ and where we shall use the fact that
\be
\label{VI.55}
\p_{x_i}\vec{n}=-\,e^{-2 (\la-\ov{\la})}\, \sum_{j=1}^2\vec{n}\cdot\p^2_{x_i x_j}\lf(e^{-\ov{\la}}\,\vec{\Phi}\rg)\ e^{-\ov{\la}}\p_{x_j}\vec{\Phi}
\ee
and that
\be
\label{VI.56}
\nabla\vec{D}=l_\sigma\ e^{-2(\la-\ov{\la})}
 \lf(\sum_{i=1}^2\vec{n}\cdot\p^2_{x_1x_i}\lf(e^{-\ov{\la}}\vec{\Phi}\rg)\,e^{-\ov{\la}}\p_{x_i}\vec{\Phi}\ ,\, \vec{n}\cdot\p^2_{x_2x_i}\lf(e^{-\ov{\la}}\vec{\Phi}\rg)\,e^{-\ov{\la}}\p_{x_i}\vec{\Phi}\rg)
\ee
Let $\vec{w}$ in $W^{1,2}_0(B_r(x_0))$ the solutions of
\be
\label{VI.56-a}
\ds\Delta\vec{w}=\nabla^\perp Y\cdot\nabla\vec{n}
\ee
Using Wente estimates we obtain
\be
\label{VI.56-b}
\ds\int_{B_r(x_0)}|\nabla\vec{w}|^2\ dx^2\le C\ \|\nabla Y\|^2_{L^{2,\infty}(D^2)}\ \int_{B_r(x_0)}|\nabla\vec{n}|^2\ dx^2
\ee
Since $\vec{v}-\vec{w}$ is harmonic for any $t\in (0,1)$ the monotonicity formula for harmonic functions give
\be
\label{VI.56-c}
\int_{B_{t\,r}(x_0)}|\nabla(\vec{v}-\vec{w})|^2\ dx^2\le t^2\ \int_{B_{r}(x_0)}|\nabla(\vec{v}-\vec{w})|^2\ dx^2
\ee
We deduce from (\ref{VI.38}), (\ref{VI.56-b}) and (\ref{VI.56-c})
\be
\label{VI.56-d}
\begin{array}{l}
\ds\int_{B_{t\,r}(x_0)}|\nabla\vec{v}|^2\ dx^2\le t^2\ \int_{B_{r}(x_0)}|\nabla\vec{v}|^2\ dx^2+C\ \sqrt{\ep}\ \frac{\|\nabla Y\|^2_{L^{2,\infty}(D^2)}}{\sqrt{\ep}} \int_{B_r(x_0)}|\nabla\vec{n}|^2\ dx^2
\end{array}
\ee
Let $T$ and $\vec{Q}$ in $W^{1,2}_0(B_{t\,r}(x_0))$ solving
\be
\label{VI.56-e}
\lf\{
\begin{array}{l}
\ds\Delta T=-\nabla\vec{n}\cdot\nabla^\perp{\vec{R}}+\mbox{div}\lf(\vec{n}\cdot\nabla\vec{v}-\,2\,(\al-\ov{\al})\ \vec{n}\cdot\nabla^\perp\vec{D}\rg)\\[5mm]
\ds\Delta\vec{Q}=\nabla\vec{n}\times\nabla^\perp\vec{R}+\nabla\vec{n}\cdot\nabla^\perp S+\mbox{div}\lf(-\,\vec{n}\times\nabla\vec{v}+2\,(\al-\ov{\al})\   (-\nabla\vec{D}+  \vec{n}\times \nabla^\perp\vec{D}\rg)
\end{array}
\rg.
\ee
Wente inequalities combined with classical elliptic estimates and the fact that $$\int_{D^2}|\nabla\vec{n}|^2\ dx^2+l^2_\sigma\,\|\al-\ov{\al}\|^2_{L^\infty(D^2)}<\ep$$ gives
\be
\label{VI.56-f}
\begin{array}{l}
\ds\int_{B_{t\,r}(x_0)}[|\nabla T|^2+|\nabla\vec{Q}|^2]\ dx^2\le C\ \sqrt{\ep}\ \lf[\int_{B_{t\,r}(x_0)}[|\nabla S|^2+|\nabla\vec{R}|^2+\delta^2\,|\nabla\vec{n}|^2\rg]\ dx^2\\[5mm]
\ds\quad\quad\quad\quad+\, C\, \int_{B_{t\,r}(x_0)}|\nabla\vec{v}|^2\ dx^2
\end{array}
\ee
where 
\[
\delta^2:=\frac{l^2_\sigma\,\|\al-\ov{\al}\|^2_{L^\infty(D^2)}+\|\nabla Y\|^2_{L^{2,\infty}(D^2)}}{\sqrt{\ep}}
\]
Since $S-T$ and $\vec{R}-\vec{Q}$ are harmonic, the monotonicity formula gives
\be
\label{VI.56-g}
\begin{array}{l}
\ds\int_{B_{t^2\,r}(x_0)}\lf[|\nabla(S-T)|^2+|\nabla(\vec{R}-\vec{Q})|^2\rg]\ dx^2\\[5mm]
\ds\quad\quad\le t^2\ \int_{B_{t\,r}(x_0)}\lf[|\nabla(S-T)|^2+|\nabla(\vec{R}-\vec{Q})|^2\rg]\ dx^2
\end{array}
\ee
Hence combining (\ref{VI.56-d}), (\ref{VI.56-f}) and (\ref{VI.56-g}) we obtain
\be
\label{VI.61}
\begin{array}{l}
\ds\int_{B_{t^2\,r}(x_0)}\lf[|\nabla S|^2+|\nabla\vec{R}|^2+|\nabla\vec{v}|^2\rg]\ dx^2\le \ t^2 \int_{B_{r}(x_0)}\lf[|\nabla S|^2+|\nabla\vec{R}|^2+|\nabla\vec{v}|^2\rg]\\[5mm]
\ds\quad\quad\quad+\,C\, \sqrt{\ep}\ \int_{B_{r}(x_0)}\lf[|\nabla S|^2+|\nabla\vec{R}|^2+\delta^2\,|\nabla\vec{n}|^2\rg]\ dx^2
\end{array}
\ee
We recall the structural equation (see \cite{Ri1})
\be
\label{VI.64-e}
\nabla\vec{n}=\nabla^\perp\vec{n}\times\vec{n}-2\, H\, \nabla\vec{\Phi}\quad.
\ee
Taking the divergence gives then
\be
\label{VI.64-f}
\Delta\vec{n}=\nabla^\perp\vec{n}\times\nabla\vec{n}-2\, \mbox{div}\lf[\,2\, H\, \nabla\vec{\Phi}\rg]
\ee
We introduce $\vec{\xi}$ to be the solution of
\be
\label{VI.57}
\lf\{
\begin{array}{l}
\ds \Delta\vec{\xi}=\nabla^\perp\vec{n}\times\nabla\vec{n}-2\, \mbox{div}\lf[\,2\, H\, \nabla\vec{\Phi}\rg]\mbox{ in } B_{r}(x_0)\\[5mm] 
\ds\vec{\xi}=0\quad\quad\quad\mbox{ on }\p B_{r}(x_0)
\end{array}
\rg.
\ee
Classical elliptic estimates combined with  the first equation of (\ref{VI.54}) gives
\be
\label{VI.58}
\int_{B_{r}(x_0)}|\nabla\vec{\xi}|^2\ dx^2\le C\, \int_{B_{r}(x_0)}\lf[|\nabla S|^2+|\nabla\vec{R}|^2+|\nabla\vec{v}|^2\rg]\ dx^2\quad.
\ee
Since $\vec{n}-\vec{\xi}$ is harmonic on $B_{r}(x_0)$ we have 
\be
\label{VI.59}
\int_{B_{t\,r}(x_0)}|\nabla(\vec{n}-\vec{\xi})|^2\ dx^2\le t^2\ \int_{B_{r}(x_0)}|\nabla(\vec{n}-\vec{\xi})|^2\ dx^2\
\ee
Hence we have
\be
\label{VI.59-a}
\begin{array}{l}
\ds\int_{B_{t\,r}(x_0)}|\nabla\vec{n}|^2\ dx^2\le t^2\ \int_{B_{r}(x_0)}|\nabla\vec{n}|^2\ dx^2\\[5mm]
\ds\quad\quad\quad+\,C\, \int_{B_{r}(x_0)}\lf[|\nabla S|^2+|\nabla\vec{R}|^2+|\nabla\vec{v}|^2\rg]\ dx^2
\end{array}
\ee
Inserting (\ref{VI.61}) in (\ref{VI.59-a}) we finally obtain, taking $\sigma$ small enough in such a way that $\delta<1$
\be
\label{VI.63}
\begin{array}{l}
\ds\int_{B_{t^2\,r}(x_0)}\lf[|\nabla S|^2+|\nabla\vec{R}|^2+|\nabla\vec{v}|^2+\delta^2\, |\nabla\vec{n}|^2\rg]\ dx^2\\[5mm]
\ds\quad\quad\le \ C\  [t^2+C\,\sqrt{\ep}+\delta^2] \int_{B_{r}(x_0)}\lf[|\nabla S|^2+|\nabla\vec{R}|^2+|\nabla\vec{v}|^2+\delta^2\,|\nabla\vec{n}|^2\rg]\ dx^2
\end{array}
\ee
We fix $t>0$ and $\ep>0$ independent of $r$, $x_0$, $\sigma<\sigma_0$ and the solution such that $C\ [t^2+C\,\sqrt{\ep}] \le 1/2$. By classical iteration argument we deduce the existence $\gamma>0$
such that
\be
\label{VI.64}
\begin{array}{l}
\ds\sup_{x_0\in D^2_{1/8}\ ; \ r<1/16}\quad r^{-\gamma}\int_{B^2_r(x_0)}\lf[|\nabla S|^2+|\nabla\vec{R}|^2+|\nabla\vec{v}|^2+\delta^2\,|\nabla^2(e^{-\ov{\la}}\vec{\Phi})|^2\rg]\ dx^2\\[5mm]
\ds\quad\quad\le \int_{D^2_{1/4}}\lf[|\nabla S|^2+|\nabla\vec{R}|^2+|\nabla\vec{v}|^2+\delta^2\,|\nabla\vec{n}|^2\rg]\ dx^2
\end{array}
\ee
Combining (\ref{VI.38}), (\ref{VI.50-b}) and (\ref{VI.64}) we obtain in particular
\be
\label{VI.64-c}
\begin{array}{l}
\ds\sup_{x_0\in D^2_{1/8}\ ; \ r<1/16}\quad r^{-\gamma/2}\lf[\int_{B^2_r(x_0)}\lf[|\nabla S|^2+|\nabla\vec{R}|^2+|\nabla\vec{v}|^2\rg]\ dx^2\rg]^{1/2}\\[5mm]
\ds\quad\quad\le   C\, \|\nabla\vec{n}\|_{L^{2}(D^2)}+C\,l_\sigma\,\|\nabla\al\|^2_{L^2(D^2)}\\[5mm]
\ds\quad\quad+C\, l_\sigma\, |\ov{\al}|\int_{D^2}e^{2\mu}\ dx^2+C\ \sigma^2\,\int_{D^2}H^4\,e^{2\la}\ dx^2+C\, l_\sigma\,\frac{A_{\vec{\Phi}}(D^2)}{A_{\vec{\Phi}}(\Sigma)}
\end{array}
\ee
Combining (\ref{VI.64-c}) and (\ref{VI.19-2}) gives
\be
\label{VI.64-d}
\begin{array}{l}
\ds\sup_{x_0\in D^2_{1/8}\ ; \ r<1/16}\quad r^{-\gamma}\,\int_{B^2_r(x_0)}H^2\lf[1+\sigma^2\, (1+H^2)]^2\ dx^2\rg]\ e^{2\la} \ dx^2\\[5mm]
\ds\quad\quad\le   C\, \int_{D^2}|\nabla\vec{n}|^2\ dx^2+C\,\lf[l_\sigma\,\int_{D^2}|\nabla\al|^2\ dx^2 \rg]^2\,\\[5mm]
\ds\quad\quad+C\, l^2_\sigma\, |\ov{\al}|^2\lf[\int_{D^2}e^{2\mu}\ dx^2\rg]^2+C\ \lf[\sigma^2\,\int_{D^2}H^4\,e^{2\la}\ dx^2\rg]^2+C\, l^2_\sigma\,\lf[\frac{A_{\vec{\Phi}}(D^2)}{A_{\vec{\Phi}}(\Sigma)}\rg]^2
\end{array}
\ee

\medskip

\noindent{\bf Step 4.} Bootstraping (\ref{VI.64-d}). Lemma~\ref{lm-L-5} applied to (\ref{VI.64-f}) implies  that
\be
\label{VI.64-g}
\begin{array}{l}
\ds\sup_{x_0\in D^2_{1/16}\ ; \ r<1/32}\quad r^{-\gamma}\,\int_{B^2_r(x_0)}|\nabla\vec{n}|^2 \ dx^2\\[5mm]
\ds\quad\quad\le   C\, \int_{D^2}|\nabla\vec{n}|^2\ dx^2+C\,\lf[l_\sigma\,\int_{D^2}|\nabla\al|^2\ dx^2 \rg]^2\, \\[5mm]
\ds\quad\quad+C\, l^2_\sigma\, |\ov{\al}|^2\lf[\int_{D^2}e^{2\mu}\ dx^2\rg]^2+C\ \lf[\sigma^2\,\int_{D^2}H^4\,e^{2\la}\ dx^2\rg]^2+C\, l^2_\sigma\,\lf[\frac{A_{\vec{\Phi}}(D^2)}{A_{\vec{\Phi}}(\Sigma)}\rg]^2
\end{array}
\ee
The Liouville equation reads as follows
\be
\label{VI.64-h}
-\Delta \al=e^{2\la}\ K\, +\Delta\mu= e^{2\la}\ K- e^{2\mu}\ K_{g_0}
\ee
Thus
\be
\label{VI.64-i}
\int_{B^2_r(x_0)}|\Delta\al|\le \int_{B^2_r(x_0)}|\nabla\vec{n}|^2 \ dx^2+4\,\pi^2\, r^2\, \|e^{2\mu}\|_{L^\infty(D^2)}
\ee
Combining (\ref{VI.64-g}), (\ref{VI.64-i}) and Adams-Morrey embedding gives that 
\be
\label{VI.64-j}
\begin{array}{l}
\forall\ p<\frac{2-\gamma}{1-\gamma}\quad\quad\|\nabla\al\|_{L^p(D^2_{1/32})}\\[5mm]
\ds\le\ C_p\, \sup_{x_0\in D^2_{1/16}\ ; \ r<1/32}\quad r^{-\gamma}\,\int_{B^2_r(x_0)}|\nabla\vec{n}|^2+C\, \|e^{2\mu}\|_{L^\infty(D^2)}+\|\nabla\al\|_{L^2(D^2)}
\end{array}
\ee
This gives in particular
\be
\label{VI.64-k}
\begin{array}{l}
\ds\lf[  \int_{D^2_{1/32}} [l_\sigma\ |\nabla\al|^{2} ]^{p/2} \rg]^{2/p}\le \ C\, l_\sigma\, \lf[\sup_{x_0\in D^2_{1/16}\ ; \ r<1/32}\quad r^{-\gamma}\,\int_{B^2_r(x_0)}|\nabla\vec{n}|^2\ dx^2\rg]^2\\[5mm]
\ds\quad\quad+C\, l_\sigma\,\|e^{2\mu}\|^2_{L^\infty(D^2)}+\, C\, l_\sigma\int_{D^2}|\nabla\al|^2\, dx^2
\end{array}
\ee
The equation (\ref{VI.3a}) gives
\be
\label{VI.65}
\begin{array}{l}
\ds\Delta \lf(2\,\vec{H} (1+2\, \sigma^2\, (1+H^2))\rg)=\mbox{div}\lf(2\, l_\sigma f'_\sigma(H)\,\nabla\vec{n}+2\ e^{-2\la}\, l_\sigma\, \nabla\vec{\Phi}\cdot\nabla\al\ \nabla\al\rg.\\[5mm]
\ds\quad\lf.-\, l_\sigma\ [-2\, f(H)+e^{-2\la}\,|\nabla\al|^2-K_{g_0}\, \al\ e^{-2\al}+K_{g_0}\,A_{\vec{\Phi}}(\Sigma)^{-1}]\ \nabla\vec{\Phi}-2\, l_\sigma\ e^{-2\la}\, (\vec{\mathbb I}\res\nabla^\perp\al)^\perp\rg)
\end{array}
\ee
This gives that $\vec{V}:=e^{\ov{\la}}\,2\,\vec{H} (1+2\,\sigma^2\, (1+H^2))$ satisfies an equation of the form
\be
\label{VI.66}
\Delta\vec{V}=\mbox{div}(\vec{I}+\vec{J})
\ee
where 
\be
\label{VI.67}
\begin{array}{l}
\ds\sup_{x_0\in D^2_{1/32}\ ; \ r<1/4}\quad r^{-\gamma}\int_{B^2_r(x_0)}|\vec{I}|\ dx^2\le C\\[5mm]
\ds\quad\quad\le   C\, \int_{D^2}|\nabla\vec{n}|^2\ dx^2+C\,\lf[l_\sigma\,\int_{D^2}|\nabla\al|^2\ dx^2 \rg]^2\,\\[5mm]
\ds\quad\quad+C\, l^2_\sigma\, |\ov{\al}|^2\ \|e^{4\mu}\|_{L^\infty(D^2)}+C\ \lf[\sigma^2\,\int_{D^2}H^4\,e^{2\la}\ dx^2\rg]^2+C\, l^2_\sigma\,\lf[\frac{A_{\vec{\Phi}}(D^2)}{A_{\vec{\Phi}}(\Sigma)}\rg]^2
\end{array}
\ee
where $\gamma$ and $C$ are independent of the solution and of $\sigma$ but depend only on the constant $C_0$ in the statement of the lemma moreover, using the fact that $2-\gamma/(1-\gamma)>2$, using (\ref{VI.64-k})
there exists $q>1$ such that
\be
\label{VI.67-a}
\begin{array}{l}
\ds\lf[\int_{D^2_{1/32}}\,|\vec{J}|^q\ dx^2\rg]^{1/q}\le\, C\, l_\sigma\, \lf[\sup_{x_0\in D^2_{1/16}\ ; \ r<1/32}\quad r^{-\gamma}\,\int_{B^2_r(x_0)}|\nabla\vec{n}|^2\ dx^2\rg]^2\\[5mm]
\ds\quad\quad\quad+C\, [l_\sigma\,\|e^{2\mu}\|_{L^\infty(D^2)}+l_\sigma\,|\ov{\al}|]\,\|e^{2\mu}\|_{L^\infty(D^2)}+\, C\, l_\sigma\int_{D^2}|\nabla\al|^2\, dx^2
\end{array}
\ee
Using (\ref{V.67-a1}) and classical Adams-Sobolev
inequalities (see \cite{Ri1}) give then the existence of $p>2$ such that
\be
\label{VI.67-b}
\begin{array}{l}
2\,\|e^{\ov{\la}}\,\vec{H} (1+2\, \sigma^2\, (1+H^2))\|^2_{L^p(D_{1/64})}=\|V\|^2_{L^p(D_{1/64})}\\[5mm]
\ds\quad\quad\le\,  C\, \int_{D^2}|\nabla\vec{n}|^2\ dx^2+C\,\lf[l_\sigma\,\int_{D^2}|\nabla\al|^2\ dx^2 \rg]^2+C\ \lf[\sigma^2\,\int_{D^2}H^4\,e^{2\la}\ dx^2\rg]^2 \\[5mm]
\ds\quad\quad\quad+C\, [l^2_\sigma\,|\ov{\al}|^2\,+ l^2_\sigma\,\|e^{4\mu}\|_{L^\infty(D^2)}]\,\|e^{4\mu}\|_{L^\infty(D^2)}+C\, l^2_\sigma\,\lf[\frac{A_{\vec{\Phi}}(D^2)}{A_{\vec{\Phi}}(\Sigma)}\rg]^2
\end{array}
\ee
where we have used (\ref{V.69-a}), (\ref{V.69-b}) .

\medskip

Bootstraping this information respectively in   the three elliptic systems (\ref{VI.64-f}), (\ref{VI.64-h}) and (\ref{VI.65}) (which are now becoming sub-critical for $\vec{V}\in L^p$ with $p>2$)
one obtains (\ref{VI.3}), (\ref{VI.3-a}) and lemma~\ref{lm-ep-reg} is proved.\hfill $\Box$

\section{The Palais-Smale Condition for the Frame Energies}
\reset
\subsection{Sequential weak compactness of weak Immersions in ${\mathcal E}_{\Sigma,2}$ with uniformly bounded Frame Energies.}

In this section we are working with the lagrangian $F^\sigma$ defined in the previous section but the parameter $\sigma$ will be \underbar{fixed} all along
the section. So, in order to simplify the presentation, we will simply work with the following lagrangian
\be
\label{VII.00}
F(\vec{\Phi}):=\int_{\Sigma}\lf[H^2+ \, \lf[ (1+H^2)^2+2^{-1}\,|d\al|_{\gP}^2+\, K_{g_0}\,\al\, e^{-2\al}\rg]\rg]\ dvol_{\gP} -2^{-1}\, K_{g_0}\,\log(\mbox{Area}(\vec{\Phi}(\Sigma))
\ee
where, as before $g_0$ is a constant Gauss curvature metric of volume 1 on $\Sigma$, $g_{\vP}=e^{2\al}\, g_0$. The following lemma holds
\begin{Lm}
\label{lm-seq-weak}
Let  $\Sigma$ be a closed surface and $\vec{\Phi}^k$ be a sequence of weak immersions in ${\mathcal E}_{\Sigma,2}$ satisfying
\be
\label{VII.1}
\limsup_{k\rightarrow +\infty}F(\vec{\Phi}^k)<+\infty\quad
\ee
then the conformal class of the associated sequence $g^k_0$ of constant scalar curvature metric of volume 1 such that $g_{\vec{\Phi}^k}=e^{2\al^k}\,g_0^k$ is pre-compact in the moduli space ${\mathcal M}(\Sigma)$ moreover,
there exists a sequence of diffeomorphisms $\Psi^k$ of $\Sigma$ such that $(\Psi^k)^\ast g^k_0$ is converging strongly in any $C^l$ topology to a limiting constant curvature metric $h$,  $\vec{\Phi}^k\circ\Psi^k$  is conformal from $(\Sigma,(\Psi^k)^\ast g^k_0)$ and is sequentially weakly pre-compact in $W^{2,4}$ and for any weakly converging
sub-sequence the limit $\vec{\xi}^\infty$ is still a weak immersion of ${\mathcal E}_{\Sigma,2}$ and
\be
\label{VII.2}
\log|d(\vec{\Phi}^k\circ\Psi^k)|^2_{(\Psi^k)^\ast g^k_0}\rightharpoonup \log|d\vec{\xi}^\infty|^2_{h}\quad\quad\mbox{ weakly in }W^{1,4}_{h}
\ee
\hfill $\Box$
\end{Lm}

\noindent{\bf Proof of lemma~\ref{lm-seq-weak}.} We are working with a Aubin Gauge in the case $K_{g_0}>0$. Using respectively  lemma~\ref{lm-bound-bel-genus} for $K_{g_0}<0$ or lemma~\ref{lm-L-infty-borne} for $K_{g_0}>0$
we have in all cases
\be
\label{VII.12-a}
\limsup_{k\rightarrow +\infty}\int_{\Sigma} (1+H_{\vec{\Phi}^k}^2)^2\ dvol_{g_{\vec{\Phi}^k}}<+\infty\quad.
\ee
Hence using again (\ref{V.55-aa}) for genus$(\Sigma)>1$ or lemma~\ref{lm-L-infty-borne} for $K_{g_0}>0$ gives in all cases
\be
\label{VII.12-b}
\limsup_{k\rightarrow +\infty}\int_{\Sigma} |d\al^k|_{g^k_0}\ dvol_{g^k_0}<+\infty\quad.
\ee
Moreover since $\Sigma$ is a closed surface we have also using Willmore-Li-Yau universal Lower Bound of the Willmore energy (see \cite{Ri1})
\be
\label{VII.12-c}
16\pi^2\le\int_{\Sigma}e^{2\al^k}\ dvol_{g^k_0}\ \int_{\Sigma} H_{\vec{\Phi}^k}^4\ dvol_{g_{\vec{\Phi}^k}}
\ee
Hence, combining (\ref{VII.12-a}) and (\ref{VII.12-c}) we then have proved in all cases
\be
\label{VII.000}
\limsup_{k\rightarrow +\infty}\lf|\log\lf(\mbox{Area}(\vec{\Phi}^k(\Sigma))\rg)\rg|<+\infty\quad.
\ee
These preliminary estimates having been established we now prove the pre-compactness of the conformal class in the non zero genus case. The case when $\Sigma$ is a torus has already been considered in \cite{MR3}. So we can restrict to the case
where genus$(\Sigma)>1$. Assume the conformal class associated to $g_{\vec{\Phi}^k}$ and hence $g^k_0$ would degenerate we have a rather precise description of this degeneration (see \cite{Hum}). It requires the formation
of at least a {\it collar} which is a subdomain of $\Sigma$ diffeomorphic to an annulus that we identify to a cylinder of the form
\[
{\mathcal C}:=\lf\{(x_1,x_2)\ ;\ \frac{2\pi}{l^k}\varphi^k<x_2<\frac{2\pi}{l^k}(\pi-\varphi^k)\ ;\ 0\le x_1\le 2\pi \rg\}
\]
where the vertical lines  $x_1=0$ and $x_1=2\pi$ are identified, $l^k$ is the length  of a closed geodesic for the hyperbolic metric $g^k_0$
\[
l^k\longrightarrow 0
 \]and  $\varphi^k:=\arctan(\sinh(l^k/2)) $. The closed geodesic of length $l^k$ is given by $x_2=\pi^2/l^k$ . On this cylinder the hyperbolic metric $g^k_0$ has the following explicit expression
\[
g^k_0=\lf(\frac{l^k}{\ds 2\pi\,\sin\lf(\frac{l^k x_2}{2\pi}\rg)}\rg)^2\ \lf[dx^2_1+dx_2^2\rg]
\]
Denote in these coordinates 
\[
g_{\vec{\Phi}^k}:=e^{2\,\la^k}\ [dx_1^2+dx_2^2]\quad\mbox{ and }\quad\vec{e}_i^{\ k}:=e^{-\la^k}\,\p_{x_i}\vec{\Phi}^k
\]
The unit vector field $\vec{e}_1^{\ k}$ is tangent to a foliation of the image in ${\R}^3$ by $\vec{\Phi}^k$ of the collar region by circles. We apply Fenchel theorem to each of these  circles. Precisely for each $t\in(\frac{2\pi}{l^k}\varphi^k,\frac{2\pi}{l^k}(\pi-\varphi^k))$ we have
\be
\label{VII.3}
2\pi\le\int_{\{x_2=t\}\cap{\mathcal C}}\lf|\frac{\p \vec{e}_1^{\ k}}{\p x_1}\rg|\ dx_1
\ee
Integrating this inequality for $x_2$ between $\frac{2\pi}{l^k}\varphi^k$ and $\frac{2\pi}{l^k}(\pi-\varphi^k)$ and using Cauchy-Schwartz gives
\be
\label{VII.4}
\frac{(2\pi)^2}{l^k}\ (\pi-2\,\varphi^k)\le \frac{2\pi}{\sqrt{l^k}}\ \sqrt{\pi-2\,\varphi^k}\ \lf[\int_{\mathcal C}|\nabla\vec{e}_1^{\ k}|^2\ dx^2\rg]^{1/2}
\ee
We have $|\nabla\vec{e}_1^{\ k}|^2\le |\nabla\vec{n}_{\vec{\Phi}^k}|^2+|\nabla\la^k|^2$ moreover $\la^k=\al^k+\mu^k$ where
\[
\mu^k=2\ \log\lf(\frac{l^k}{2\pi}\rg)-2\ \log\lf( \sin\lf(\frac{l^k\ x_2}{2\pi}\rg)\rg)
\]
Thus
\[
|\nabla\mu^k|^2=\frac{(l^k)^2}{\pi^2}\ \frac{\cos^2\lf(\frac{l^k\ x_2}{2\pi}\rg)}{\sin^2\lf(\frac{l^k\ x_2}{2\pi}\rg)}
\]
Hence
\be
\label{VII.5}
\int_{\mathcal C}|\nabla\mu|^2\ dx^2\le \frac{(l^k)^2}{\pi^2}\int_{\frac{2\pi\varphi^k}{l^k}}^{\frac{2\pi\,(\pi-\varphi^k)}{l^k}}\frac{dx_2}{\sin^2\lf(\frac{l^k\ x_2}{2\pi}\rg)}\le\ C\ (l^k)^2\ \int^\infty_{\frac{2\pi\varphi^k}{l^k}}\frac{dx_2}{(l^k)^2\, x_2^2}\le C
\ee
and since from (\ref{VII.12-b}) we have that
\[
\limsup_{k\rightarrow +\infty}\int_{\mathcal C}|\nabla\al|^2\ dx^2=\limsup_{k\rightarrow +\infty}\int_{\mathcal C}|d\al|^2_{g_{\vec{\Phi}^k}}\ dvol_{g_{\vec{\Phi}^k}}<+\infty
\]
we deduce that 
\be
\label{VII.6}
\limsup_{k\rightarrow +\infty}\int_{\mathcal C}|\nabla\vec{e}_1^{\ k}|^2\ dx^2<+\infty
\ee
Combining (\ref{VII.4}) and (\ref{VII.6}) gives that $l^k$ is bounded from below by a positive number which contradicts the formation of a collar and the degeneracy of the conformal class of $[g_0^k]$. Modulo the 
composition with isometries $g_0^k$ is strongly converging in every Banach space $C^l(\Sigma)$. In order to simplify
the presentation we assume that $g_0^k$ is fixed. We cover the riemannian surface $(\Sigma,g_0)$ by finitely many conformal charts $\phi_j\, :\, D^2\ \rightarrow\ \phi_j(D^2)$ for $j\in J$ such that $\Sigma\subset \cup_{j}\phi_j(D^2_{1/2})$.
Denote again $\phi_j^\ast g_{\vec{\Phi}^k}=e^{2\la^k_j}\,[dx_1^2+dx_2^2]$ the expression of $g_{\vec{\Phi}^k}$  in the chart $\phi_j$. We have
\be
\label{VII.7}
\limsup_{k\rightarrow +\infty}\int_{D^2}|\nabla\vec{n}_{\vec{\Phi}^k\circ\phi_j}|^2+|\nabla\la^k_j|^2\ dx^2<+\infty
\ee
For $i=1,2$ we denote $\vec{e}_{j,i}^{\ k}:=e^{-\la^k_j}\p_{x_i}(\vec{\Phi}^k\circ\phi_j)$ and Liouville equation gives
\[
-\Delta \la^k_j=(\nabla\vec{e}_{j,1}^{\ k}\cdot\nabla^\perp\vec{e}_{j,2}^{\ k})\quad\quad\mbox{  in }D^2\quad.
\]
Inequality (\ref{VII.7}) implies
\be
\label{VII.8}
\limsup_{k\rightarrow +\infty}\int_{D^2}|\nabla\vec{e}_{j,i}^{\ k}|^2\ dx^2<+\infty
\ee
Combining (\ref{VII.7}) and (\ref{VII.8}) together with Wente estimates we obtain the existence of $\ov{\la^{\ k}_j}\in {\R}$ such that
\be
\label{VII.9}
\limsup_{k\rightarrow +\infty}\|\la^{\ k}_j-\ov{\la^{\ k}_j}\|_{L^\infty(D^2_{3/4})}<+\infty
\ee
Due to the connectedness of $\Sigma$ we deduce that
\be
\label{VII.10}
\sup_{j\ne l}\limsup_{k\rightarrow +\infty}\|\ov{\la_l^{\ k}}-\ov{\la^{\ k}_j}\|_{L^\infty(D^2_{3/4})}<+\infty
\ee
We have moreover
\be
\label{VII.11}
\limsup_{k\rightarrow+\infty}\int_{D^2} e^{2\,\la^k_j}\ dx^2<+\infty
\ee
We have
\be
\label{VII.12}
\begin{array}{l}
\ds4\pi\le\sum_{j\in J}\int_{D^2_{1/2}}H_{\vec{\Phi}^k\circ\phi_j}^2\ e^{2\,\la^k_j}\ dx^2\\[5mm]
\ds\quad\ \le \sum_{j\in J}\lf(\int_{D^2_{1/2}} e^{2\, \la_j^k}\ dx^2\rg)^{-1/2}\ \lf(\int_{D^2_{1/2}}H_{\vec{\Phi}^k\circ\phi_j}^4\ e^{2\,\la^k_j}\ dx^2\rg)^{1/2}
\end{array}
\ee
Since
\be
\label{VII.13}
\limsup_{k\rightarrow +\infty}\int_{D^2_{1/2}}H_{\vec{\Phi}^k\circ\phi_j}^4\ e^{2\,\la^k_j}\ dx^2<+\infty
\ee
we deduce from (\ref{VII.10}), (\ref{VII.11}) and (\ref{VII.12}) that
\be
\label{VII.14}
\max_{j\in J}\limsup_{k\rightarrow +\infty}\|\la^k_{j}\|_{L^\infty(D^2_{3/4})}<+\infty
\ee
Using
\[
\Delta(\vec{\Phi}^k\circ\phi_j)=e^{2\,\la_j^k}\ H_{\vec{\Phi}^k\circ\phi_j}
\]
together with (\ref{VII.13}) we deduce that $\vec{\Phi}^k\circ\phi_j$ is sequentially weakly pre-compact in $W^{2,4}(D^2_{1/2})$. Bootstrapping this information with (\ref{VII.14}) gives (\ref{VII.2}) and lemma~\ref{lm-seq-weak}
is proved.\hfill $\Box$
\subsection{The Palais-Smale Condition.}
The aim of the present subsection is to establish the following lemma
\begin{Lm}
\label{lm-palais}
Let  $\Sigma$ be a closed surface and $\vec{\Phi}^k$ be a sequence of weak immersions in ${\mathcal E}_{\Sigma,2}$ satisfying
\be
\label{VII.15}
\limsup_{k\rightarrow +\infty}F(\vec{\Phi}^k)<+\infty\quad,
\ee
where $F$ is given by (\ref{VI.0-0})
and such that we have 
\be
\label{VII.16}
\lim_{k\rightarrow +\infty}\sup_{\|\vec{w}\|_{\vec{\Phi}^k}\le 1}\p F(\vec{\Phi}^k)\cdot\vec{w}=0
\ee
where for any $\vec{\Phi}\in {\mathcal E}_{\Sigma,2}$ we denote
\[
\|\vec{w}\|_{\vec{\Phi}}^4:=\int_\Sigma\lf[|\nabla^{\gP} d\vec{w}|_{\gP}^4+|d\vec{w}|_{\gP}^4+|w|^4\rg]\ dvol_{\gP}
\]
then, modulo extraction of a subsequence, there exists a sequence of parametrization $\Psi^k$ such that $\vec{\Phi}^k\circ\Psi^k$ strongly converges in $W^{2,4}$ towards 
a critical point of $F^\sigma$ in ${\mathcal E}_{\Sigma,2}$
\hfill $\Box$
\end{Lm}
\noindent{\bf Proof of lemma~\ref{lm-palais}.} We are taking an Aubin Gauge for the $S^2$ case. We denote $f(t):=t^2+(1+t^2)^2$. We consider the charts $\phi^j$ given by the previous lemma and we omit to mention the index $j$. We also skip writing the index $k$ when it is not absolutely needed.
In each of this chart the information given by (\ref{VII.16}) is saying that 
\be
\label{palais-2}
\begin{array}{l}
\ds\vec{G}:=\mbox{div}\lf[\nabla (f'(H)\, \vec{n})-2\,  f'(H)\,\nabla\vec{n}-2\ e^{-2\la}\,  \nabla\vec{\Phi}\cdot\nabla\al\ \nabla\al\rg.\\[5mm]
\ds\quad+\lf.\ [-2\, f(H)+e^{-2\la}\,|\nabla\al|^2-K_{g_0}\, \al\ e^{-2\al}+K_{g_0}\,A_{\vec{\Phi}}(\Sigma)^{-1}]\ \nabla\vec{\Phi}-2\, \ e^{-2\la}\, \lf(\vec{\mathbb I}\res\nabla^\perp\al\rg)^\perp\rg]\\[5mm]
\ds\quad\quad\longrightarrow 0\quad\quad\mbox{strongly in }W^{-2,4/3}(D^2)=(W^{2,4}_0(D^2))^\ast
\end{array}
\ee
Let 
\[
\vec{M}:=\vec{G}\star\frac{1}{2\pi}\log |x| 
\]
This convolution is justified by a density argument mixed with Calderon Zygmund Theory and gives that
\be
\label{VII.17}
\vec{M}\longrightarrow 0\quad\quad\mbox{ strongly in }L^{4/3}(D^2)\quad.
\ee
Applying Poincar\'e Lemma we obtain the existence of $\vec{L}$ such that
\be
\label{VII.18}
\begin{array}{l}
\ds\nabla\vec{L}:=\nabla^\perp (f'(H)\, \vec{n})-2\,  f'(H)\,\nabla^\perp\vec{n}-2\ e^{-2\la}\, \sigma^2 \nabla\vec{\Phi}\cdot\nabla\al\ \nabla^\perp\al\\[5mm]
\ds\quad+\,\sigma^2\ [-2\, f(H)+e^{-2\la}\,|\nabla\al|^2-K_{g_0}\, \al\ e^{-2\al}+K_{g_0}\,A_{\vec{\Phi}}(\Sigma)^{-1}]\ \nabla^\perp\vec{\Phi}\\[5mm]
\ds\quad+2\, \ e^{-2\la}\, \vec{\mathbb I}\res\nabla^\perp\al-\nabla^\perp\vec{M}
\end{array}
\ee
 Equation (\ref{V.44}) gives
\be
\label{VII.19}
\nabla\vec{\Phi}\cdot\nabla^{\perp}\vec{L}-\mbox{div}\lf[ \nabla\vec{\Phi}\cdot \vec{M}\rg]=2\,\ e^{2\la}\, (2\,f(H)-\,H\, f'(H)+K_{g_0}\, \al\, e^{-2\al}-\,K_{g_0}\,A_{\vec{\Phi}}(\Sigma)^{-1})-\vec{M}\cdot\Delta\vec{\Phi}
\ee
Let $Y$ be the solution of
\be
\label{VII.20}
\lf\{
\begin{array}{l}
\ds\Delta Y=2\, e^{2\la}\, (2\,f(H)-\,H\, f'(H)+K_{g_0}\, \al\, e^{-2\al}-K_{g_0}\,A_{\vec{\Phi}}(\Sigma)^{-1})-\vec{M}\cdot\Delta\vec{\Phi}\quad\mbox{ in }D^2\\[5mm]
\ds Y=0\quad\quad\mbox{ on }\p D^2
\end{array}
\rg.
\ee
Observe that $2\,f(H)-\,H\, f'(H)=2\, (1-H^4)$. So $Y$ satisfies
\be
\label{VII.21}
\lf\{
\begin{array}{l}
\ds-\Delta Y=4\, e^{2\la}\,(1-H^4)+2\,K_{g_0}\, \al\, e^{2\mu}-2\,K_{g_0}\,A_{\vec{\Phi}}(\Sigma)^{-1}-\vec{M}\cdot\Delta\vec{\Phi}\quad\mbox{ in }D^2\\[5mm]
\ds Y=0\quad\quad\mbox{ on }\p D^2
\end{array}
\rg.
\ee
Since $\Delta\vec{\Phi}$ is uniformly bounded in $L^4$ we have $\vec{M}\cdot\Delta\vec{\Phi}\rightarrow 0$ strongly in $L^1$. Hence
\be
\label{VII.22}
\limsup_{k\rightarrow +\infty}\|\nabla Y_k\|_{L^{2,\infty}(D^2)}<+\infty
\ee
Using Poincar\'e theorem we deduce the existence of a function $S$ such that
\be
\label{VII.23}
\nabla S=\vec{L}\cdot\nabla\vec{\Phi}-\vec{M}\cdot\nabla^\perp\vec{\Phi}+\nabla^\perp Y
\ee
The equation (\ref{V.53}) in conformal coordinates gives
\be
\label{VII.24}
\nabla\vec{\Phi}\times\lf[\nabla^\perp\vec{L}-\nabla\vec{M}\rg]=-\,\nabla^\perp\vec{\Phi}\cdot\nabla f'(H)+2\,\nabla\al\cdot\nabla^\perp\vec{D}
\ee
where
\be
\label{VII.25}
\nabla\vec{D}=\lf(\ e^{-2\la} \sum_{i=1}^2{\mathbb I}_{1i}\,\p_{x_i}\vec{\Phi}\ ,\ e^{-2\la} \sum_{i=1}^2{\mathbb I}_{2i}\,\p_{x_i}\vec{\Phi}\rg)
\ee
Let $\vec{W}$ be the solution of
\be
\label{VII.26}
\lf\{
\begin{array}{l}
\ds\Delta\vec{W}=\vec{M}\times\Delta\vec{\Phi}\quad\quad\mbox{ in } D^2\\[5mm]
\ds\vec{W}=0\quad\quad\quad\mbox{ on }\p D^2
\end{array}
\rg.
\ee
Since $\Delta\vec{\Phi}$ is uniformly bounded in $L^4$ we have $\vec{M}\times\Delta\vec{\Phi}\rightarrow 0$ strongly in $L^1$. Hence
\be
\label{VII.28}
\lim_{k\rightarrow +\infty}\|\nabla \vec{W}_k\|_{L^{2,\infty}(D^2)}=0
\ee
Using again Poincar\'e Lemma on $D^2$ we obtain the existence of $\vec{V}$ such that
\be
\label{VII.27}
\nabla\vec{V}:=\vec{L}\times\nabla\vec{\Phi}-\vec{M}\times\nabla^\perp\vec{\Phi}+f'(H)\,\nabla\vec{\Phi}-2\, (\al-\ov{\al})\, \nabla\vec{D}+\nabla^\perp\vec{W}
\ee
Using the explicit expression of $\nabla\vec{D}$ given by (\ref{V.54}) we obtain
\be
\label{VII.29}
\begin{array}{l}
\vec{n}\cdot(\nabla\vec{V}-\nabla^\perp\vec{W})=\vec{n}\cdot(\vec{L}\times\nabla\vec{\Phi}-\vec{M}\times\nabla^\perp\vec{\Phi})\\[5mm]
\quad\quad=\vec{L}\cdot\nabla^\perp\vec{\Phi}+\vec{M}\cdot\nabla\vec{\Phi}=\nabla^\perp S+\nabla Y
\end{array}
\ee
We have also
\be
\label{VII.30}
\begin{array}{l}
\ds\vec{n}\times(\nabla\vec{V}-\nabla^\perp\vec{W})=-\, (\vec{L}\cdot\vec{n})\,\nabla\vec{\Phi}+(\vec{n}\cdot\vec{M})\nabla^\perp\vec{\Phi}\\[5mm]
\ds\quad\quad-\,\sigma^2\,f_\sigma'(H)\,\nabla^\perp\vec{\Phi}-2\,(\al-\ov{\al})\,\vec{n}\times\nabla\vec{D} 
\end{array}
\ee
Denote $\pi_T(\nabla^\perp \vec{V}+\nabla\vec{W})$ the tangential projection of $\nabla^\perp\vec{V}+\nabla\vec{W}$ , we have
\be
\label{VII.31}
\begin{array}{l}
\ds\pi_T(\nabla^\perp \vec{V}+\nabla\vec{W})= (\vec{L}\cdot\vec{n})\,\nabla\vec{\Phi}-(\vec{n}\cdot\vec{M})\nabla^\perp\vec{\Phi}\\[5mm]
\ds\quad\quad+\,f'(H)\,\nabla^\perp\vec{\Phi}-2\, (\al-\ov{\al})\, \nabla^\perp\vec{D}
\end{array}
\ee
Hence
\be
\label{VII.32}
\vec{n}\times(\nabla\vec{V}-\nabla^\perp\vec{W})=-\nabla^\perp \vec{V}-\nabla\vec{W}-\,2\, (\al-\ov{\al})\ \lf( \nabla^\perp\vec{D}+ \vec{n}\times\nabla\vec{D}\rg)-\vec{n}\,(\nabla S-\nabla^\perp Y)
\ee
Let $\vec{v}$ be the unique solution to
\be
\label{VII.33}
\lf\{
\begin{array}{l}
\ds\Delta\vec{v}=\nabla^\perp Y\cdot\nabla\vec{n}\quad\quad\mbox{ in } D^2\\[5mm]
\ds\vec{v}=0 \quad\quad\mbox{ on }\p D^2
\end{array}
\rg.
\ee
Using one more time Poincar\'e lemma we obtain the existence of $\vec{u}$ such that
\be
\label{VII.34}
\vec{n}\,\nabla^\perp{Y}=\nabla\vec{v}+\nabla^\perp \vec{u}
\ee
Finally, let $\vec{R}:=\vec{V}-\vec{u}$. We have
\be
\label{VII.35}
\vec{n}\times\nabla\vec{V}=\vec{n}\times\nabla\vec{R}+\vec{n}\times\nabla\vec{u}=\vec{n}\times\nabla\vec{R}+\vec{n}\times\nabla^\perp\vec{v}
\ee
Hence (\ref{VII.32}) becomes
\be
\label{VII.36}
\begin{array}{l}
\vec{n}\times(\nabla\vec{R}-\nabla^\perp\vec{W})+\vec{n}\times\nabla^\perp\vec{v}=-\nabla^\perp\vec{R}-\nabla\vec{W}+\nabla\vec{v}\\[5mm]
\ds\quad\quad-\vec{n}\,\nabla S-\,2\, (\al-\ov{\al})\ \lf( \nabla^\perp\vec{D}+ \vec{n}\times\nabla\vec{D}\rg)
\end{array}
\ee
which gives
\be
\label{VII.37}
\lf\{
\begin{array}{l}
\ds\Delta S=-\nabla\vec{n}\cdot\nabla^\perp{\vec{R}}+\mbox{div}\lf(\vec{n}\cdot\nabla\vec{z}-\,2\,(\al-\ov{\al})\ \vec{n}\cdot\nabla^\perp\vec{D}\rg)\\[5mm]
\ds\Delta\vec{R}=\nabla\vec{n}\times\nabla^\perp\vec{R}+\nabla\vec{n}\cdot\nabla^\perp S+\mbox{div}\lf(-\,\vec{n}\times\nabla\vec{z}+2\,(\al-\ov{\al})\   (-\nabla\vec{D}+  \vec{n}\times \nabla^\perp\vec{D}\rg)
\end{array}
\rg.
\ee
where $\vec{z}:=\vec{v}-\vec{W}$. Let $\vec{U}:=\vec{R}+2\,(\al-\ov{\al})\, \vec{D}$, with this notation (\ref{VII.37}) becomes
\be
\label{VII.38}
\lf\{
\begin{array}{l}
\ds\Delta S=-\nabla\vec{n}\cdot\nabla^\perp{\vec{U}}+\mbox{div}\lf(\vec{n}\cdot\nabla\vec{z}+\,2\, \vec{n}\cdot\vec{D}\ \nabla^\perp\al\rg)\\[5mm]
\ds\Delta\vec{U}=\nabla\vec{n}\times\nabla^\perp\vec{U}+\nabla\vec{n}\cdot\nabla^\perp S-\mbox{div}\lf(\vec{n}\times\nabla\vec{z}-\,2\,\nabla\al\ \vec{D}+  \vec{n}\times \vec{D}\ \nabla^\perp\al\rg)
\end{array}
\rg.
\ee
From lemma~\ref{lm-seq-weak} we know that $\nabla\vec{n}$ is weakly sequentially pre-compact in $L^4(D^2)$. The factor $\al$ satisfies the following Liouville type equation
\[
e^{-2\mu}\Delta\al=e^{2\al}\ K_{\gP}-K_{g_0}
\]
Since $\al$ is uniformly bounded in $L^\infty$ as well as $\mu$ on $D^2$ we deduce that $\Delta\al$ is uniformly bounded in $L^2$ and hence, $\nabla\al$ is \underbar{strongly} pre-compact in the
spaces $L^p_{loc}(D^2)$ for any $p<+\infty$. We know that $\nabla\vec{W}$ converges to zero strongly in $L^{2,\infty}$ and that $\nabla\vec{v}$ is strongly precompact in any space $L^q(D^2)$ for any $q<4$ hence, since also $\nabla\vec{D}$ is weakly pre-compact in $L^4$ and since we have chosen $\vec{D}$ to be of average zero on $D^2$, it is precompact in any $L^p_{loc}(D^2)$ for any $p<+\infty$. We can then apply lemma~\ref{lm-compcomp} and deduce that $\nabla S$ and $\nabla\vec{U}$ are \underbar{strongly pre-compact} in $L^{4/3}_{loc}(D^2)$. 

\medskip
 
\noindent Taking now the vectorial product between (\ref{VII.27}) and $\nabla^\perp\vec{\Phi}$ we obtain
\be
\label{VII.39}
\begin{array}{l}
\ds(\nabla\vec{V}-\nabla^\perp\vec{W})\times\nabla^\perp\vec{\Phi}\\[5mm]
\ds\quad=(\vec{L}\cdot\nabla^\perp\vec{\Phi}-\vec{M}\cdot\nabla\vec{\Phi})\cdot\nabla\vec{\Phi}-\,2\, f'(H)\, e^{2\la}\,\vec{n}-2\, (\al-\ov{\al})\, \nabla\vec{D}\times\nabla^\perp\vec{\Phi}+2\,\vec{M}\, e^{2\la}\\[5mm]
\ds\quad=\nabla^\perp S\cdot\nabla\vec{\Phi}+\nabla Y\cdot\nabla\vec{\Phi}-2\,(\vec{M}\cdot\nabla\vec{\Phi})\cdot\nabla\vec{\Phi} -2\, f'(H)\, e^{2\la}\,\vec{n}\\[5mm]
\ds\quad\quad-2\, (\al-\ov{\al})\, \nabla\vec{D}\times\nabla^\perp\vec{\Phi}+2\,\vec{M}\, e^{2\la}\quad.
\end{array}
\ee
We have also
\be
\label{VII.40}
\begin{array}{l}
\ds\nabla\vec{V}\times\nabla^\perp\vec{\Phi}=\nabla\vec{R}\times\nabla^\perp\vec{\Phi}+\nabla\vec{u}\times\nabla^\perp\vec{\Phi}\\[5mm]
\ds\quad=\nabla\vec{R}\times\nabla^\perp\vec{\Phi}+\nabla^\perp\vec{v}\times\nabla^\perp\vec{\Phi}+\nabla Y\cdot (\vec{n}\times\nabla^\perp\vec{\Phi})\\[5mm]
\ds\quad=\nabla\vec{R}\times\nabla^\perp\vec{\Phi}+\nabla\vec{v}\times\nabla\vec{\Phi}+\nabla Y\cdot\nabla\vec{\Phi}\quad.
\end{array}
\ee
Combining (\ref{VII.39}) and (\ref{VII.40}) gives
\be
\label{VII.41}
\begin{array}{l}
2\,  f'(H)\, e^{2\la}\,\vec{n}=\nabla^\perp S\cdot\nabla\vec{\Phi}-2\, (\al-\ov{\al})\, \nabla\vec{D}\times\nabla^\perp\vec{\Phi}-\nabla\vec{R}\times\nabla^\perp\vec{\Phi}-\nabla\vec{v}\times\nabla\vec{\Phi}\\[5mm]
\ds\quad\quad+2\,\vec{M}\ e^{2\la}-2\,(\vec{M}\cdot\nabla\vec{\Phi})\cdot\nabla\vec{\Phi}+\nabla\vec{W}\times\nabla\vec{\Phi}
\end{array}
\ee
Which implies
\be
\label{VII.42}
\begin{array}{l}
2\, f'(H)\, e^{2\la}\,\vec{n}=\nabla^\perp S\cdot\nabla\vec{\Phi}-\nabla\vec{U}\times\nabla^\perp\vec{\Phi}-\nabla\vec{v}\times\nabla\vec{\Phi}+2\, \vec{D}\times\nabla^\perp\vec{\Phi}\cdot\nabla\al\\[5mm]
\ds\quad\quad+2\,\vec{M}\ e^{2\la}-2\,(\vec{M}\cdot\nabla\vec{\Phi})\cdot\nabla\vec{\Phi}+\nabla\vec{W}\times\nabla\vec{\Phi}
\end{array}
\ee
From the results we established above we have then that $f'(H)$ is strongly pre-compact in $L^{4/3}_{loc}(D^2)$. We have explicitly $\, f'(H)=2\, H\,(3 +H^2))$. Denote
\be
\label{VII.43}
\vec{J}^\infty:=\lim_{k\rightarrow +\infty}\vec{H}^k\,(3+|H^k|^2)\quad\mbox{ strongly in }L^{4/3}_{loc}(D^2)
\ee

\medskip

Since $\nabla\al^k$ is strongly pre-compact in any $L^p_{loc}(D^2)$ for $p<+\infty$, this is also the case for $\nabla\la^k$. Moreover $\Delta\vec{\Phi}^k$ is uniformly bounded in $L^4$ so, modulo extraction of a subsequence we have
\[
\vec{H}^k\rightharpoonup \vec{H}^\infty:=2^{-1}\ e^{-2\la^\infty}\Delta\vec{\Phi}^\infty\mbox{ weakly in }L^{4}(D^2_{3/4})
\]
and $\vec{\Phi}^\infty$ is a conformal immersion of the disc $D^2_{1/2}$ in ${\mathcal E}_{D^2_{1/2},2}$. Observe that 
\[
\vec{H}^k\cdot \vec{H}^k\,(3+|H^k|^2)\rightharpoonup \vec{H}^\infty\cdot\vec{J}^\infty\in L^1(D^2_{1/2})\quad\mbox{ weakly in }{\mathcal D}'(D^2_{1/2})
\]
Hence the sequence of non-negative measures  $|H^k|^2\,(3+|H^k|^2)$ does not concentrate with respect to the Lebesgue measure :
\be
\label{VII.44}
\forall\ep>0\quad\exists\,\delta>0\quad\mbox{ s. t. }\forall\, E\subset D^2_{1/2}\ \,\mbox{ measur. }\quad |E|\le \delta\quad\Rightarrow\quad \limsup_{k\rightarrow +\infty}\int_E|H^k|^4\, dx^2\le \ep
\ee
From the strong convergence (\ref{VII.43}) we deduce an almost everywhere convergence
\[
\vec{H}^k\longrightarrow\ \vec{I}^\infty\quad\quad\mbox{ a. e. in }D^2
\]
Using Egorov theorem, for any $\delta>0$ there exists $E_\delta\subset D^2$ such that $|E_\delta|<\delta$ and 
\[
\vec{H}^k\longrightarrow\ \vec{I}^\infty\quad\quad\mbox{ uniformly in }D^2\setminus E_\delta
\]
Hence for any test function $\varphi\in C^\infty_0(D^2_{1/2})$ we have
\be
\label{VII.45}
\limsup_{k\rightarrow 0}\lf|\int_{D^2}\varphi\, \vec{H}^k\ dx^2-\int_{D^2\setminus E_\delta}\varphi\, \vec{I}^\infty\ dx^2\rg|\le\|\varphi\|_\infty\ \limsup_{k\rightarrow +\infty}\int_{E_\delta}|H^k| dx^2
\ee
Combining (\ref{VII.44}) and (\ref{VII.45}) we deduce that
\be
\label{VII.46}
\vec{I}^\infty=\vec{H}^\infty
\ee
and, using (\ref{VII.43}), that $\vec{H}^k$ converges strongly towards $\vec{H}^\infty$ in $L^4$. Hence $\vec{\Phi}^k$ is \underbar{strongly} pre-compact in $L^4(D^2_{1/2})$. Inserting this information in (\ref{palais-2}) gives that the limiting immersion satisfies the Euler Lagrange equation of $F$ which concludes the proof of lemma~\ref{lm-palais}.\hfill $\Box$
\subsection{Min-Max procedures for Frame Energies on $S^2$}

\subsubsection{The free case.}
Let ${\mathcal P}({\mathcal E}_{S^2})$ the space of subsets of the space of weak immersions ${\mathcal E}_{S^2}$. 

\begin{Dfi}
\label{df-admissible}
A non empty subset ${\mathcal A}$ of ${\mathcal P}({\mathcal E}_{\Sigma,2})$ is called admissible if for any homeomorphism $\varphi\in $ Hom$({\mathcal E}_{S^2,2})$
isotopic to the identity we have that
\[
\forall A\in {\mathcal A}\quad\quad\varphi(A)\in {\mathcal A}
\]
Moreover, there exists a topological set $X$ such that for any $A\in {\mathcal A}$ there exists $\vec{\Phi}^A$ in $C^0(X,{\mathcal E}_{S^2,2})$ such that
\[
A=\vec{\Phi}^A(X)\quad.
\]
\hfill $\Box$
\end{Dfi} 
Let now ${\mathcal A}$ admissible such that
\be
\label{VII.47}
\beta(0):=\inf_{A\in {\mathcal A}}\max_{\vec{\Phi}\in A} W(\vec{\Phi})>0
\ee
Since for any arbitrary fixed $\vec{\Phi}$ the map $\sigma\rightarrow F^\sigma(\vec{\Phi})$ is increasing
we can use a beautiful argument initially introduced by Michael Struwe in \cite{Str} and follow word by word the proof of theorem 6.4 in \cite{MiR} together with the Palais Smale property of $F^\sigma$ established in lemma~\ref{lm-palais}
in order to deduce the following lemma which was the main goal of the present subsection.
\begin{Lm}
\label{lm-sequence} Let ${\mathcal A}$ be an admissible family. There exists a sequence $\sigma^k\rightarrow 0$ and a sequence of critical points there exists a critical point
$\vec{\Phi}^{\sigma^k}$ of $F^{\sigma^k}$ such that
\be
\label{VII.48}
\beta(\sigma^k)=F^{\sigma^k}(\vec{\Phi}^{\sigma^k})\quad\mbox{and }\quad\p_\sigma F^{\sigma^k}(\vec{\Phi}^{\sigma^k})=o\lf(\frac{1}{\sigma^k\,\log(\sigma^k)^{-1}}\rg)
\ee
so in particular we have that
\be
\label{VII.49}
\lim_{k\rightarrow +\infty}W(\vec{\Phi}^{\sigma^k})=\beta(0)\quad.
\ee
\hfill $\Box$
\end{Lm}
\subsubsection{The Area Constrained Case.}

We are now adapting the previous case to the situation when we fix the area to be equal to 1. Precisely we define
\[
{\mathcal E}_{\Sigma,2}^1:=\lf\{ \vec{\Phi}\in {\mathcal E}_{\Sigma,2}\quad;\quad \int_\Sigma dvol_{\gP}=1\rg\}
\]
This is not difficult to check that it defines a Finsler Manifold structure based on the $W^{2,4}$ topology. The notion of admissible set
is defined identically as above but replacing general homeomorphisms of ${\mathcal E}_{\Sigma,2}$ by homeomorphisms of ${\mathcal E}_{\Sigma,2}^1$.
The construction of the pseudo-gradient restricted to ${\mathcal E}_{\Sigma,2}^1$ applies and we can follow again word by word the proof of theorem 6.4 in \cite{MiR} and the Palais Smale property of $F^\sigma$ established in lemma~\ref{lm-palais}
in order to deduce the same statement as lemma~\ref{lm-sequence} but under the area constraint $\mbox{Area}(\vec{\Phi}^{\sigma^k})=1$. 
We shall now establish the following lemma which is a consequence of the scaling invariance of the Willmore energy in ${\R}^m$
\begin{Lm}
\label{lagrange-multip}
Let $\vec{\Phi}^{\sigma^k}$ be a critical point of $F^{\sigma^k}$ under the constraint $\mbox{Area}(\vec{\Phi}^{\sigma^k})=1$ and satisfying
\be
\label{VII.49-a}
\p_\sigma F^{\sigma^k}(\vec{\Phi}^{\sigma^k})=o\lf(\frac{1}{\sigma^k\,\log(\sigma^k)^{-1}}\rg)
\ee
Then it satisfies the equation (for any choice of gauge $\al^k$ in the case $\Sigma= S^2$)
\be
\label{VII.50}
\begin{array}{l}
d\lf[\ast_{{g_{\vec{\Phi}^k}}} d[l_{\sigma^k} f_{\sigma^k}'(H^k)\,\vec{n}_{\vec{\Phi}^k}]-\, 2\,l_{\sigma^k}\,f_{\sigma^k}'(H^k)\, \ast_{{g_{\vec{\Phi}^k}}}d\vec{n}_{\vec{\Phi}^k}\rg.\\[5mm]
\ds\quad\quad\lf.+\,l_{\sigma^k}\, \lf[-\, 2\,f_{\sigma^k}(H^k)+|d\al^k|_{{g_{\vec{\Phi}^k}}}^2-\,K_{g_0}\,\al^k\,e^{-2\al^k}+\,K_{g_0}\,[A_{\vec{\Phi}^k}(\Sigma)]^{-1} \rg]\, \ast_{{g_{\vec{\Phi}^k}}} d\vec{\Phi}^k\rg.\\[5mm]
\ds\quad\quad\lf.-\,2\,l_{\sigma^k}\,\lf<d\vec{\Phi}^k,d\al^k\rg>_{{g_{\vec{\Phi}^k}}}\ \ast_{{g_{\vec{\Phi}^k}}}d\al^k+\,2\,l_{\sigma^k}\,\vec{\mathbb I}^{\,k}\res_{{g_{\vec{\Phi}^k}}}(\ast_{{g_{\vec{\Phi}^k}}}d\al^k)\rg]= C^k\ d\lf[\ast_{{g_{\vec{\Phi}^k}}}\, d\vec{\Phi}^k\rg]
\end{array}
\ee
where 
\be
\label{VII.51}
C^k=2\, (\sigma^k)^2 \int_{\Sigma}(1-|H^k|^4)+l_{\sigma^k}\,\int_\Sigma K_{g_0}\al\ dvol_{g_0}- l_{\sigma^k}\, K_{g_0}
\ee
Hence (for the choice of an Aubin gauge in the case $\Sigma=S^2$) we have
\be
\label{VII.52}
|C^k|=o(1)
\ee
\hfill $\Box$
\end{Lm}
\begin{Rm}
\label{rm-gauge-lagrange}
Observe that a-priori $C^k$ depends on the choice of gauge $\al^k$.\hfill $\Box$
\end{Rm}
\noindent{\bf Proof of Lemma~\ref{lagrange-multip}.}
We omit to write the index $k$. The fact that equation (\ref{VII.50}) is satisfied comes from (\ref{V.41}) and classical lagrange multiplier theory bearing in mind
that the first derivative of the fixed area constraint is proportional to $d\lf[\ast_{{g_{\vec{\Phi}^k}}}\, d\vec{\Phi}^k\rg]$ which cannot be zero since
obviously there is no compact minimal immersion in ${\R}^m$ (hence the constraint is non degenerate). We take the scalar product between (\ref{VII.50})
and $\vec{\Phi}$ and we integrate the resulting 2-form over the closed surface $\Sigma$. This gives
\be
\label{VII.53}
{C}=l_\sigma\int_{\Sigma} 2\,f_\sigma(H)-H\,f_\sigma'(H)+K_{g_0}\,\int_\Sigma\,\al\, e^{-2\al}\ dvol_{g_{\vec{\Phi}}}+\, 4\pi \ l_{\sigma}\,(\mbox{genus}(\Sigma)-1)
\ee
Since $2\,f_\sigma(H)-H\,f_\sigma'(H)=2\, \sigma^2\ l_\sigma^{-1}(1-H^4)$ we obtain (\ref{VII.51}). For the choice of an Aubin gauge we have due to (\ref{VII.49-a}) combined with  theorem~\ref{dirichlet-bound} and (\ref{V.69-a})
\be
\label{VII.54}
l_\sigma\|\al\|_{L^\infty(S^2)}=o(1)\quad.
\ee
This implies (\ref{VII.52}) for $\Sigma=S^2$. For $\Sigma\ne S^2$, lemma~\ref{lm-bound-bel-genus} combined with the assumption (\ref{VII.49-a}) implies
also
\be
\label{VII.55}
l_\sigma\|\al\|_{L^\infty(\Sigma)}=o(1)
\ee
 and lemma~\ref{lagrange-multip} is proved in any case. \hfill $\Box$

\medskip

For area constrained critical point we then have the following $\epsilon-$regularity lemma whose proof is following step by step the proof of 
lemma~\ref{lm-ep-reg} since (\ref{VII.52}) holds.
\begin{Lm}
\label{lm-ep-reg-constr} {\bf [uniform $\ep-$regularity under area constraint.]}
For any $C_0>0$, there exists $\ep>0$ and $\sigma_0>0$ such that for any $\sigma\in(0,\sigma_0)$ and any critical point $\vec{\Phi}$ of $F^\sigma$ under the constraint Area$(\vec{\Phi})=1$ satisfying
\be
\label{VI.0-constr}
F^\sigma(\vec{\Phi})\le C_1\quad\mbox{ and }\quad\p_\sigma F^{\sigma}(\vec{\Phi})\le \frac{\ep}{\sigma \log(\sigma)^{-1}}
\ee
assume moreover
\be
\label{VI.2-constr}
\int_{D^2}|\nabla\vec{n}|^2<\ep
\ee
then for any $j\in {\N}$ the estimates (\ref{VI.3}) and (\ref{VI.3-a}) hold in any case for $\Sigma\ne S^2$ and for any Coulomb gauge in case $\Sigma=S^2$.
\hfill $\Box$
\end{Lm}

\section{The passage to the limit $\sigma\rightarrow 0$.}
\reset
We shall give two results regarding the passage to the limit in the equation. A subsection will be devoted to each of the two results.
\subsection{The limiting immersions.}
In the sequel we shall denote by ${\mathcal M}^+(S^2)$ the non-compact M\"obius group of positive conformal diffeomorphism of the 2-sphere $S^2$.
\begin{Lm}
\label{lm-VIII-1}
Let $\sigma^k\rightarrow 0$ and a sequence of weak immersions $\vec{\Phi}^{k}\in {\mathcal E}_{S^2}$ which are critical points of $F^{\sigma_k}$ under area constraint equal to 1 and such that
\be
\label{VIII.1}
\limsup_{k\rightarrow +\infty}F^{\sigma^k}(\vec{\Phi}^{k})<+\infty\quad\mbox{and }\quad\p_\sigma F^{\sigma^k}(\vec{\Phi}^{k})=o\lf(\frac{1}{\sigma^k\,\log(\sigma^k)^{-1}}\rg)
\ee
Then modulo translation there exists
a subsequence that we still denote $\vec{\Phi}^k$, there exists a family of bilipschitz homeomorphism $\Psi^k$, there exists 
a finite family of sequences $(f_i^k)_{i=1\cdots N}$ of elements in ${\mathcal M}^+(S^2)$,
there exists a finite family of natural integers $(N_i)_{i=1\cdots N}$ and for each $i\in\{1\cdots N\}$ there exists finitely many points of $S^2$, $b_{i,1}\cdots b_{i,N_i}$ such that
\be
\label{VIII.5}
\vec{\Phi}^k\circ\Psi^k\longrightarrow \vec{f}^\infty\quad\mbox{ strongly in }C^0(S^2,{\R}^m)\quad,
\ee
where $\vec{f}^\infty\in W^{1,\infty}(S^2,{\R}^m)$, moreover
\be
\label{VIII.6}
\vec{\Phi}^k\circ f_i^k\rightharpoonup \vec{\xi}^\infty_i\quad\mbox{ strongly in } C^l_{loc}(S^2\setminus\{b_{i,1}\cdots b_{i,N_i}\})\quad,
\ee
for any $l\in {\N}$ where $\vec{\xi}^\infty_j$ is a Willmore conformal possibly branched immersion of $S^2$.  In addition we have 
\be
\label{VIII.7}
\ds\vec{f}^\infty(S^2)=\bigcup_{i=1}^N\vec{\xi}^\infty_i(S^2)\quad,
\ee
moreover
\be
\label{VIII.6b}
A(\vec{\Phi}_k)=\int_{S^2}1\ dvol_{g_{\vec{\Phi}_k}}\longrightarrow A(\vec{f}^\infty)=\sum_{i=1}^NA(\vec{\xi}^\infty_i)\quad,
\ee
and finally
\be
\label{VIII.8}
(\vec{f}^\infty)_\ast[S^2]=\sum_{i=1}^N(\vec{\xi}^\infty_i)_\ast[S^2]\quad,
\ee
where for any Lipschitz mapping $\vec{a}$ from $S^2$ into ${\R}^m$, $(\vec{a})_\ast[S^2]$ denotes the current
given by the push-forward by $\vec{a}$ of the current of integration over $S^2$ : for any smooth two-form
$\om$ on ${\R}^m$
\[
\lf<(\vec{a})_\ast[S^2],\om\rg>:=\int_{S^2}(\vec{a})^\ast\om\quad.
\] 
\hfill $\Box$
\end{Lm}
\begin{Rm}
\label{rm-VIII-1}
Lemma~\ref{lm-VIII-1} is ``detecting'' the Willmore spheres ``visible'' at the limit but is ignoring the possible ``asymptotic Willmore spheres'' which are shrinking and
disappearing at the limit. The detection of these asymptotic Willmore spheres is the purpose of the next subsection. Finally the detection of the possible
 loss of energy in the so called ``neck regions'' and the energy quantization effect is going to be investigated in section VII of the paper.
\end{Rm}
\noindent{\bf Proof of lemma~\ref{lm-VIII-1}.}
We are working modulo extraction of subsequences. Consider the various diffeomorphisms of $S^2$ $f^k_i$ given by theorem I.2 of \cite{MR1}. We choose the gauges given by $f^k_i$ that is the pairs $(f^k_i,\al^k_i)$ satisfying
\[
g_{\vec{\Phi}^k\circ f^k_i}=e^{2\al^k_i}\ g_{S^2}
\]
From the analysis in \cite{MR1} we have the existence of $N_i$ points of $S^2$, $b_{i,1}\cdots b_{i,N_i}$ such that
\be
\label{VIII.9}
\limsup_{k\rightarrow +\infty}\|\al^k_i\|_{L^\infty_{loc}(S^2\setminus\{b_{i,1}\cdots b_{i,N_i}\})}+\|\nabla\al^k_i\|_{L^2_{loc}(S^2\setminus\{b_{i,1}\cdots b_{i,N_i}\})}<+\infty
\ee
Hence we have 
\be
\label{VIII.10}
\lim_{k\rightarrow +\infty}l_{\sigma^k}\|\al^k_i\|_{L^\infty_{loc}(S^2\setminus\{b_{i,1}\cdots b_{i,N_i}\})}+\sqrt{l_{\sigma^k}}\, \|\nabla\al^k_i\|_{L^2_{loc}(S^2\setminus\{b_{i,1}\cdots b_{i,N_i}\})}=0
\ee
The assumptions (\ref{VIII.1}) implies moreover
\be
\label{VIII.11}
(\sigma^k)^2\ \int_{S^2}(1+H^2_{\vec{\Phi}^k})^2\ dvol_{\vec{\Phi}^k}=o(l_{\sigma^k})
\ee
Again from the analysis of \cite{MR1} we have that the density of energy $$|d\vec{n}_{\vec{\Phi}^k\circ f^k_i}|^2_{g_{\vec{\Phi}^k\circ f^k_i}}\ dvol_{g_{\vec{\Phi}^k\circ f^k_i}}$$ remains uniformly absolutely continuous with respect to Lebesgue in $S^2\setminus \cup_{j=1}^{N_i}B_\delta(b_{i,j})$. Hence all the assumptions which permit to apply the uniform $\epsilon-$regularity lemma~\ref{lm-ep-reg} are fulfilled and we deduce the strong convergence
(\ref{VIII.6}) towards a Willmore sphere that can be possibly branched at the $b_{i,j}$ but which is a \underbar{true} Willmore in the sense that the Willmore residues are zero
\[
\int_{\p B_{\delta}(b_{i,j})}\p_\nu \vec{H}_{\vec{\xi}^\infty_i}-\,2\,H_{\vec{\xi}^\infty_i}\,\p_\nu\vec{n}_{\vec{\xi}^\infty_i}-2\,H_{\vec{\xi}^\infty_i}^2\,\p_{\nu}\vec{\xi}^\infty_i\ dl=0
\] 
This concludes the proof of lemma~\ref{lm-VIII-1}.\hfill $\Box$
\subsection{Bubble Detection.}
\begin{Lm}
\label{lm-bubble-detection} {\bf [Bubble Detection Lemma.]}
Let $\Sigma$ be a closed surface, let $C_0>0$ and $\ep(C_0)>0$ given  by lemma~\ref{lm-ep-reg-constr} . Let $\sigma^k\rightarrow 0$ and $\vec{\Phi}^{\,k}$ be a sequence of $F^{\sigma^k}$ under the constraint Area$(\vec{\Phi}^k)=1$ and satisfying
\be
\label{VIII.12-add}
F^{\sigma^k}(\vec{\Phi}^k)\le C_0\quad\mbox{ and }\quad\p_\sigma F^{\sigma}(\vec{\Phi})=o\lf( \frac{1}{\sigma \log(\sigma)^{-1}}\rg)
\ee
assume when $\Sigma\ne S^2$ that the sequence of constant Gauss curvature metric $g_0^k$ such that $g_{\vec{\Phi}^k}=e^{2\al^k}\ g_0^k$ is 
pre-compact in any $C^l(\Sigma)$ topology. Let $B_{\rho^k}(p^k)$ be a geodesic ball for the metric $g_0^k$ such that  
\be
\label{VIII.13}
\int_{B_{\rho^k}(p^k)}|d\vec{n}|_{g_{\vec{\Phi}^k}}^2\ dvol_{g_{\vec{\Phi}^k}}<\ep
\ee
Compose $\vec{\Phi}^{\, k}$ in a sequence of conformal charts $(x^k_1,x^k_2)$  from $D^2$ into $(B_{\rho^k}(x^k),g_0^k)$ such that there exists $\ov{\mu}_k\in{\R}$ satisfying $g^k_0=e^{2\mu^k}\ [dx_1^2+dx_2^2]$ and $\mu^k-\ov{\mu}^k$ is strongly converging. Let $\ov{\al}^k$ be the average of $\al^k$ on $D^2_{1/2}$.
Then
\[
e^{-\ov{\mu}^k-\ov{\al}^k}\ \lf[\vec{\Phi}^{\,k}(x)-\vec{\Phi}^{\,k}(0)\rg]
\]
is strongly converging in $D^2_{1/2}$ in any $C^l$ norm towards a Willmore disc.
\hfill $\Box$
\end{Lm}
\begin{Rm}
\label{rm-VIII-1-add}
We are mostly interested with the balls $B_{\rho^k}(p^k)$ for which
\be
\label{VIII.14-add}
\limsup_{k\rightarrow 0}\int_{B_{\rho^k}(p^k)}|d\vec{n}|_{g_{\vec{\Phi}^k}}^2\ dvol_{g_{\vec{\Phi}^k}}>0
\ee
where there is indeed a non flat Willmore bubble forming which is swallowing  part of the energy. The lemma however applies
also when (\ref{VIII.14-add}) does not hold.\hfill $\Box$
\end{Rm}
\noindent{\bf Proof of lemma~\ref{lm-bubble-detection}.}
We keep denoting $\vec{\Phi}^{\, k}$ the composition of $\vec{\Phi}^{\, k}$ with the given chart and we work on $D^2$. Because of lemma~\ref{lagrange-multip} $\vec{\Phi}$ satisfies the equation (where we omit to write explicitly the upper index $k$)
\be
\label{VIII.15-add}
\begin{array}{l}
\ds e^{\ov{\la}}\ \mbox{div}\lf[\nabla\vec{H}_{\vec{\Phi}}-2\,H_{\vec{\Phi}}\,\nabla\vec{n}_{\vec{\Phi}}-2\,H^2_{\vec{\Phi}}\,\nabla\vec{\Phi}\rg]\\[5mm]
\ds\quad=\mbox{div}\lf[-\,2\,e^{\ov{\la}}\,\sigma^2\, \nabla\lf[  \vec{H}_{\vec{\Phi}}\, (1+H^2_{\vec{\Phi}})\rg]+4\,\sigma^2\ e^{\ov{\la}}\, H_{\vec{\Phi}}\, (1+H^2_{\vec{\Phi}})\ \nabla\vec{n}_{\vec{\Phi}}\rg]\\[5mm]
\ds\quad+\,\mbox{div}\lf[ \sigma^2\, (1+H^2_{\vec{\Phi}})^2\, e^{\ov{\la}}\ \nabla\vec{\Phi}+\,l_\sigma \ e^{\ov{\la}-2\la}\ \lf[ ({\mathbb I}_{\vec{\Phi}}\res\nabla^\perp\al)^\perp+\nabla\vec{\Phi}\cdot\nabla\al\ \nabla\al\rg]  \rg]\\[5mm]
\ds\quad-\,\mbox{div}\lf[2^{-1}\, l_\sigma\, \lf[ e^{-2\la}\,|\nabla\al|^2-\,K_{g_0}\, \al\ e^{-2\al}+K_{g_0}   \rg]\ e^{\ov{\la}}\, \rg]+\,2^{-1}\, C\, e^{\ov{\la}}\, \vec{H}_{\vec{\Phi}}
\end{array}
\ee
where $\ov{\la}$ is the average of $\la$ on $D^2_{1/2}$ and is equal to $\ov{\mu}+\ov{\al}$. The uniform $\epsilon-$regularity under area constraint, lemma~\ref{lm-ep-reg-constr}, combined with our assumption (\ref{VIII.12-add}) and
 lemma~\ref{lagrange-multip} imply that the right hand side of (\ref{VIII.15-add}) is converging towards 0 in any norm.
 Denote $\vec{\xi}:=e^{-\ov{\la}}(\vec{\Phi}(x)-\vec{\Phi}(0))$ we have
 \[
 \vec{n}_{\vec{\xi}}=\vec{n}_{\vec{\Phi}}\quad,\quad\vec{H}_{\vec{\xi}}= e^{\ov{\la}}\, \vec{H}_{\vec{\Phi}}\quad\mbox{ and }\quad\nabla\vec{\xi}= e^{-\ov{\la}}\,\nabla\vec{\Phi}
 \]
Hence we have
\[
\mbox{div}\lf[\nabla\vec{H}_{\vec{\xi}^{\,k}}-2\,H_{\vec{\xi}^{\,k}}\,\nabla\vec{n}_{\vec{\xi}^{\, k}}-2\,H^2_{\vec{\xi}^{\, k}}\,\nabla\vec{\xi}^{\,k}\rg]\longrightarrow 0\quad\mbox{ in }C^l(D^2_{1/2})\quad\forall\,l\in{\N}
\] 
and
\[
\limsup_{k\rightarrow +\infty }\int_{D^2_{1/2}}|\nabla\vec{n}_{\vec{\xi}^{\, k}}|^2<+\infty\quad\mbox{ and }\quad\limsup_{k\rightarrow +\infty}\log|\nabla\vec{\xi}^{\,k}|^2<+\infty
\] 
hence, since it is passing also through the origin, adapting the arguments in \cite{Ri2} to this perturbed case gives the strong convergence of $\vec{\xi}^{\,k}$
to a limiting Willmore immersion in $C^l(D^2_{1/2})\quad\forall\,l\in{\N}$ and the lemma~\ref{lm-bubble-detection} is proved. \hfill $\Box$

\subsection{The energy quantization}
\reset
With lemma~\ref{lm-bubble-detection} at hand, in order to prove our main result theorem~\ref{th-minmax}, following the scheme of the paper \cite{BR3}
it remains only to establish the vanishing of the energy in the so called  {\it neck region}. We are now restricting to the sphere case exclusively.
Precisely we are going to prove the following lemma. 
\begin{Lm}
\label{lm-energy-quantiz}
Let $\sigma^k\rightarrow 0$ and a sequence of weak immersions $\vec{\Phi}^{k}\in {\mathcal E}_{S^2}$ which are critical points of $F^{\sigma_k}$ 
under area constraint and such that
\be
\label{VIII.12}
\limsup_{k\rightarrow +\infty}F^{\sigma^k}(\vec{\Phi}^{k})<+\infty\quad\mbox{and }\quad\p_\sigma F^{\sigma^k}(\vec{\Phi}^{k})=o\lf(\frac{1}{\sigma^k\,\log(\sigma^k)^{-1}}\rg)
\ee
Then there exists $\vec{\xi}_1\cdots\vec{\xi}_n$, finitely  Willmore immersions of $S^2$ minus finitely many points, such that
\be
\label{VIII.14}
\lim_{k\rightarrow +\infty} W(\vec{\Phi}^k)=\sum_{i=1}^n W(\vec{\xi}_i)-\,4\pi\, N
\ee
where $N\in {\N}$. \hfill $\Box$
\end{Lm}
\noindent{\bf Proof of lemma~\ref{lm-energy-quantiz}.}
We shall work with an Aubin gauge $(\Psi^k,\al^k)$ satisfying
\[
g_{\vec{\Phi}^k\circ\Psi^k}=e^{2\,\al^k}\ \frac{g_{S^2}}{4\pi}\quad\mbox{ and }\quad\forall i=1,2,3\quad\int_{S^2}x_i\ e^{2\,\al^k}\ dvol_{S^2}=0
\]
where $g_{S^2}$ and $dvol_{S^2}$ are respectively the standard metric and standard associated volume form on $S^2$. In order to simplify the notation we omit
to write $\Psi_k$ that we assume to be the identity. The assumption (\ref{VIII.12}) reads as follows
\be
\label{VIII.15}
\begin{array}{l}
\ds o\lf(\frac{1}{\sigma^k\,\log(\sigma^k)^{-1}}\rg)
=2\,\sigma^k\int_{S^2}(1+H_{\vec{\Phi}^k}^2)^2\ e^{2\al_k}\ dvol_{g_{S^2}}\\[5mm]
\ds\quad+\frac{1}{\sigma^k\,(\log\sigma^k)^2}\lf[ \int_{S^2}2^{-1}|d\al^k|_{g_{S^2}}\ dvol_{g_{S^2}}+\int_{S^2}\al^k\ dvol_{g_{S^2}}-2\pi\,\log \mbox{Area}(\vec{\Phi}_k)     \rg]
\end{array}
\ee
The non negativity of both the first term and the second term in the r.h.s. of (\ref{VIII.15}) gives respectively
\be
\label{VIII.16}
(\sigma^k)^2\,\int_{S^2}(1+H_{\vec{\Phi}^k}^2)^2\ e^{2\al_k}\ dvol_{g_{S^2}}=o\lf(\frac{1}{\log(\sigma^k)^{-1}}\rg)
\ee
and
\be
\label{VIII.17}
\frac{1}{\log(\sigma^k)^{-1}}\lf[ \int_{S^2}2^{-1}|d\al^k|_{g_{S^2}}\ dvol_{g_{S^2}}+\int_{S^2}\al^k\ dvol_{g_{S^2}}-2\pi\,\log \mbox{Area}(\vec{\Phi}_k)     \rg]=o(1)
\ee
We keep using the notation from previous sections that is $l_{\sigma^k}:=1/\log(\sigma^k)^{-1}$. Using lemma~\ref{dirichlet-bound} we get
for this Aubin gauge
\be
\label{VIII.18}
l_{\sigma^k}\,\int_{S^2}|d\al^k|^2_{g_{S^2}}\ dvol_{g_{S^2}}=o(1)
\ee
and
\be
\label{VIII.19}
l_{\sigma^k}\,\|\al^k\|_{L^\infty(S^2)}=o(1)
\ee
Hence in order to apply the uniform $\epsilon$-regularity lemma~\ref{lm-ep-reg-constr} on any geodesic ball $B_r(x_0)$ for the $S^2$ metric it suffices
to assume that $\sigma$ is small enough and 
\[
\int_{B_r(x_0)} |d\vec{n}|^2_{g_{S^2}}\ dvol_{g_{S^2}}<\ep
\]
then in particular
\be
\label{VIII.20}
\begin{array}{l}
\ds r^4\, \|e^\la\, \nabla (\vec{H}\, (1+2\sigma^2(1+H^2))) \|^2_{L^\infty(B_{r/2}(x_0))}+r^2\, \sigma^4\, \|e^\la\,H (1+H^2)\|^2_{L^\infty(B_{r/2}(x_0))}\\[5mm]
\ds\quad+r^2\, \|\nabla \vec{n}\|^2_{L^\infty(B_{r/2}(x_0))}\le C\,\int_{B^2_r(x_0)}|\nabla\vec{n}|^2\ dx^2+ 
C\ \lf[\int_{B_{r}(x_0)}[l_{\sigma^k}\,e^{2\al^k}+l_{\sigma^k}\,\|\al^k\|_{L^\infty(S^2)} ]\ dvol_{g_{S^2}}\rg]^2\\[5mm]
\ds\quad+ C\,\lf[\int_{B_r^2(x_0)}(\sigma^k)^2H^4\ e^{2\la}\ dx^2\rg]^2+\,C\,\lf[\int_{B_r^2(x_0)}l_{\sigma^k}\ |\nabla\al^k|^2\ dx^2\rg]^2
\end{array}
\ee
and
\be
\label{VIII.21}
\begin{array}{l}
\ds r^2\ l_{\sigma^k}\,\|\nabla\al^k\|_{L^\infty(B_{r/2}(x_0))}^2\le\, C\, l_{\sigma^k}\,\int_{B^2_r(x_0)}|\nabla\al^k|^2\ dx^2+ \,C\,\lf[\int_{B^2_r(x_0)}|\nabla\vec{n}|^2\ dx^2\rg]^2\\[5mm]
\ds\quad+C\,\lf[\int_{B_r^2(x_0)}(\sigma^k)^2H^4\ e^{2\la}\ dx^2\rg]^4+ 
C\, \lf[\int_{B_{r}(x_0)}(\sqrt{l_{\sigma^k}}+l_{\sigma^k}\,e^{2\al^k})\ dvol_{g_{S^2}}\rg]^2
\end{array}
\ee
Recall that the fact that $\vec{\Phi}^k$ is a critical point of $ F^{\sigma^k}$ is equivalent to the existence of $\vec{L}$ such that
\be
\label{VIII.22-a}
\begin{array}{l}
\ds\nabla\vec{L}:=2\,\nabla^\perp \lf(\vec{H}\,(1+2\,\sigma^2\,(1+H^2)\rg)-4\,H\,(1+2\,\sigma^2\,(1+H^2))\,\nabla^\perp\vec{n}\\[5mm]
\ds\quad-2\ e^{-2\la}\, l_\sigma\, \nabla\vec{\Phi}\cdot\nabla\al\ \nabla^\perp\al+2\, l_\sigma\ e^{-2\la}\, \vec{\mathbb I}\res\nabla^\perp\al
-\,2\, (H^2+\sigma^2(1+H^2)^2)\, \nabla^\perp\vec{\Phi}\\[5mm]
\ds\quad+l_\sigma\ \lf[e^{-2\la}\,|\nabla\al|^2-K_{g_0}\, \al\ e^{-2\al}+4\pi\,A_{\vec{\Phi}}(S^2)^{-1}\rg]\ \nabla^\perp\vec{\Phi}
\end{array}
\ee
Since $2\,\sigma^2\ (1+H^2) |H|\, |\nabla\vec{n}|\, e^{\la}\le |\nabla\vec{n}|^2\, +\sigma^4\, H^2 (1+H^2)^2\,e^{2\la}$,
the previous estimates imply
\be
\label{VIII.22}
\begin{array}{l}
\ds r^2\|e^\la\,\nabla \vec{L}\|_{L^\infty(B_{r/2}(x_0))}\le\, C\, \lf[\int_{B^2_r(x_0)}|\nabla\vec{n}|^2\ dx^2\rg]^{1/2}\, +\,C\, l_{\sigma^k}\,\int_{B^2_r(x_0)}|\nabla\al^k|^2\ dx^2\\[5mm]
\ds\quad+ \int_{B_r^2(x_0)}(\sigma^k)^2H^4\ e^{2\la}\ dx^2+l_\sigma\, \int_{B_r(x_0)}\, dvol_{g_{\vec{\Phi}}}+l_\sigma\,\|\al\|_\infty\ \int_{B_r(x_0)}\, dvol_{g_{S^2}}
\end{array}
\ee
We now follow step by step the arguments from section VI of \cite{BR3} and check how each estimate is slightly modified by the viscous terms.
We consider a neck region that is an annulus (for the $S^2$ metric) of the form $B^2_{R^k}(0)\setminus B^2_{r^k}(0)$ where
\[
\lim_{k\rightarrow+\infty}R^k=0\quad\quad\mbox{and}\quad\quad\lim_{k\rightarrow+\infty}\frac{r^k}{R^k}=0
\]
and such that 
\be
\label{VIII.23}
\lim_{k\rightarrow+\infty}\sup_{r^k<s<R^k/4}\int_{B^2_{4s}(0)\setminus B^2_s(0)}|\nabla\vec{n}_{\vec{\Phi}^k}|^2\ dx^2=0
\ee
we shall omit very often when it is not necessary to write explicitly the superscript $k$. Denote for $s<R^k/4$
\be
\label{VIII.24}
\begin{array}{l}
\ds s\ \delta(s):=\lf[\int_{B^2_{2s}(0)\setminus B^2_{s/2}(0)}|\nabla\vec{n}|^2\ dx^2\rg]^{1/2}\, +\,C\,\lf[ l_{\sigma^k}\,\int_{B^2_{2s}(0)\setminus B^2_{s/2}(0)}|\nabla\al^k|^2\ dx^2\rg]^{1/2}\\[5mm]
\ds\quad\quad\quad+ \lf[\int_{B^2_{2s}(0)\setminus B^2_{s/2}(0)}(\sigma^k)^2 (1+H^2)^2\ e^{2\la}\ dx^2\rg]^{1/2}+\lf[l_\sigma\, \int_{B^2_{2s}(0)\setminus B^2_{s/2}(0)}\, dvol_{g_{\vec{\Phi}}}\rg]^{1/2}\\[5mm]
\ds\quad\quad\quad+\lf[l_\sigma\, \|\al\|_\infty\ \int_{{B^2_{2s}(0)\setminus B^2_{s/2}(0)}}\, dvol_{g_{S^2}}\rg]^{1/2}=o(1)
\end{array}
\ee
where $o(1)$ is chosen in such a way that, from (\ref{VIII.22}), we have in the neck region
\be
\label{VIII.25}
|x|^2 |\nabla\vec{L}|(x)\le \delta(|x|)\ e^{-\la(x)}
\ee
and such that
\be
\label{VIII.26}
\limsup_{k\rightarrow +\infty}\int_r^{R/4}\delta^2(s)\ s\ ds<+\infty
\ee
We keep up with the notations in \cite{BR3} and introduce
\[
\vec{L}_t:=\frac{1}{|\p B_t(0)|}\int_{\p B_t(0)}\vec{L}\quad\mbox{ and }\quad\la(t):=\frac{1}{|\p B_t(0)|}\int_{\p B_t(0)}\la.
\]
Using (\ref{VIII.22}) we have
\be
\label{VIII.27}
\begin{array}{l}
\ds\frac{d\vec{L}_t}{dt}=\frac{1}{2\pi}\int_{0}^{2\pi}\frac{\p \vec{L}}{\p t}(t,\theta)\ d\theta=-\frac{2}{\pi}\int_{0}^{2\pi}H\,(1+2\,\sigma^2\,(1+H^2))\,\frac{1}{t}  \frac{\p\vec{n}}{\p \theta}\ d\theta\\[5mm]
\ds \quad-\frac{1}{\pi}\int_{0}^{2\pi} e^{-2\la}\, l_\sigma\,\nabla\vec{\Phi}\cdot\nabla\al\ \frac{1}{t}  \frac{\p\al}{\p \theta}+\frac{1}{\pi}\int_{0}^{2\pi} e^{-2\la}\, l_\sigma\,\lf(\vec{\mathbb I}\res\nabla^\perp\al\rg)\cdot\frac{\p}{\p r}\ d\theta\\[5mm]
\ds\quad\quad-\frac{1}{\pi}\int_{0}^{2\pi}(H^2+\sigma^2\,(1+H^2)^2)\,\frac{1}{t}  \frac{\p\vec{\Phi}}{\p \theta}\ d\theta\\[5mm]
\ds\quad\quad+\frac{1}{2\pi}\int_{0}^{2\pi} l_\sigma\ \lf[e^{-2\la}\,|\nabla\al|^2-K_{g_0}\, \al\ e^{-2\al}+4\pi\,A_{\vec{\Phi}}(S^2)^{-1}\rg]\ \,\frac{1}{t}  \frac{\p\vec{\Phi}}{\p \theta}\ d\theta
\end{array}
\ee
This gives using again $\sigma^2 (1+H^2)^2\le 2\, \sigma^4 (1+H^2)^2 H^2+3\, (1+H^2)$
\be
\label{VIII.28}
\begin{array}{l}
\ds e^{\la(t)}\,\lf|\frac{d\vec{L}}{dt}\rg|(t)\le C\ \sigma^4\|e^{\la} H\, (1+H^2)\|^2_{L^\infty(\p B_t(0))}+\,C\, \|\nabla\vec{n}\|^2_{L^\infty(\p B_t(0))}\\[5mm]
\ds\quad\quad+C\,l_\sigma\, \|\nabla\al\|^2_{L^\infty(\p B_t(0))}+ l_\sigma\,e^{2\,\la(t)} +l_\sigma\,\|\al\|_\infty
\end{array}
\ee
Combining the previous inequality with (\ref{VIII.20}) and (\ref{VIII.24}) we finally obtain
\be
\label{VIII.29}
e^{\la(t)}\,\lf|\frac{d\vec{L}}{dt}\rg|(t)\le C\, \delta^2(t)
\ee
Following the arguments in \cite{BR3} we can choose a normalization in such a way that
\be
\label{VIII.30}
e^{\la(x)}|\vec{L}|(x)\le C\ |x|^{-1}\quad\mbox{on }B_{R/4}\setminus B_{2r}
\ee
where $C$ is independent of $k$. We adopt the notations of the proof of lemma~\ref{lm-ep-reg}. Let $Y$ satisfying
\be
\label{VIII.31}
\lf\{
\begin{array}{l}
\ds-\Delta Y=4\, e^{2\la}\, \sigma^2\, (1-H^4)+2\,l_\sigma\,K_{g_0}\, \al\, e^{2\mu}-8\pi\,l_\sigma\,e^{2\la}\,A_{\vec{\Phi}}(S^2)^{-1}\quad\mbox{ in }B_R(0)\\[5mm]
\ds Y=0\quad\quad\mbox{ on }\p B_R(0)
\end{array}
\rg.
\ee
Inequality gives (\ref{VI.38})
\be
\label{VIII.32}
\begin{array}{l}
\ds\|\nabla Y\|_{L^{2,\infty}(B_R(0))}\le C\ \int_{B_R(0)} \sigma^2\, [1+H^4]\ e^{2\la}\ dx^2+\\[5mm]
\ds\quad\quad\quad\quad\quad+C\, l_\sigma\, \|\al\|_\infty\int_{B_R(0)}e^{2\mu}\ dx^2+C\,l_\sigma\,\frac{A_{\vec{\Phi}}(B_R(0))}{A_{\vec{\Phi}}(S^2)}=o(1)
\end{array}
\ee
On $B_{2t}\setminus B_{t/2}$  we have
\be
\label{VIII.33}
\|\Delta Y\|_{L^\infty(B_{2t}\setminus B_{t/2})}\le C\, \delta^2(t)
\ee
Hence we deduce using standard interpolation theory
\be
\label{VIII.34}
t^2\,\|\nabla Y\|^2_{L^\infty(\p B_t(0))}\le C\,\|\nabla Y\|_{L^{2,\infty}(B_R(0))}\, t^2\, \|\Delta Y\|_{L^\infty(B_{2t}\setminus B_{t/2})}+\|\nabla Y\|^2_{L^{2,\infty}(B_R(0))}
\ee
Combining (\ref{VIII.32}), (\ref{VIII.33}) and (\ref{VIII.34}) give
\be
\label{VIII.35}
\|\nabla Y\|_{L^\infty(\p B_t(0))}\le o(1)\ [\delta(t)+t^{-1}]=o(1)\, t^{-1}
\ee
Using Poincar\'e theorem we deduce the existence of a function $S$ such that
\be
\label{VIII.36}
\nabla S=\vec{L}\cdot\nabla\vec{\Phi}+\nabla^\perp Y\quad\mbox{ on }B_R(0)
\ee
Combining (\ref{VIII.30}) and (\ref{VIII.35}) gives
\be
\label{VIII.37}
\|\nabla S\|_{L^\infty(\p B_t(0))}\le\ C\, t^{-1}\quad\quad\mbox{ for }t\in [2r,R/2]
\ee
Let $\vec{v}$ satisfying
\be
\label{VIII.38}
\lf\{
\begin{array}{l}
\ds\Delta\vec{v}=\nabla^\perp Y\cdot\nabla\vec{n}\quad\mbox{ in }B_R(0)\\[5mm]
\ds \vec{v}=0\quad\quad\mbox{ on }\p B_R(0)
\end{array}
\rg.
\ee
This gives using Wente estimates (see \cite{Ri1})
\be
\label{VIII.39}
\int_{B_R(0)}|\nabla \vec{v}|^2\ dx^2\le C\,\|\nabla Y\|^2_{2,\infty}\ \int_{B_R(0)}|\nabla\vec{n}|^2\ dx^2=o(1)
\ee
and
\be
\label{VIII.40}
t^2\,\|\nabla\vec{v}\|^2_{L^\infty(\p B_t(0))}\le C\, \|\nabla\vec{v}\|_{L^2}\ t^2\ \|\Delta\vec{v}\|_{L^\infty(B_{2t}\setminus B_{t/2})}+\int_{B_{2t}\setminus B_{t/2}}|\nabla \vec{v}|^2\ dx^2
\ee
Let 
\[
\eta_1(t):=t^{-1}\ \lf[\int_{B_{2t}\setminus B_{t/2}}|\nabla \vec{v}|^2\ dx^2\rg]^{1/2}=o(1)\, t^{-1}
\]
We have
\be
\label{VIII.40-a}
\int_r^R\, \eta_1^2(t)\, t\ dt=o(1)
\ee
Hence using  estimates (\ref{VIII.40}) we have
\be
\label{VIII.41}
\|\nabla \vec{v}\|_{L^\infty(\p B_t(0))}\le o(1)\ \delta(t)+\eta_1(t)=o(1)\, t^{-1}
\ee
Let now $\vec{u}$ such that
\be
\label{VIII.42}
\vec{n}\,\nabla^\perp Y=\nabla\vec{v}+\nabla^\perp\vec{u}
\ee
It satisfies
\be
\label{VIII.43}
\|\nabla \vec{u}\|_{L^\infty(\p B_t(0))}\le o(1)\, t^{-1}
\ee
Let
\be
\label{VIII.44}
\nabla\vec{D}=\lf(l_\sigma\ e^{-2\la} \sum_{i=1}^2{\mathbb I}_{1i}\,\p_{x_i}\vec{\Phi}\ ,\ l_\sigma\ e^{-2\la} \sum_{i=1}^2{\mathbb I}_{2i}\,\p_{x_i}\vec{\Phi}\rg)
\ee
Using again Poincar\'e Lemma on $B_R(0)$ we obtain the existence of $\vec{V}$ such that
\be
\label{VIII.45}
\nabla\vec{V}:=\vec{L}\times\nabla\vec{\Phi}+\, 2\, H (1+2\,\sigma^2\,(1+H^2))\,\nabla\vec{\Phi}-2\, \al\, \nabla\vec{D}
\ee
Using again $2\, |H|\,\sigma^2\,(1+H^2))\le 2\, \sqrt{H^2\, e^{2\la}\,\sigma^4\ (1+H^2)^2}$ together with (\ref{VIII.20}) and using also
(\ref{VIII.30})  we obtain
\be
\label{VIII.46}
\|\nabla\vec{V}\|_{L^\infty(\p B_t(0))}\le C\, [t^{-1}+\delta(t)]\le 2\, C\, t^{-1}
\ee
We denote $\vec{R}:=\vec{V}-\vec{u}$ and (\ref{VI.16}) implies
\be
\label{VIII.47}
\vec{n}\times\nabla\vec{R}+\vec{n}\times\nabla^\perp\vec{v}=-\nabla^\perp\vec{R}+\nabla\vec{v}-\vec{n}\,\nabla S-\,2\, \al\ \lf( \nabla^\perp\vec{D}+ \vec{n}\times\nabla\vec{D}\rg)
\ee
Combining (\ref{VIII.43}) and (\ref{VIII.46}) we have 
\be
\label{VIII.47-a}
\|\nabla\vec{R}\|_{L^\infty(\p B_t(0))}\le  C\, t^{-1}
\ee
Let $\vec{E}$ be the solution of 
\be
\label{VIII.48}
\lf\{
\begin{array}{l}
\ds\Delta\vec{E}=2\,\nabla\al\cdot\nabla^\perp\vec{D}\quad\mbox{on }B_R(0)\\[5mm]
\ds\vec{E}=0\quad\quad\mbox{ on }\p \vec{D}
\end{array}
\rg.
\ee
We have using Wente estimates (see \cite{Ri1})
\be
\label{VIII.49}
\|\nabla\vec{E}\|_{L^{2}(B_R(0))}\le\|\nabla\vec{E}\|_{L^{2,1}(B_R(0))}\le C\, l_\sigma\, \|\nabla\al\|_{L^2(B_R(0))}\, \|\nabla\vec{n}\|_{L^2(B_R(0))}=o(1)\quad.
\ee
Denote
\[
 \eta_2(t):=t^{-1}\ \lf[\int_{B_{2t}\setminus B_{t/2}}|\nabla \vec{E}|^2\ dx^2\rg]^{1/2}=o(1)\, t^{-1}
\]
From (\ref{VIII.49}) we obtain
\be
\label{VIII.49-a}
\int_{2r}^{R/2}\eta_2^2(t)\ t\, dt=o(1)
\ee
Interpolation inequalities give again
\be
\label{VIII.50}
\begin{array}{l}
\ds t^2\, \|\nabla\vec{E}\|^2_{L^\infty(\p B_t(0))}\le C\, \|\nabla\vec{E}\|_{L^2(B_R(0))}\ t^2\ \|\Delta\vec{E}\|_{L^\infty(B_{2t}\setminus B_{t/2})}+ C\, \int_{B_{2t}\setminus B_{t/2}}|\nabla \vec{E}|^2\ dx^2\\[5mm]
\ds \quad\quad\le C\,t^2\  o(1)\, l_\sigma \|\nabla\al\|_{L^\infty(B_{2t}\setminus B_{t/2})}\ \|\nabla\vec{n}\|_{L^\infty(B_{2t}\setminus B_{t/2})}+t^2\ \eta_2^2(t)\, \\[5mm]
\ds\quad\quad\le C\, t^2\ [o(1)\ \delta^2(t)+\eta_2^2(t)]=o(1)
\end{array}
\ee
Hence we deduce
\be
\label{VIII.51}
\|\nabla\vec{E}\|_{L^\infty(\p B_t(0))}\le C\ [o(1)\, \delta(t)+\eta_2(t)]=o(1)\, t^{-1}
\ee
Let $\vec{F}$ such that 
\[
2\,\al\,\nabla^\perp \vec{D}=\nabla^\perp\vec{F}+\nabla\vec{E}
\]
Observe that we have in one hand
\be
\label{VIII.52}
\|\nabla\vec{F}\|_{L^2(B_R(0))}\le 2\, l_\sigma\,\|\al\|_\infty\ \|\nabla\vec{n}\|_{L^2(B_R(0))}+\|\nabla\vec{E}\|_{L^2(B_R(0))}=o(1)
\ee
and in the other hand for $t\in [4r,R/4]$
\be
\label{VIII.53}
\begin{array}{l}
\ds\|\nabla\vec{F}\|_{L^\infty(\p B_t(0))}\le 2\, l_\sigma\,\|\al\|_\infty\ \|\nabla\vec{n}\|_{L^\infty(\p B_t(0))}+\|\nabla\vec{E}\|_{L^\infty(\p B_t(0))}\\[5mm]
\ds\quad\quad=C\ [o(1)\, \delta(t)+\eta_2(t)]=o(1)\, t^{-1}
\end{array}
\ee
Let $\vec{X}:=\vec{R}+\vec{F}$ we have using (\ref{VIII.47-a})
\be
\label{VIII.53-a}
\|\nabla\vec{X}\|_{L^\infty(\p B_t(0))}\le C\, t^{-1}
\ee
and it satisfies
\be
\label{VIII.54}
\vec{n}\times\nabla\vec{X}+\vec{n}\times\nabla^\perp\vec{v}=-\nabla^\perp\vec{X}+\nabla\vec{v}-\vec{n}\,\nabla S-\nabla\vec{E}+ \vec{n}\times\nabla^\perp\vec{E}
\ee
Let $\vec{w}$ be the solution of
\be
\label{VIII.55}
\lf\{
\begin{array}{l}
\ds\Delta\vec{w}=\nabla\vec{n}\cdot\nabla^\perp (\vec{v}-\vec{E})\quad\mbox{on }B_R(0)\\[5mm]
\ds\vec{w}=0\quad\quad\mbox{ on }\p B_R(0)
\end{array}
\rg.
\ee
Using Wente estimates we obtain
\be
\label{VIII.56}
\|\nabla\vec{w}\|_{L^{2}(B_R(0))}\le\|\nabla\vec{w}\|_{L^{2,1}(B_R(0))}\le C\,  \|\nabla(\vec{v}-\vec{E})\|_{L^2(B_R(0))}\, \|\nabla\vec{n}\|_{L^2(B_R(0))}=o(1)\quad.
\ee
Denote
\[
 \eta_3(t):=t^{-1}\ \lf[\int_{B_{2t}\setminus B_{t/2}}|\nabla \vec{w}|^2\ dx^2\rg]^{1/2}=o(1)\, t^{-1}
\]
From (\ref{VIII.56}) we obtain
\be
\label{VIII.56-a}
\int_{2r}^{R/2}\eta_3^2(t)\ t\, dt=o(1)
\ee
Interpolation inequalities give again using (\ref{VIII.41}) we obtain
\be
\label{VIII.57}
\|\nabla\vec{w}\|_{L^\infty(\p B_t(0))}\le C\ [o(1)\, \delta(t)+\eta_3(t)]=\ o(1)\ t^{-1}
\ee
Let $\vec{Z}$ such that
\[
\vec{n}\times\nabla^\perp(\vec{v}-\vec{E})=\nabla^\perp\vec{Z}+\nabla\vec{w}\quad.
\]
Clearly from the above we have also
\be
\label{VIII.58}
\begin{array}{l}
\ds\|\nabla\vec{Z}\|_{L^{2}(B_R(0))}=o(1)\\[5mm]
\ds\mbox{ and }\quad \|\nabla\vec{Z}\|_{L^\infty(\p B_t(0))}\le C\ \lf[o(1)\ \delta(t)+\sum_{i=1}^3\eta_i(t)\rg]=o(1)\ t^{-1}
\end{array}
\ee
Let $\vec{T}:=\vec{X}+\vec{Z}$. We have
\be
\label{VIII.59}
\nabla\vec{T}=-\nabla^\perp\vec{v}+\nabla^\perp\vec{E}+\vec{n}\ \nabla^\perp S+\vec{n}\times\nabla^\perp\vec{X}
\ee
which implies in particular
\be
\label{VIII.60}
\nabla^\perp S+\vec{n}\cdot\nabla^\perp(\vec{E}-\vec{v})=\vec{n}\cdot\nabla\vec{T}
\ee
Regarding the estimates we have from (\ref{VIII.53-a}) and (\ref{VIII.58}) we have
\be
\label{VIII.60-a}
\|\nabla\vec{T}\|_{L^\infty(\p B_t(0))}\le C\ t^{-1}\quad\quad\Longrightarrow\quad\quad \|\nabla\vec{T}\|_{L^{2,\infty}(B_{R}\setminus B_r)}\le C
\ee
where the constant $C$ (as for all constants $C$ above) is independent of $k$. Let $B$ be the solution of
\be
\label{VIII.61}
\lf\{
\begin{array}{l}
\ds\Delta B=\nabla\vec{n}\cdot\nabla^\perp(\vec{E}-\vec{v})\quad\mbox{ in }B_R(0)\\[5mm]
\ds B=0
\end{array}
\rg.
\ee
Similarly as above we obtain
\be
\label{VIII.62}
\|\nabla B\|_{L^{2,1}(B_R(0))}=o(1)\quad\mbox{ and }\quad\|\nabla B\|_{L^\infty(\p B_t(0))}\le C\ o(1)\ t^{-1}
\ee
Denote
\[
 \eta_4(t):=t^{-1}\ \lf[\int_{B_{2t}\setminus B_{t/2}}|\nabla B|^2\ dx^2\rg]^{1/2}=o(1)\, t^{-1}
\]
From (\ref{VIII.62}) we obtain
\be
\label{VIII.62-a}
\int_{2r}^{R/2}\eta_4^2(t)\ t\, dt=o(1)
\ee
Interpolation inequalities give again 
\be
\label{VIII.62-b}
\|\nabla B\|_{L^\infty(\p B_t(0))}\le C\ [o(1)\, \delta(t)+\eta_4(t)]   \le o(1)\ t^{-1}
\ee
Let $D$ such that
\be
\label{VIII.63}
\vec{n}\cdot\nabla^\perp(\vec{E}-\vec{v})=\nabla B+\nabla^\perp D
\ee
Similarly as above we obtain
\be
\label{VIII.64}
\begin{array}{l}
\ds\|\nabla D\|_{L^{2}(B_R(0))}=o(1)\quad\\[5mm]
\ds\mbox{ and }\quad\|\nabla D\|_{L^\infty(\p B_t(0))}\le C\ \lf[o(1)\, \delta(t)+\eta_1(t)+\eta_2(t)+\eta_4(t)\rg]=o(1)\ t^{-1}
\end{array}
\ee
Denoting $U:=S+D$ we have form (\ref{VIII.37}) and (\ref{VIII.64})
\be
\label{VIII.64-a}
\|\nabla U\|_{L^\infty(\p B_t(0))}\le C\ t^{-1}\quad\quad\Longrightarrow\quad\quad \|\nabla U\|_{L^{2,\infty}(B_{R}\setminus B_r)}\le C
\ee
The pair $(U,\vec{T})$ satisfies the following system
\be
\label{VIII.65}
\lf\{
\begin{array}{l}
\ds\nabla U=\nabla^\perp C-\vec{n}\cdot\nabla^\perp\vec{T}\\[5mm]
\ds\nabla\vec{T}=-\nabla^\perp\vec{v}+\nabla^\perp\vec{E}+\vec{n}\ \nabla^\perp S+\vec{n}\times\nabla^\perp\vec{X}
\end{array}
\rg.
\ee
Let $U_t:=|\p B_t|^{-1}\,\int_{\p B_t} U$ and $\vec{T}_t:=|\p B_t|^{-1}\,\int_{\p B_t} \vec{T}$. The system (\ref{VIII.65}) implies
\be
\label{VIII.66}
\lf\{
\begin{array}{l}
\ds \frac{d U_t}{dt}=-\frac{1}{2\pi}\int_0^{2\pi}[\vec{n}-\vec{n}_t]\cdot\frac{1}{t}\frac{\p\vec{T}}{\p \theta}\, d\theta\\[5mm]
\ds\frac{d\vec{T}_t}{dt}=\frac{1}{2\pi}\int_0^{2\pi}[\vec{n}-\vec{n}_t]\  \frac{1}{t}\frac{\p S}{\p \theta}\, d\theta+\frac{1}{2\pi}\int_0^{2\pi}[\vec{n}-\vec{n}_t]\times\frac{1}{t}\frac{\p\vec{X}}{\p \theta}\, d\theta
\end{array}
\rg.
\ee
where $\vec{n}_t:=|\p B_t|^{-1}\,\int_{\p B_t} \vec{n}$. We have
\[
\|\vec{n}-\vec{n}_t\|_{L^\infty(\p B_t)}\le C\ t\ \delta(t)
\]
Hence using the estimates above we obtain
\be
\label{VIII.67}
\lf|\frac{d U_t}{dt}\rg|+\lf|\frac{d \vec{T}_t}{dt}\rg|\le C\, \delta(t)
\ee
which implies
\be
\label{VIII.68}
\int_r^R\lf|\frac{d U_t}{dt}\rg|^2+\lf|\frac{d \vec{T}_t}{dt}\rg|^2\ t\ dt\le C
\ee
where $C$ is again independent of $k$. The system (\ref{VIII.65}) implies
\be
\label{VIII.69}
\lf\{
\begin{array}{l}
\ds\Delta U=-\nabla\vec{n}\cdot\nabla^\perp\vec{T}\\[5mm]
\ds\Delta\vec{T}=\nabla\vec{n}\cdot \nabla^\perp S+\nabla\vec{n}\times\nabla^\perp\vec{X}
\end{array}
\rg.
\ee
We can then make use of lemma 10 of \cite{LR1} we obtain that $\nabla U$ and $\nabla \vec{T}$ are uniformly bounded in $L^2$
\be
\label{VIII.70}
\int_{B_R\setminus B_r}|\nabla U|^2+|\nabla\vec{T}|^2\ dx^2\le C
\ee
which itself implies using (\ref{VIII.64}) and (\ref{VIII.58})
\be
\label{VIII.71}
\int_{B_R\setminus B_r}|\nabla S|^2+|\nabla\vec{X}|^2\ dx^2\le C
\ee
where $C$ is independent of $k$. We bootstrap this information in (\ref{VIII.66}) as follows
\be
\label{VIII.72}
\begin{array}{l}
\ds \int_r^R\lf|\frac{d U_t}{dt}\rg|+\lf|\frac{d \vec{T}_t}{dt}\rg|\ dt\le\int_r^R\ \delta(t)\ \int_{\p B_t}|\nabla\vec{T}|+|\nabla\vec{X}|+|\nabla S|\ \\[5mm]
\ds\quad\le \lf[\int_r^R\ \delta^2(t)\ t\ dt\rg]^{1/2}\  \lf[\int_r^R\  t^{-1}\lf[   \int_{\p B_t}|\nabla\vec{T}|+|\nabla\vec{X}|+|\nabla S| \rg]^2\ dt\rg]^{1/2}\\[5mm]
\ds\quad\le \lf[\int_r^R\ \delta^2(t)\ t\ dt\rg]^{1/2}\ \lf[\int_{B_R\setminus B_r}|\nabla S|^2+|\nabla\vec{T}|^2+|\nabla\vec{X}|^2\ dx^2\rg]^{1/2}\le C
\end{array}
\ee
We can chose $S$ and $\vec{R}$ which were fixed modulo the addition of an arbitrary constant in such a way that
\be
\label{VIII.73}
0=U_r= |\p B_r|^{-1}\,\int_{\p B_r} U\quad\quad\mbox{ and }\quad\quad 0=\vec{T}_r= |\p B_r|^{-1}\,\int_{\p B_r} \vec{T}\
\ee
combining this choice with (\ref{VIII.72}) we obtain
\be
\label{VIII.74}
|U_t|_{L^\infty([r,R])}+|\vec{T}_t|_{L^\infty([r,R])}\le C
\ee 
We can then make use of lemma 8 of \cite{LR1} and deduce
\be
\label{VIII.75}
\|\nabla U\|_{L^{2,1}(B_R\setminus B_r)}+\|\nabla \vec{T}\|_{L^{2,1}(B_R\setminus B_r)}\le C
\ee
We have from (\ref{VI.19-2})
\be
\label{VIII.76}
\begin{array}{l}
\ds 2\,(1+2\, \sigma^2\, (1+H^2)-l_\sigma\,\al)\ e^{2\,\la}\vec{H}=\nabla^\perp S\cdot\nabla\vec{\Phi}-\nabla\vec{R}\times\nabla^\perp\vec{\Phi}-\nabla\vec{v}\times\nabla\vec{\Phi}\\[5mm]
=\nabla^\perp U\cdot\nabla\vec{\Phi}-\nabla\vec{T}\times\nabla^\perp\vec{\Phi}-\nabla\vec{v}\times\nabla\vec{\Phi}-\nabla^\perp D\cdot\nabla\vec{\Phi}+\nabla(\vec{F}+\vec{Z})\times\nabla^\perp\vec{\Phi}
\end{array}
\ee
Hence (\ref{VIII.75}) implies
\be
\label{VIII.77}
\|h_\sigma^{-1}\ e^{\la}\vec{H}+[\nabla\vec{v}+\nabla^\perp(\vec{F}+\vec{Z})]\times\nabla\vec{\Phi}\, e^{-\la}+\nabla^\perp D\cdot\nabla\vec{\Phi}\, e^{-\la}\|_{L^{2,1}(B_R\setminus B_r)}\le C
\ee
where $h_\sigma^{-1}:=2\,(1+2\, \sigma^2\, (1+H^2)-l_\sigma\,\al)$. For any $\ep>0$ we choose $r^k$ and $R^k$ such that
\[
\|s\,\delta^k(s)\|_{L^\infty([r^k,R^k])}\le \ep
\]
which implies in particular using (\ref{VIII.20}) that
\be
\label{VIII.78}
\|(h^k_\sigma)^{-1}\ e^{\la^k}\vec{H^k}\|_{L^{2,\infty}(B_{R^k}\setminus B_{r^k})}\le C\, \ep
\ee
Using (\ref{VIII.39}),  (\ref{VIII.52}), (\ref{VIII.58}) and (\ref{VIII.64}) , for $k$ large enough we have
\be
\label{VIII.79}
\|[\nabla\vec{v}^{\,k}+\nabla^\perp(\vec{F}^{\,k}+\vec{Z}^{\,k})]\times\nabla\vec{\Phi}^{\, k}\, e^{-\la^k}+\nabla^\perp D^k\cdot\nabla\vec{\Phi}^{\,k}\, e^{-\la^k}\|_{L^{2}(B_{R^k}\setminus B_{r^k})}\le\, \ep
\ee
Combining (\ref{VIII.78}) and (\ref{VIII.79}) we obtain in particular that
\be
\label{VIII.80}
\|(h^k_\sigma)^{-1}\ e^{\la^k}\vec{H^k}+[\nabla\vec{v}^{\,k}+\nabla^\perp(\vec{F}^{\,k}+\vec{Z}^{\,k})]\times\nabla\vec{\Phi}^{\, k}\, e^{-\la^k}+\nabla^\perp D^k\cdot\nabla\vec{\Phi}^{\,k}\, e^{-\la^k}\|_{L^{2,\infty}(B_{R^k}\setminus B_{r^k})}\le\, C\,  \ep
\ee
Combining (\ref{VIII.77}) and (\ref{VIII.80}) we obtain
\be
\label{VIII.81}
\|(h^k_\sigma)^{-1}\ e^{\la^k}\vec{H^k}+[\nabla\vec{v}^{\,k}+\nabla^\perp(\vec{F}^{\,k}+\vec{Z}^{\,k})]\times\nabla\vec{\Phi}^{\, k}\, e^{-\la^k}+\nabla^\perp D^k\cdot\nabla\vec{\Phi}^{\,k}\, e^{-\la^k}\|_{L^{2}(B_{R^k}\setminus B_{r^k})}\le\, C\,  \sqrt{\ep}
\ee
Combining (\ref{VIII.79}) and (\ref{VIII.81}) we then obtain
\be
\label{VIII.82}
\|(h^k_\sigma)^{-1}\ e^{\la^k}\vec{H^k}\|_{L^{2}(B_{R^k}\setminus B_{r^k})}\le\, C\,  [\sqrt{\ep}+\ep]
\ee
which implies that the Willmore energy is as small as we want in any neck region for $k$ large enough. We deduce lemma~\ref{lm-energy-quantiz} from this fact and the final arguments of \cite{BR3}
\section{Appendix}
\begin{Lma}
\label{lm-wente}
There exists $C>0$ such that for any $\sigma\in (0,1)$ and any $a\in L^{2,\infty}(D^2)$ and $b\in L^{2,\infty}(D^2)$, denoting $\varphi$ the $W^{1,1}$ solution
to the following equation
\be
\label{A-1}
\lf\{
\begin{array}{l}
\ds-\Delta\varphi=\p_{x_1}a\,\p_{x_2}b-\p_{x_2}a\,\p_{x_1}b\quad\mbox{in }D^2\\[5mm]
\ds\varphi=0\quad\quad\mbox{ on }\p D^2
\end{array}
\rg.
\ee
then the following inequality holds
\be
\label{A-2}
\|\nabla\varphi\|_{L^{2,\infty}}\le\, C\, \|\nabla a\|_{L^{2,\infty}\cap\sigma^{1/2} L^4(D^2)}\ \|\nabla b\|_{L^{2,\infty}+\sigma^{-1/2}L^{4/3}(D^2)}
\ee
and
\be
\label{A-2-b}
\|\nabla\varphi\|_{L^{2,\infty}+\sigma^{-1/2}L^{4/3}(D^2)}\le\, C\, \|\nabla a\|_{L^{2,\infty}}\ \|\nabla b\|_{L^{2,\infty}+\sigma^{-1/2}L^{4/3}(D^2)}
\ee
where
\[
\|f\|_{L^{2,\infty}\cap\sigma^{1/2} L^4(D^2)}:=\|f\|_{L^{2,\infty}(D^2)}+\sigma^{1/2}\,\|f\|_{L^4(D^2)}\quad,
\]
 and 
 \[
 \|f\|_{L^{2,\infty}+\sigma^{-1/2} L^{4/3}(D^2)}:=\inf\lf\{\|f_1\|_{L^{2,\infty}(D^2)}+\sigma^{-1/2}\|f_2\|_{L^{4/3}(D^2)}\quad;\quad f=f_1+f_2\rg\}\quad.
\]
\end{Lma}
\noindent{\bf Proof of lemma~\ref{lm-wente}.}
Let $X_1$ and $X_2$ such that
\be
\label{A-2-2}
\nabla b= X_1+X_2\quad\|X_1\|_{L^{2,\infty}(D^2)}+\sigma^{-1/2}\|X_2\|_{L^{4/3}(D^2)}\le 2\, \|\nabla b\|_{L^{2,\infty}+\sigma^{-1/2}L^{4/3}(D^2)}
\ee
Let $c_i$ for $i=1,2$ satisfying
\be
\label{A-3}
\lf\{
\begin{array}{l}
\ds-\Delta c_i=\mbox{ div}(X_i^\perp)\quad\mbox{ in }D^2\\[5mm]
\ds c_i=0\quad\quad\quad\mbox{ on }\p D^2
\end{array}
\rg.
\ee
There exists a constant $C$ independent of $\sigma$ such that
\be
\label{A-4}
\|\nabla c_1\|_{L^{2,\infty}(D^2)}\le\, C\,  \|X_1\|_{L^{2,\infty}(D^2)}\quad\mbox{ and }\quad\|\nabla c_2\|_{L^{4/3}(D^2)}\le\, C\,  \|X_2\|_{L^{4/3}(D^2)}\quad.
\ee
Applying Poincar\'e Lemma we obtain the existence of $b_i$  such that $\int_{D^2}b_i=0$ and
\be
\label{A-5-a}
X_i^\perp+\nabla c_i=\nabla^\perp b_i\quad\Longleftrightarrow\quad X^i=\nabla b_i+\nabla^\perp c_i
\ee
and we have
\be
\label{A-5}
\|\nabla b_1\|_{L^{2,\infty}(D^2)}\le\, (C+1)\,  \|X_1\|_{L^{2,\infty}(D^2)}\quad\mbox{ and }\quad\|\nabla b_2\|_{L^{4/3}(D^2)}\le\, (C+1)\,  \|X_2\|_{L^{4/3}(D^2)}\quad.
\ee
Observe that
\be
\label{A-6}
\nabla b=\nabla b_1+\nabla b_2 +\nabla^\perp c_1+\nabla^\perp c_2\quad\quad \Delta (c_1+c_2)=0
\ee
Since $c_1+c_2=0$ on $\p D^2$ we have $c_1+c_2\equiv 0$ on $D^2$ and
\be
\label{A-7}
\nabla b=\nabla b_1+\nabla b_2 
\ee
Let
\be
\label{A-8}
\lf\{
\begin{array}{l}
\ds-\Delta\varphi_i=\p_{x_1}a\,\p_{x_2}b_i-\p_{x_2}a\,\p_{x_1}b_i\quad\mbox{in }D^2\\[5mm]
\ds\varphi_i=0\quad\quad\mbox{ on }\p D^2
\end{array}
\rg.
\ee
Using Wente estimates (see for instance \cite{Ri1}) we obtain respectively
\be
\label{A-9}
\|\nabla\varphi_1\|_{L^{2,\infty}(D^2)}\le C\ \|\nabla a\|_{L^{2,\infty}(D^2)}\, \|\nabla b_1\|_{L^{2,\infty}(D^2)}\quad,
\ee
and
\be
\label{A-10}
\|\nabla\varphi_2\|_{L^{2}(D^2)}\le C\ \|\nabla a\|_{L^{4}(D^2)}\, \|\nabla b_2\|_{L^{4/3}(D^2)}
\ee
Combining (\ref{A-9}) and (\ref{A-10}) we obtain (\ref{A-2}). We now write
\be
\label{A-11}
\lf\{
\begin{array}{l}
\ds-\Delta\varphi_2=\mbox{div}(b_2\,\nabla^\perp a)\quad\mbox{in }D^2\\[5mm]
\ds\varphi=0\quad\quad\mbox{ on }\p D^2
\end{array}
\rg.
\ee
Sobolev-Lorentz embedding theorem gives
\be
\label{A-12}
\|b_2\|_{L^{4,4/3}(D^2)}\le C\, \|\nabla b_2\|_{L^{4/3}(D^2)}
\ee
Hence we have using fundamental properties of Lorentz spaces (see \cite{Graf1})
\be
\label{A-13}
\|b_2\,\nabla^\perp a\|_{L^{4/3}(D^2)}\le C\ \|\nabla a\|_{L^{2,\infty}(D^2)}\ \|b_2\|_{L^{4,4/3}(D^2)}\le C\, \|\nabla a\|_{L^{2,\infty}(D^2)}\ \|\nabla b_2\|_{L^{4/3}(D^2)}
\ee
Hence using classical elliptic estimates we have
\be
\label{A-14}
\|\nabla\varphi_2\|_{L^{4/3}(D^2)}\le C\ \|\nabla a\|_{L^{2,\infty}(D^2)}\ \|\nabla b_2\|_{L^{4/3}(D^2)}
\ee
Combining (\ref{A-9}) and (\ref{A-14}) we obtain using (\ref{A-2-2})
\be
\label{A-15}
\begin{array}{l}
\|\nabla\varphi_1\|_{L^{2,\infty}(D^2)}+\sigma^{-1/2}\|\nabla\varphi_2\|_{L^{4/3}(D^2)}\\[5mm]
\ds\quad\quad\le C\ \|\nabla a\|_{L^{2,\infty}(D^2)}\ \lf[\|\nabla b_1\|_{L^{2,\infty}(D^2)}+\sigma^{-1/2}\|\nabla b_2\|_{L^{4/3}(D^2)}\rg]\\[5mm]
\ds\quad\quad\le C\,\|\nabla a\|_{L^{2,\infty}(D^2)}\ \lf[\|X_1\|_{L^{2,\infty}(D^2)}+\sigma^{-1/2}\|X_2\|_{L^{4/3}(D^2)}\rg]\\[5mm]
\ds\quad\quad\le C\,\|\nabla a\|_{L^{2,\infty}(D^2)}\  \|\nabla b\|_{L^{2,\infty}+\sigma^{-1/2}L^{4/3}(D^2)}
\end{array}
\ee 
This implies (\ref{A-2-b}) and lemma~\ref{lm-wente} is proved.\hfill $\Box$
\begin{Lma}
\label{lm-compcomp}
Let $m\in {\N}^\ast$, and $1\le p<+\infty$. There exists $\ep(m,p)>0$ such that for any sequence of maps ${\mathcal A}_k\in W^{1,2p}(D^2, M_m({\R}))$ satisfying
\be
\label{A-15-a}
\int_{D^2}|\nabla{\mathcal A}_k|^{2p}\le\ep(m,p)\quad
\ee
and weakly converging to ${\mathcal A}_\infty$, for any sequence of maps $\vec{\varphi}_k$ weakly converging in $W^{1,2p/2p-1}(D^2,{\R}^m)$ and any sequence of maps
$\vec{F}_k$ strongly converging to $\vec{F}_\infty$ in $L^{2p/2p-1}(D^2,{\R}^2\otimes{\R}^m)$ to a limit $\vec{F}_\infty$ and satisfying
\[
-\Delta\vec{\varphi}_k=\nabla{\mathcal A}_k\cdot\nabla^\perp\vec{\varphi}_k+\mbox{div}\ \vec{F}_k\quad\quad\mbox{ in }\quad{\mathcal D}'(D^2)\quad,
\]
then $\vec{\varphi}_k$ strongly converges in $W_{{loc}}^{1,2p/2p-1}(D^2,{\R}^m)$. If one assume further $p>1$ the smallness condition (\ref{A-15-a}) is not needed for the
same result to hold.\hfill $\Box$ 
\end{Lma}

\begin{Lma}
\label{lm-L-4} For any $C_0>0$ there exists $\ep>0$ such that for any conformal weak immersion in ${\mathcal E}_{\Sigma,2}(D^2)$ satisfying
\[
\|\nabla\la\|_{L^{2,\infty}(D^2)}\le C_0\quad\mbox{ and }\quad\int_{D^2}|\nabla\vec{n}|^2\ dx^2<\ep
\]
then
\be
\label{A-16}
\sigma^2\int_{D^2_{1/2}}|\nabla\vec{n}|^4\ e^{-2\la}\, dx^2\le C\ \sigma^2\,\int_{D^2}H^4\,e^{2\la}\ dx^2+C\, \sigma^2\ e^{-2\ov{\la}}\lf[\int_{D^2}|\nabla\vec{n}|^2\ dx^2\rg]^2
\ee
where $e^\la=|\p_{x_1}\vec{\Phi}|=|\p_{x_2}\vec{\Phi}|$ and $\ov{\la}=|D^2_{1/2}|^{-1}\,\int_{D^2_{1/2}}\la(x)\ dx^2$. \hfill $\Box$
\end{Lma}
\noindent{\bf Proof of lemma~\ref{lm-L-4}.} Arguing as in the beginning of the proof of lemma~\ref{lm-ep-reg}  we have 
\be
\label{A-17}
\|\al-\ov{\al}\|_{L^\infty(D^2_{5/6})}\le C
\ee
We have also
\be
\label{VI.25}
\Delta\vec{\Phi}=2\, e^{2\la} \ \vec{H}=2\, e^{2\ov{\la}}\, e^{2\, (\la-\ov{\la})}\ \vec{H}\quad\mbox{ in }D^2
\ee
This gives that $\nabla^2\vec{\Phi}\in L^4(D^2_{3/4})$ and hence $\nabla\vec{n}\in L^4(D^2_{3/4})$
Let $\vec{a}$ satisfying
\be
\label{VI.26}
\lf\{
\begin{array}{l}
\ds\Delta\vec{a}=\mbox{div}(\vec{n}\times\nabla\vec{n})\quad\mbox{ in }D^2\\[5mm]
\ds\vec{a}=0\quad\quad\mbox{ on }\p D^2
\end{array}
\rg.
\ee
and let $\vec{b}$ such that $\vec{n}\times\nabla\vec{n}=\nabla\vec{a}+\nabla^\perp\vec{b}$. Using classical elliptic estimates we have
\be
\label{VI.27}
\int_{D^2}|\nabla\vec{a}|^2+|\nabla\vec{b}|^2\ dx^2\le C\ \int_{D^2}|\nabla\vec{n}|^2\ dx^2\le C\, \ep
\ee
Let $\rho\in [1/2,3/4]$ such that
\be
\label{VI.28 }
\int_{\p D^2_{\rho}}|\nabla\vec{a}|^2+|\nabla\vec{b}|^2\ dl\le 4 \int_{D^2}|\nabla\vec{n}|^2\ dx^2\le 4\ep
\ee
Observe that $W^{1,2}(\p D^2_\rho) \hookrightarrow     W^{1-1/4,4}(\p D^2_\rho)$. Hence we have
\be
\label{VI.29}
\|\vec{a}\|_{W^{1-1/4,4}(\p D^2_\rho)}+\|\vec{b}\|_{W^{1-1/4,4}(\p D^2_\rho)}\le C\, \lf[  \int_{D^2}|\nabla\vec{n}|^2\ dx^2\rg]^{1/2}
\ee
Recall now from \cite{Ri1} the following general formula
\be
\label{VI.30}
\vec{n}\times\nabla\vec{n}=2\, H\, \nabla^\perp\vec{\Phi}+\nabla^\perp\vec{n}
\ee
Hence we have in particular on $D^2_\rho$
\be
\label{VI.31}
\Delta\vec{a}= 2\,\mbox{div}(H\, \nabla^\perp\vec{\Phi})
\ee
Classical elliptic estimates give then, using (\ref{VI.29}),
\be
\label{VI.32}
\|\nabla\vec{a}\|_{L^4(D^2_\rho)}\le C\, e^{\ov{\la}}\lf[\int_{D^2_\rho}H^4\, dx^2\rg]^{1/4}+C\, \lf[  \int_{D^2}|\nabla\vec{n}|^2\ dx^2\rg]^{1/2}
\ee
We have also on $D^2_\rho$
\be
\label{VI.33}
\Delta\vec{b}=\nabla^\perp\vec{n}\times\nabla\vec{n}
\ee
\be
\label{VI.34}
\begin{array}{l}
\ds\|\nabla\vec{b}\|_{L^4(D^2_\rho)}\le C\, \lf[\int_{D^2_\rho}|\nabla\vec{n}|^{4/3}\, |\nabla\vec{n}|^{4/3}\ dx^2\rg]^{3/4}\\[5mm]
\ds\quad\quad\quad\quad\quad+C\, \lf[  \int_{D^2}|\nabla\vec{n}|^2\ dx^2\rg]^{1/2}
\end{array}
\ee
Hence
\be
\label{VI.35}
\|\nabla\vec{b}\|_{L^4(D^2_\rho)}\le C\, \lf[  \int_{D^2_{\rho}}|\nabla\vec{n}|^2\ dx^2\rg]^{1/2}\ \|\nabla\vec{n}\|_{L^4(D^2_\rho)}+C\, \|\nabla\vec{n}\|_{L^2(D^2)}
\ee
Combining (\ref{VI.32}) with (\ref{VI.35}) together with the fact that $\|\nabla\vec{n}\|_{L^2(D^2)}<\ep$ gives
\be
\label{VI.36}
\begin{array}{l}
\ds\|\nabla\vec{n}\|_{L^4(D^2_\rho)}\le \, C\, e^{\,\ov{\la}}\lf[\int_{D^2_\rho}H^4\, dx^2\rg]^{1/4}+C\, \ep\ \|\nabla\vec{n}\|_{L^4(D^2_\rho)}\\[5mm]
\ds\quad\quad\quad +C\, \|\nabla\vec{n}\|_{L^2(D^2)}
\end{array}
\ee
Hence for $\ep$ small enough we finally obtain
\be
\label{VI.37}
\sigma^2\int_{D^2_{\rho}}|\nabla\vec{n}|^4\ e^{-2\la}\, dx^2\le C\ \sigma^2\,\int_{D^2_\rho}H^4\,e^{2\la}\ dx^2+\sigma^2\ e^{-2\ov{\la}}\lf[\int_{D^2}|\nabla\vec{n}|^2\ dx^2\rg]^2
\ee
\begin{Lma}
\label{lm-L-5}
There exists $\ep>0$ such that for any $\vec{\phi}$ in $W^{1,2}(D^2)$ satisfying
\be
\label{A-18}
\Delta\vec{\phi}=\nabla^\perp\vec{n}\times\nabla\vec{\phi}+\mbox{div} \,\vec{F}
\ee
where
\be
\label{A-19}
\int_{D^2}|\nabla\vec{n}|^2\ dx^2\le\ep\quad\mbox{ and }\quad\sup_{B_r(x)\subset D^2}r^{-\gamma}\,\int_{B_r(x)}|\vec{F}|^2\ dx^2<+\infty
\ee
for some $\gamma\in (0,1)$ then
\be
\label{A-20}
\sup_{B_r(x)\subset D^2_{1/2}}r^{-\gamma}\,\int_{B_r(x)}|\nabla\vec{\phi}|^2\ dx^2\le C_\gamma\ \lf[\sup_{B_r(x)\subset D^2}r^{-\gamma}\,\int_{B_r(x)}|\vec{F}|^2\ dx^2+\,\int_{D^2}|\nabla\vec{\phi}|^2\ dx^2 \rg]
\ee
where $C_\gamma$ depends only on $\gamma\in (0,1)$.
\hfill $\Box$
\end{Lma}
\renewcommand{\theequation}{A.\arabic{equation}}
\renewcommand{\theTh}{A.\arabic{Th}}
\renewcommand{\theProp}{A.\arabic{Prop}}
\renewcommand{\theLma}{A.\arabic{Lma}}
\renewcommand{\theCo}{A.\arabic{Co}}
\renewcommand{\theRm}{A.\arabic{Rm}}
\renewcommand{\theequation}{A.\arabic{equation}}
\setcounter{equation}{0} 
\reset
\appendix

\end{document}